\documentclass[10pt, reqno]{amsart}
\usepackage{amsthm, amsmath, amsfonts, amssymb, enumerate, textcmds, tensor, enumitem, mathdots}
\usepackage{geometry}
\usepackage{pst-node}
\usepackage{tikz-cd} 
\usepackage{appendix}
\usepackage{graphicx,tikz}
\usepackage{nicefrac}
\usepackage{faktor}
\usepackage[colorlinks=true]{hyperref}
\usepackage{mathtools}
\usepackage[style=alphabetic,maxbibnames=99,maxalphanames=99]{biblatex}
\usepackage{xcolor}
\usepackage{comment}

\hypersetup{
	colorlinks,
	citecolor=black,
	filecolor=black,
	linkcolor=black,
	urlcolor=black
}
\usepackage{cleveref}

\usepackage{physics}
\usepackage{enumerate}
\usepackage{pinlabel}
\usepackage[colorinlistoftodos]{todonotes}
\newcommand{\R}{\mathbb R}

\newcommand{\Icpx}{\mathbb{I}}
\newcommand{\Jcpx}{\mathbb{J}}
\newcommand{\Kcpx}{\mathbb{K}}

\newcommand{\gpara}{\mathbf{g}}

\newcommand{\Ipara}{\mathbf{I}}
\newcommand{\Jpara}{\mathbf{J}}
\newcommand{\Kpara}{\mathbf{K}}

\newcommand{\Ads}{\mathbb{A}\mathrm{d}\mathbb{S}}
\newcommand{\PSL}{\mathrm{PSL}}

\newcommand{\K}{\mathrm{K}}

\newtheorem{thm}{Theorem}[section]
\newtheorem{lem}[thm]{Lemma}

\newtheorem{prop}[thm]{Proposition}
\newtheorem{cor}[thm]{Corollary}

\newtheorem*{thma'}{Theorem A'}

\newtheorem*{thmb'}{Theorem B'}

\DeclareMathOperator{\id}{id}

\DeclareMathOperator{\C}{\mathbb{C}}
\DeclareMathOperator{\QD}{\mathrm{QD}}

\DeclareMathOperator{\Real}{Re}
\DeclareMathOperator{\Ima}{Im}

\DeclareMathOperator{\Diff}{\mathrm{Diff}}
\DeclareMathOperator{\CP}{\mathbb{C}\mathbb{P}}

\DeclareMathOperator{\SL}{\mathrm{SL}}

\DeclareMathOperator{\Hess}{\mathrm{Hess}}

\DeclareMathOperator{\inners}{\langle \cdot, \cdot \rangle}
\DeclareMathOperator{\Hyp}{\mathbb H}

\newcommand{\ccpair}{(J_1, -J_2)}

\newcommand{\hpair}{\vb*h{(J_1, -J_2)}}
\DeclareMathOperator{\phee}{\varphi}
\DeclareMathOperator{\eps}{\varepsilon}

\newcommand{\redd}[1]{\textcolor{red}{#1}}

\newcommand{\putref}{\textcolor{red}{(PUT REF)}}
\newcommand{\inner}[1]{\langle{#1}\rangle}

\usepackage{physics}
\usepackage{soul}

\usepackage{scalerel}

\usepackage{stmaryrd}
\usepackage{upgreek}

\newcommand{\CSCS}{\mathcal C(S)\times \mathcal C(\overline S)}
\newcommand{\TSTS}{\mathcal T(S)\times \mathcal T(\overline S)}

\newcommand{\AT}[1]{{\textcolor{teal}{  AT: #1}}}
\newcommand{\NR}[1]{{\textcolor{blue}{  NR: #1}}}

\usepackage{mathrsfs}
\usepackage{soul}

\newcommand{\Is}{\mathtt I}
\newcommand{\II}{\mathtt {II}}
\newcommand{\QF}{\mathcal{QF}}
\newcommand{\Hom}{\mathrm{Hom}}

\theoremstyle{definition}

\newtheorem{defn}[thm]{Definition}
\newtheorem{remark}[thm]{Remark}

\makeatletter
\newcommand\longleftrightarrowfill@{%
  \arrowfill@\leftarrow\relbar\rightarrow}
\makeatother

\author[C. El Emam]{Christian El Emam}
\address[Christian El Emam]{Department of Mathematics ``Giuseppe Peano'', University of Turin, Via Carlo Alberto 10, 10123 Torino, Italy.}
\email{christian.elemam@unito.it}

\author[F. Mazzoli]{Filippo Mazzoli}
\address[Filippo Mazzoli]{Department of Mathematics, University of California Riverside, United States}
\email{filippo.mazzoli@ucr.edu}

\author[N. Rungi]{Nicholas Rungi}
\address[Nicholas Rungi]{Instituto de Ciencias Matemáticas (CSIC-UAM-UC3M-UCM), Nicolás Cabrera 13--15, Cantoblanco, 28049 Madrid, Spain.}
\email{nicholas.rungi@unito.it}

\author[A. Seppi]{Andrea Seppi}
\address[Andrea Seppi]{Department of Mathematics ``Giuseppe Peano'', University of Turin, Via Carlo Alberto 10, 10123 Torino, Italy.}
\email{andrea.seppi@unito.it}

\author[A. Tamburelli]{Andrea Tamburelli}
\address[Andrea Tamburelli]{Department of Mathematics, University of Pisa, Italy}
\email{andrea.tamburelli@unipi.it}

\begin{document}
\title[Complex hyperk\"ahler metric]{Extending Donaldson's HyperKähler Structure via Minimal Surfaces into $\PSL(2,\C)$ 
}

\begin{abstract}

   We construct a deformation space for minimal immersions into $\PSL(2,\C)$ with its holomorphic Killing form, and define an explicit, mapping class group invariant, complex hyperKähler structure on such deformation space. Such complex hyperKähler structure restricts to: Donaldson's hyperKähler structure on a totally real copy of almost-Fuchsian space, Mazzoli-Seppi-Tamburelli's para-hyperKähler structure on a totally real copy of maximally globally hyperbolic $\Ads^3$ manifolds, and Loustau-Sanders' complex bi-Kähler structure on a holomorphic copy of quasi-Fuchsian space. Moreover, all these loci are totally geodesic.
   
   We  discuss other geometric aspects of this space, such as the Hamiltonian property of a natural $\mathbb C^*$-action, involving the Goldman symplectic form and a complex-valued area functional.
\end{abstract}

\maketitle

\setcounter{tocdepth}{1}        
\setcounter{secnumdepth}{3}     
\hypersetup{bookmarksdepth=3}   
\tableofcontents

\section{Introduction}

Given a closed oriented surface $S$ of genus $\geq 2$, its Teichm\"uller space $\mathcal T(S)$ can be equivalently defined as the space of complex structures on $S$ compatible with the orientation of $S$, or the space of hyperbolic metrics on $S$, in both cases up to the action of the group $\mathrm{Diff}_0(S)$ of diffeomorphsims of $S$ isotopic to the identity. While the complex structure of $\mathcal T(S)$ had been introduced by Teichm\"uller already in the 1930s, in 1958 Weil introduced in \cite{weil58} a natural Hermitian metric on $\mathcal T(S)$, now known as the Weil-Petersson metric, and suggested that this would be a K\"ahler structure, which was soon after proved by Ahlfors in \cite{ahlfors61}. 

Since $\mathcal T(S)$ is topologically a ball of dimension $3|\chi(S)|$, its topology gives no obstruction to the existence of a K\"ahler metric; it is the naturality of such metric, which has to be intended as the invariance under the action of the mapping class group $\mathrm{MCG}(S):=\mathrm{Diff}_+(S)/\mathrm{Diff}_0(S)$, that makes the problem of existence of a Kähler metric non-trivial. Likewise the mapping class group invariance will be a crucial property satisfied by all the constructions that will be mentioned and tackled in this work.

Before moving on, we recall that the Teichm\"uller space has yet another interpretation as a connected component of the character variety for $G=\PSL(2,\R)$. Recall that the character variety of a Lie group $G$ is defined as 
$$\chi(\pi_1(S),G):=\mathrm{Hom}\big(\pi_1(S),G\big)/\!\!/G~,$$
where the action of $G$ on $\mathrm{Hom}\big(\pi_1(S),G\big)$ is by conjugation and the ``double'' quotient denotes the space of closed orbits, which is a Hausdorff topological space with the induced topology (see \cite{sikora2012}). In the 1980s, Goldman proved two fundamental results. First, in his PhD thesis \cite{goldmanphd}, he showed that the holonomy map is a diffeomorphism between $\mathcal T(S)$ and a (smooth) connected component of $\chi(\pi_1(S),\PSL(2,\R))$. Second, in \cite{goldman1984}, he proved that the imaginary part of the Weil-Petersson Hermitian product corresponds under such diffeomorphism to a symplectic form that can be constructed, following Atiyah-Bott \cite{atiyah-bott}, on $\chi(\pi_1(S),G)$ for very general Lie groups $G$ (here applied to the particular case of $G=\PSL(2,\R)$), now also called the Goldman form.

\subsection{Known extensions of the Weil-Petersson metric} 
A longstanding line of research consists in studying extensions of the Teichm\"uller space --- that is, deformation spaces of ``geometric structures'' (either on $S$ itself, or on other manifolds related to $S$) that generalize hyperbolic structures on $S$ --- and, in parallel,  constructing ``special metrics'' on such larger deformation spaces. Without attempting to be exhaustive, we mention here  three instances that are particularly relevant for our results.

\begin{enumerate}[leftmargin=2em]
    \item The deformation space of \emph{quasi-Fuchsian manifolds}, that is, the convex cocompact hyperbolic three-manifolds diffeomorphic to $S\times\R$, is denoted by $\QF(S)$, and is diffeomorphic to $\mathcal T(S)\times\mathcal T(\overline S)$ by the Bers' Double Uniformization theorem. Bridgemann-Taylor introduced a semi-definite metric on $\QF(S)$ in \cite{bridgeman-taylor}, whose degeneracy was studied in \cite{bridgeman2010}. A natural K\"ahler metric on $\QF(S)$ was constructed in \cite{kim-wan-zhang25}. However, the most important result for this work is the construction of a \emph{complex bi-K\"ahler metric} by Loustau-Sanders \cite{loustau-sanders}. The holomorphic Riemannian metric involved in their construction has also been studied from a different perspective in \cite{ElE}.
    \item Among quasi-Fuchsian manifolds, an important (open) subclass is that of \emph{almost-Fuchsian manifolds}, inspired by the pioneering work \cite{uhlenbeck1983} of Uhlenbeck of 1983, whose definition will be recalled in Section~\ref{subsec:intro analytic cont} of this introduction. In \cite{donaldson2003}, Donaldson introduced a \emph{hyperK\"ahler metric} on the deformation space of almost-Fuchsian manifolds, that we denote by $\mathcal{AF}(S)\subset\QF(S)$. Such hyperK\"ahler metric was studied  in more depth in \cite{hodge2005hyperkahler, trautwein2018infinite}.
    \item \emph{Maximally globally hyperbolic Anti-de Sitter manifolds} are Lorentzian counterparts of quasi-Fuchsian manifolds, and largely studied in relation with Teichm\"uller theory after the seminal work of Mess from 1990 \cite{mess07}, which also showed that their deformation space $\mathcal{MGH}(S)$ is naturally identified with $\mathcal T(S)\times\mathcal T(S)$. The second, fourth and fifth author constructed in \cite{MazSepTamb} a \emph{para-hyperKähler metric}\footnote{The notion of para-hyperKähler metric was also known in the literature under the names of hypersymplectic quotient or neutral hyperKähler structure, see \cite{hitchinhyper,fino04}} on $\mathcal{MGH}(S)$. We mention incidentally that this is an instance of special metrics on higher Teichmüller spaces, see \cite{li16,labourie,kim-zhang,RT25,ESrank2} for others.
\end{enumerate}

The three constructions above have in common that they extend, in a suitable sense, the Weil-Petersson metric, when one includes in a natural way $\mathcal T(S)$ as a subset of $\QF(S)$, $\mathcal{AF}(S)$ and $\mathcal{MGH}(S)$ respectively. In fact, at the level of character varieties, $\QF(S)$ and $\mathcal{AF}(S)$ are open subsets of $\chi(\pi_1(S),\PSL(2,\C))$, that contain $\mathcal T(S)$ as the real locus, while $\mathcal{MGH}(S)$ is a connected component of $\chi(\pi_1(S),\PSL(2,\R)\times \PSL(2,\R))$, containing $\mathcal T(S)$ as the diagonal.

\subsection{Main result} Let us focus on the main result of this article. Broadly speaking, it shows that the three aforementioned constructions, namely the complex bi-K\"ahler metric on $\QF(S)$, the hyperK\"ahler metric on $\mathcal{AF}(S)$, and the para-hyperKähler metric on $\mathcal{MGH}(S)$, have a common holomorphic extension on a larger space, having twice the dimension of each of them. 

It turns out that the appropriate notion is that of a \emph{complex hyperK\"ahler metric}, that has been introduced, from a quite different perspective than ours, in \cite{Bridgeland2021}. We will provide more details in Section~\ref{sec:intro alg} below. The definition of the appropriate space, denoted  $\mathbb{M}_0(S)$ (the subscript serves to denote that it is a connected component of a certain bigger space), is more technical and crucially involves the notion of minimal surface in the complex Lie group $\PSL(2,\C)$, on which we will elaborate in Section~\ref{sec:intro geom} below. 

For the moment, it suffices to know that this space $\mathbb{M}_0(S)$ is a complex manifold of complex dimension $6|\chi(S)|$, locally biholomorphic to $\chi(\pi_1(S),\PSL(2,\C)\times \PSL(2,\C))$, and it contains a biholomorphic copy of an open subset of  $\chi(\pi_1(S),\PSL(2,\C)\times \PSL(2,\C))$ which in turn contains:
\begin{enumerate}[leftmargin=2em]
    \item A copy of $\QF(S)$, embedded as a \emph{complex} half-dimensional submanifold, that corresponds at the level of character varieties to the diagonal embedding  \[\Delta:\PSL(2,\C)\to\PSL(2,\C)\times \PSL(2,\C)~.\]
    \item A copy of $\mathcal{AF}(S)$, embedded as a \emph{totally real} half-dimensional submanifold, now corresponding to the ``twisted diagonal'' embedding, that we denote \[\Delta^\dagger:\PSL(2,\C)\to\PSL(2,\C)\times \PSL(2,\C)~,\] sending $A$ to $(A,A^\dagger)$, where $A^\dagger$ denotes the inverse conjugate transpose.
    \item A copy of $\mathcal{MGH}(S)$, again a \emph{totally real} half-dimensional submanifold, that corresponds to the obvious inclusion \[\iota_\R:\PSL(2,\R)\times \PSL(2,\R)\to \PSL(2,\C)\times \PSL(2,\C)~.\]
\end{enumerate}

These embeddings, at the level of Lie groups (left) and deformation spaces (right) are summarized in the diagrams below.\vspace{0.5em}\newline
\begin{tikzcd}
                                                            & {\mathrm{PSL}(2,\mathbb{R})} \arrow[d, "\iota_\mathbb R"'] \arrow[ld, "\iota_\mathbb R"'] \arrow[rd, "\Delta"] &                                                                        &                                          & \mathcal{T}(S) \arrow[d] \arrow[rd] \arrow[ld] &                                          \\
{\mathrm{PSL}(2,\mathbb{C})} \arrow[rd, "\Delta"'] & {\mathrm{PSL}(2,\mathbb{C})} \arrow[d, "\Delta^\dagger"']                                                            & {\mathrm{PSL}(2,\mathbb{R})^2} \arrow[ld, "\iota_\mathbb R"] & \mathcal{QF}(S) \arrow[rd, "\mathbb C"'] & \mathcal{AF}(S) \arrow[d, "\mathbb R"']        & \mathcal{MGH}(S) \arrow[ld, "\mathbb R"] \\
                                                            & {\mathrm{PSL}(2,\mathbb{C})^2}                                                                                 &                                                                        &                                          & \mathbb M_0(S)                                 &                                         
\end{tikzcd}

Let us now state the main result of this article.
\begin{thm}\label{thm:main}
    There exists a mapping class group invariant complex hyperK\"ahler structure $(g,\mathbf i, \mathbb I, \mathbb J, \mathbb K)$ on $\mathbb{M}_0(S)$, where $\mathbf i$ is the almost-complex structure of $\mathbb{M}_0(S)$, such that:
    \begin{enumerate}[leftmargin=2em]
    \item The tensors $(-(i/2)g,\mathbf i, \mathbb I, \mathbf i\mathbb I)$ restrict to the Loustau-Sanders complex bi-K\"ahler structure on the complex embedding of $\QF(S)$.
        \item  The tensors $(g,\mathbb I, \mathbb J, \mathbb K)$ restrict to Donaldson's hyperK\"ahler structure on the totally real embedding of $\mathcal{AF}(S)$.
        \item  The tensors $(g,\mathbb I, \mathbf i\mathbb J, \mathbf i\mathbb K)$ restrict to the para-hyperK\"ahler structure of Mazzoli-Seppi-Tamburelli on the totally real embedding of $\mathcal{MGH}(S)$. 
    \end{enumerate}
    Moreover, $\QF(S)$, $\mathcal{AF}(S)$ and $\mathcal{MGH}(S)$ are all totally geodesic in $\mathbb{M}_0(S)$.
\end{thm}

The Teichmüller space $\mathcal{T}(S)$ is, as said above, a common subspace of $\QF(S)$, $\mathcal{AF}(S)$ and $\mathcal{MGH}(S)$, and $g$ restricts to twice the Weil-Petersson metric on $\mathcal{T}(S)$.
The factor $-(i/2)$ that multiplies the metric $g$ in item $(1)$ is due to a different normalization  on $\mathcal{T}(S)$ adopted by \cite{loustau-sanders}. We also operated several sign changes in the tensors $(\mathbb I, \mathbb J, \mathbb K)$ with respect to the definitions in \cite{hodge2005hyperkahler,trautwein2018infinite,MazSepTamb} (Remarks \ref{rem: convention on I} and \ref{rmk: convention I'}), in order to keep a consistent notation.

In the next two sections, Section \ref{sec:intro alg} and Section \ref{sec:intro geom}, we explain more carefully the notions involved in the statement of Theorem~\ref{thm:main}. Then, in Sections \ref{subsec:intro tensors} and \ref{subsec:intro analytic cont} we describe additional information and properties on the space $\mathbb M_0(S)$ and its complex hyperKähler structure, introducing along the way some of the techniques of the proof of Theorem \ref{thm:main}. The introduction will be concluded in Section \ref{subsec:intro area} by related result, Theorem \ref{thm:area intro}, concerning a local $\mathbb C^*$-action, that we prove to be Hamiltonian with respect to a complex-valued area functional.

\subsection{Complex hyperK\"ahler and biquaternion (sub)algebras} \label{sec:intro alg}

A complex hyperK\"ahler structure on a complex manifold $\mathbb X$ consists of the data $(g,\mathbf i, \mathbb I, \mathbb J, \mathbb K)$ as follows. First, we denote by $\mathbf i$ the (integrable) almost-complex structure of $\mathbb X$. Second, $g$ is a holomorphic Riemannian metric, that is, a complex-valued non-degenerate bilinear form on $T^{1,0}\mathbb X$, whose coefficients in complex charts are holomorphic. Finally, $\mathbb I$, $\mathbb J$, $\mathbb K$ are $\C$--linear endomorphisms of $T^{1,0}\mathbb X$ (equivalently, one could think of them as endomorphisms of the real tangent space of $\mathbb X$ that commute with $\mathbf i$) that are $g$--isometries, are $g$--parallel and satisfy the usual multiplication rules of the imaginary quaternion units, that is: $\mathbb I^2= \mathbb J^2=\mathbb K^2=\mathbb{IJK}=-\mathrm{Id}$. A complex hyperK\"ahler structure can only live on manifolds of complex dimension $4n$.

It is instructive to think at the complex hyperK\"ahler structure and Theorem~\ref{thm:main} in terms of the biquaternions (also known as complex quaternions). Recall that the algebra of  biquaternions $\mathbb H\otimes_\R \C$ has the basis over $\R$ given by $1,\sqrt{-1},i,\sqrt{-1}i,j,\sqrt{-1}j,k,\sqrt{-1}k$, where $i,j,k$ satisfy the usual quaternion identities $i^2=j^2=k^2=ijk=-1$, while they all commute with $\sqrt{-1}$. So, $\mathbf i, \mathbb I, \mathbb J, \mathbb K$ can be thought  as endomorphisms of $T^{1,0}\mathbb X$ ($\mathbf i$ being just multiplication by $\sqrt{-1}$) or equivalently of $T_\R\mathbb X$ (in which case $\mathbf i$ is the almost-complex structure), and they satisfy the same algebraic identities of the biquaternions, playing the roles of $\sqrt{-1},i,j,k$ respectively.

From this perspective, the complex bi-K\"ahler, hyperK\"ahler and para-hyperKähler structures as in the three items of Theorem~\ref{thm:main} parallel the three (up to conjugation) real 4-dimensional  semi-simple subalgebras of $\mathbb H\otimes_\R \C$:
\begin{enumerate}[leftmargin=2em]
\item The commutative algebra $\C\oplus\C$, generated by $1,\sqrt{-1},i,\sqrt{-1}i$, which corresponds to the algebraic rules of a complex bi-Kähler structure $(g,\mathbf i,\mathbb I,\mathbb P)$: one has $\mathbf i^2=\mathbb I^2=-\mathrm{Id}$, $\mathbb P=\mathbf i\mathbb I=\mathbb I\mathbf i$, hence $\mathbb P^2=\mathrm{Id}$ (i.e. $\mathbb P$ is a para-complex structure).
\item The algebra $\mathbb H$ of ordinary quaternions, generated by $1,i,j,k$, that corresponds to the rules of a hyperKähler structure $(g,\mathbb I, \mathbb J, \mathbb K)$: $\mathbb I^2= \mathbb J^2=\mathbb K^2=\mathbb{IJK}=-\mathrm{Id}$.
\item The algebra $\mathrm{M}_2(\R)$ of split quaternions, generated by $1,i,\sqrt{-1}j,\sqrt{-1}k$, thus corresponding to a para-hyperKähler structure $(g,\mathbb I, \mathbb P, \mathbb Q)$, where $\mathbb I$ is an almost-complex structure and $\mathbb P=\mathbf i\mathbb J,\mathbb Q=\mathbf i\mathbb K$ are para-complex structures.
\end{enumerate}

All these subalgebras in turn contain a copy of $\C$, which thus corresponds to the Kähler structure of $\mathcal T(S)$. In fact, $\mathcal T(S)$ is simultaneously embedded in $\QF(S)$, $\mathcal{AF}(S)$ and $\mathcal{MGH}(S)$.

To conclude this discussion, the statements from Theorem~\ref{thm:main} that certain tensors restrict to the complex bi-Kähler/hyperKähler/para-hyperKähler structures have to be understood as follows. For the metric $g$, this means that the restriction of $g$ to $\QF(S)$/$\mathcal{AF}(S)$/$\mathcal{MGH}(S)$ equals precisely the metric involved in the corresponding structure: so, it is again a holomorphic Riemannian metric on $\QF(S)$, while in the other two cases it is a \emph{real-valued} metric: Riemannian for $\mathcal{AF}(S)$ and of neutral signature for $\mathcal{MGH}(S)$. For the endomorphisms, when we say that they restrict to certain data of the the complex bi-Kähler/hyperKähler/para-hyperKähler structures, we mean that they preserve the tangent space of $\QF(S)$/$\mathcal{AF}(S)$/$\mathcal{MGH}(S)$ (for the latter two cases, this makes sense after canonically identifying $T_\R\mathbb{M}_0(S)$ with $T^{1,0}\mathbb{M}_0(S)$), and when restricted to such tangent spaces, they are equal to the listed tensors that are part of the complex bi-Kähler/hyperKähler/para-hyperKähler structures.
  
\subsection{Minimal immersions into \texorpdfstring{$\mathrm{PSL}(2,\mathbb C)$}{PSL(2,C)}} \label{sec:intro geom}

Let us now focus on the description of the space $\mathbb{M}_0(S)$ on which the complex hyperKähler metric of Theorem~\ref{thm:main} is defined. 

The constructions of the hyperKähler metric on $\mathcal{AF}(S)$ and the para-hyperKähler metric on $\mathcal{MGH}(S)$ have a central feature in common: they rely on the identification of $\mathcal{AF}(S)$ (resp. $\mathcal{MGH}(S)$) with a space of equivariant minimal (resp. maximal) immersions of $\widetilde S$ into the hyperbolic space $\mathbb H^3$ (resp. the Anti-de Sitter space $\Ads^3$). Here both minimal (in $\mathbb H^3$) and maximal (in $\Ads^3$) mean that the trace of the second fundamental form vanishes. By the work of Uhlenbeck \cite{uhlenbeck1983} (resp. Krasnov-Schlenker \cite{kirill07}), $\mathcal{AF}(S)$ (resp. $\mathcal{MGH}(S)$) is identified with an open subset of $T^*\mathcal T(S)$, the cotangent bundle of the Teichmüller space --- which we think here as the bundle over $\mathcal T(S)$ with fiber the space of holomorphic quadratic differentials. (In the case of $\mathcal{MGH}(S)$, this open subset is in fact the entire $T^*\mathcal T(S)$.) Such identification comes from the fundamental observation that the second fundamental form of a minimal/maximal surface is the real part of a  quadratic differential on the surface, which is holomorphic with respect to the complex structure induced by the first fundamental form. Since minimal/maximal surfaces in almost-Fuchsian/maximal globally hyperbolic manifolds are unique, taking the quotient by the group $\Diff_0(S)$, one can define in this way a map from $\mathcal{AF}(S)$ (resp. $\mathcal{MGH}(S)$) to $T^*\mathcal T(S)$ which is a diffeomorphism onto its image (resp. a global diffeomorphism). The complex structures $\mathbb I$ in the hyperKähler structure of $\mathcal{AF}(S)$ and the para-hyperKähler structure of $\mathcal{MGH}(S)$ are then simply the pull-backs of the natural complex structure of $T^*\mathcal T(S)$.

Inspired by these considerations, and using notions introduced in \cite{BEE} of extrinsic geometry of surfaces into the complex Lie group $\PSL(2,\C)$  endowed with a holomorphic Riemannian metric induced by the complex Killing form, we consider a space of equivariant minimal (i.e. of vanishing complex mean curvature) immersions of $\widetilde S$ into $\PSL(2,\C)$. 

In this setting, to an equivariant minimal immersion in $\PSL(2,\C)$, satisfying some additional technical hypothesis, one can associate a 4-tuple $(J_1,q_1, -J_2, \overline{q_2})$, where $q_1$ and $\overline{q_2}$ are \emph{two} quadratic differentials on $S$, holomorphic with respect to \emph{two} complex structures $J_1$ and $-J_2$, one compatible with the fixed orientation on $S$, and the other compatible with the opposite orientation. 
There are, however, some important caveats. 
\begin{itemize}[leftmargin=2em]
    \item First, a 4-tuple $(J_1,q_1, -J_2, \overline{q_2})$ as above does not determine uniquely a minimal immersion in $\PSL(2,\C)$ (up to ambient isometries). In a completely analogous manners as for minimal immersions in $\mathbb H^3$, the immersion is uniquely determined once an additional conformal factor $e^{2u}$, which  is now complex-valued and solves a complex-valued Gauss equation (see Equation~\eqref{eq: Gauss PSL(2,C)}, is given. For this reason, we shall instead consider 5-tuples of the form $(J_1,q_1, -J_2, \overline{q_2}, u)$, for $u:S\to\C$ a smooth function.
    \item Second, as mentioned above, the datum of the 4-tuple $(J_1,q_1, -J_2, \overline{q_2})$ does not determine $u$, and $u$ might not be even locally determined, which generates technical difficulties. For this reason, we consider the space, that we denote $\mathrm M^*(S)$, of 5-tuples $(J_1,q_1, -J_2, \overline{q_2}, u)$ such that the linearization of the Gauss equation with respect to $u$ is an invertible operator.
    \item Last, and perhaps most important, since the group $\Diff_0(S)$ acts on pairs $(J,q)$ by pull-back, there is an action of $\Diff_0(S)\times \Diff_0(S)$ on the space of 4-tuples $(J_1,q_1, -J_2, \overline{q_2})$ whose quotient is $T^*\mathcal T(S)\times T^*\mathcal T(\overline{S})$, but this action is \emph{not} well-defined on $\mathrm M^*(S)$. On $\mathrm M^*(S)$, there is only a natural $\Diff_0(S)$--action, acting diagonally on $(J_1,q_1, -J_2, \overline{q_2})$ and by pre-composition on $u$. As a consequence, there are \emph{infinitely many} equivariant minimal immersions that induce the same element in $T^*\mathcal T(S)\times T^*\mathcal T(\overline{S})$, and also \emph{infinitely many} that induce the same holonomy representation from $\pi_1(S)$ in $\PSL(2,\C)\times \PSL(2,\C)$, the isometry group of $\PSL(2,\C)$.
\end{itemize}

The fundamental step for our constriction is the definition, and the study of the properties, of a ``good'' space 
$\mathbb{M}_0(S)$ where the complex hyperKähler structure lives, summarized as follows. %

\begin{thm}\label{thm: local biholos intro}
    There exists a quotient space $\mathbb{M}(S)$ of $\mathrm M^*(S)$, which is a complex manifold of complex dimension $6|\chi(S)|$, such that the following two maps :
    \begin{itemize}[leftmargin=2em]
        \item the map $\mathrm M^*(S)\to T^*\mathcal T(S)\times T^*\mathcal T(\overline{S})$ associating to the 5-tuple
        $(J_1,q_1, -J_2, \overline{q_2}, u)$
        the pairs $[J_1,q_1]\in T^*\mathcal T(S)$ and $[-J_2,\overline{q_2}]\in T^*\mathcal T(\overline{S})$, and
        \item the holonomy map $\mathrm M^*(S)\to\chi\big(\pi_1(S), \PSL(2,\C)\times \PSL(2,\C)\big)$
    \end{itemize} 
    descend to  well-defined and locally biholomorphic maps on $\mathbb{M}(S)$.
    \end{thm}

That is, we have the following commutative diagram:
\begin{center}
   \begin{tikzcd}
 & \mathrm M^*(S) \arrow[ld] \arrow[d] \arrow[dr] \\
T^*\mathcal T(S)\times T^*\mathcal T(\overline{S}) & \mathbb M(S) \arrow[l, "\pi"'] \arrow[r, "\mathrm{hol}"] & \chi\big(\pi_1(S), \PSL(2,\C)^2\big)
\end{tikzcd}
\end{center}
The two maps $\pi$ and $\mathrm{hol}$ on the bottom line are local biholomorphisms (actually, the spaces are defined so that the right arrow takes value in the smooth part of the character variety). 

The most remarkable aspect of Theorem~\ref{thm: local biholos intro} is that \emph{the same} quotient space of $\mathrm M^*(S)$ has well-defined maps to the two different targets, $T^*\mathcal T(S)\times T^*\mathcal T(\overline{S})$ and the character variety. The statement that the map $\pi:\mathbb M(S)\to T^*\mathcal T(S)\times T^*\mathcal T(\overline{S})$ (on the left of the diagram) is well-defined and a local biholomorphism was already proven in \cite[Proposition 6.3]{ESrank2}. In fact, the tools for Theorem~\ref{thm: local biholos intro} and more generally for the study of the space $\mathbb M(S)$  where introduced in  \cite{ESrank2}, which provided us with the fundamental technical and conceptual background for our work. 
Our contribution here, which is nonetheless essential towards the proof of Theorem~\ref{thm:main}, is that the \emph{holonomy map} $\mathrm{hol}:\mathbb M(S)\to \chi\big(\pi_1(S), \PSL(2,\C)^2\big)$ gives a local biholomorphism in the quotient.

The space $\mathbb{M}_0(S)$ that appears in Theorem~\ref{thm:main} is then a connected component of $\mathbb{M}(S)$. We explain in Section~\ref{subsec:intro analytic cont} below the importance of the connectedness condition.

\subsection{Explicit construction of the tensors}\label{subsec:intro tensors}

We have now the right set-up to discuss the construction of the tensors $(g,\mathbf i, \mathbb I, \mathbb J, \mathbb K)$ that constitute the complex hyperKähler structure of Theorem \ref{thm:main}. In fact, these have a very concrete and explicit description in terms of the maps of Theorem \ref{thm: local biholos intro}. 

First of all, the integrable almost-complex structure $\mathbf i$ is, as a consequence of Theorem \ref{thm: local biholos intro}, equal both to the pull-back via $\pi$ of the standard product almost-complex structure on
$T^*\mathcal T(S)\times T^*\mathcal T(\overline{S})$, that we denote by $(J_{T^*\mathcal T},J_{\overline{T^*\mathcal T}})$ and corresponds essentially to multiplication by $i$ on both factors, and to the pull-back via $\mathrm{hol}$ of the almost-complex structure on $\chi\big(\pi_1(S), \PSL(2,\C)^2\big)$, that we denote by $(J_\chi,J_\chi)$ since it acts as the standard almost-complex structure of the $\PSL(2,\C)$-character variety on both factors. In symbols, 

\begin{equation}
\mathbf i=\pi^*(J_{T^*\mathcal T},J_{\overline{T^*\mathcal T}})=\mathrm{hol}^*(J_\chi,J_\chi)~.
\end{equation}

We can now describe the endomorphisms $\mathbb I,\mathbb J$ and $\mathbb K$. For the first two, they are defined by the formulae:
\begin{equation}
\mathbb I=\pi^*(J_{T^*\mathcal T},-J_{\overline{T^*\mathcal T}})\qquad \qquad \mathbb J=\mathrm{hol}^*(-J_\chi,J_\chi)~.
\end{equation}

For the third endomorphism $\mathbb K$, it is of course obtained as the composition $\mathbb K=\mathbb I\mathbb J$, but it also has a geometric interpretation as follows. In Section \ref{sec:geometric interpretation of K}, we define a map $\Theta$ from an open subset in $\mathbb M_0(S)$ containing both $\mathcal{AF}(S)$ and $\mathcal{MGH}(S)$, taking values again in $\chi\big(\pi_1(S), \PSL(2,\C)^2\big)$, but different from the holonomy map. When restricted to $\mathcal{MGH}(S)$, this map takes values in $\chi\big(\pi_1(S), \PSL(2,\R)^2\big)$ and the output are the holonomy representations of the Riemannian metrics of constant curvature $-2$ obtained as the first fundamental form of two surfaces in $\Ads^3$ equidistant from a maximal surface, a construction introduced in \cite{kirill07} and studied also in \cite{BMS2013I, BMS2015II, MazSepTamb}. The map $\Theta$ then holomorphically extends this map on $\mathcal{MGH}(S)$ to an open set in $\mathbb M_0(S)$, and the output is now a pair of quasi-Fuchsian representations associated, via the Bers' Simultaneous Uniformization Theorem, to a pair of complex metrics whose curvature, which is now a priori a complex number, is identically equal to $-2$, again obtained as the first fundamental form of ``equidistant'' surfaces inside $\mathrm{PSL}(2,\C)$. On such open set, the tensor $\mathbb K$ satisfies the identity

\begin{equation}\label{eq:pullbackK}
\mathbb K=\Theta^*(-J_\chi,J_\chi)~.
\end{equation}

Now, having described the tensors $\mathbb I, \mathbb J$ and $\mathbb K$, we move on to the complex, $\mathbf i$-holomorphic symplectic forms defined by the pairings:
$$\omega_{\mathbb I}=g(\mathbb I\cdot,\cdot)\qquad \omega_{\mathbb J}=g(\mathbb J\cdot,\cdot)\qquad \omega_{\mathbb K}=g(\mathbb K\cdot,\cdot)~.$$

The metric $g$ can then be recovered from the knowledge of any tensor $\mathbb I, \mathbb J, \mathbb K$ and the corresponding complex symplectic form $\omega_{\mathbb I}, \omega_{\mathbb J}$ or $\omega_{\mathbb K}$. For the first and the last, we have the identities:

\begin{equation} \label{eq:omegaIK intro}
\omega_{\mathbb I}:= \mathrm{hol}^*\left(  \pi_r^* \omega_G+ \pi_l^*\omega_G\right)\qquad \qquad \omega_{\mathbb K}:= i\mathrm{hol}^*\left(  \pi_r^* \omega_G- \pi_l^*\omega_G\right)
\end{equation}
where $\omega_G$ here denotes the Goldman complex symplectic form on $\chi\big(\pi_1(S), \PSL(2,\C)\big)$, and $\pi_l,\pi_r$ the projections from $\chi\big(\pi_1(S), \PSL(2,\C)^2\big)$ on the left and right $\PSL(2,\C)$-factor. For $\omega_{\mathbb J}$, we instead have the following:

\begin{equation}\label{eq:omegaJ intro}
\omega_{\mathbb J}=-i\pi^*(\mathrm{pr}_l^*\Omega_{T^*\mathcal T}-\mathrm{pr}_r^*\Omega_{\overline{T^*\mathcal T}})
\end{equation}
where $\mathrm{pr}_l,\mathrm{pr}_r$ now denote  projections onto the left and right factor in 
$T^*\mathcal T(S)\times T^*\mathcal T(\overline{S})$, $\Omega_{T^*\mathcal T}=-\mathrm d\lambda$ is the complex symplectic form of ${T^*\mathcal T(S)}$ ($\lambda$ being the complex Liouville form), and $\Omega_{\overline{T^*\mathcal T}}=\overline{\Omega_{T^*\mathcal T}}$ the one for $T^*\mathcal T(\overline S)$.
Incidentally, \eqref{eq:omegaIK intro} and \eqref{eq:omegaJ intro} show that $\omega_{\mathbb I}, \omega_{\mathbb J}$ and $\omega_{\mathbb K}$ are all closed and non-degenerate, hence are complex symplectic forms.

We conclude this discussion by mentioning that, using the bi-complex algebra $\mathbb{C}_\tau=\{a+\tau b\,|\,a,b\in \C\}$ where $\tau$ is an ``imaginary'' unit such that $\tau^2=1$, we can package some of the identities together as follows. Since the algebra $\C_\tau$ is isomorphic to $\C\oplus\C$, the character variety $\chi\big(\pi_1(S), \PSL(2,\C)^2\big)$ is isomorphic to $\chi\big(\pi_1(S), \PSL(2,\C_\tau)\big)$, and this allows to define a bi-complex (i.e. $\C_\tau$-valued) Goldman symplectic form that we denote by $\omega_{G_\tau}$. In Section \ref{sec:bi-complex}, we then prove that $\omega_{G_\tau}$ is bi-holomorphic and we obtian the
following identities:

	\begin{equation}\label{eq:bicomplex identities}
		2\mathrm{hol}^* \omega_{G_\tau} = \omega_{\mathbb{I}} + i\tau \omega_{\mathbb{K}}
	\qquad \qquad
		2\Theta^* \omega_{G_\tau} = \omega_{\mathbb{I}} + i\tau \omega_{\mathbb{J}}~.
	\end{equation}
where the second holds, as before, on an open set where $\Theta$ is defined.

\subsection{Analytic continuation}\label{subsec:intro analytic cont}
We mentioned already, and implicitly used when we stated Theorem~\ref{thm:main}, that the space $\mathbb{M}_0(S)$ contains embedded copies of $\QF(S)$, $\mathcal{AF}(S)$ and $\mathcal{MGH}(S)$. (In fact, $\mathbb{M}_0(S)$ is the connected component of $\mathbb{M}(S)$ containing the images of such embeddings.) Let us now spell out these embeddings:
\begin{enumerate}[leftmargin=2em]
    \item A holomorphically embedded copy of $\QF(S)$ is included in $\mathbb{M}_0(S)$, with image the equivalence classes of 5-tuples of the form $(J_1,0,-J_2,0,0)$. It turns out that the isotopy classes of $J_1$ and $-J_2$ are precisely the Bers' parameters given by the Simultaneous Uniformization Theorem. Geometrically, 5-tuples of this form correspond to minimal immersions in a totally geodesic holomorphic submanifold of $\PSL(2,\C)$, which is a copy of the space of unoriented geodesics of $\mathbb H^3$ -- equivalently, a 2-to-1 quotient of the complement of the diagonal in $\CP^1\times\CP^1$. The minimal immersions then have the property that the projection to the first and second $\CP^1$--factor is holomorphic with respect to $J_1$ and $-J_2$ respectively.
    \item A totally real embedding of $\mathcal{AF}(S)$ is obtained as follows. Recall that a quasi-Fuchsian manifold diffeomorphic to $S\times\R$ is almost-Fuchsian if it contains a closed minimal surface (which is a posteriori unique) such that the norm of its second fundamental form is at most 2. Almost-Fuchsian manifolds are then parameterized by  triples $(J,q,u)$ up to isotopy, where $u$ is a real-valued function on $S$ solving the Gauss equation in $\mathbb H^3$ with respect to the complex structure $J$ and the holomorphic quadratic differential $q$. These are included in $\mathbb{M}_0(S)$ by sending such triple $(J,q,u)$ to the 5-tuple $(J,q,-J,\overline{q},u)$, and geometrically correspond to minimal surfaces in a totally geodesic copy of $\mathbb H^3$ in $\PSL(2,\C)$.
    \item Similarly, a totally real embedding of $\mathcal{MGH}(S)$ is obtained by sending a triple $(J,q,u)$ (with again $u$ a real-valued smooth function), that now determines a unique maximal surface in a maximal globally hyperbolic manifold (here without any assumption on the norm of the second fundamental form), to the 5-tuple $(J,q,-J,-\overline{q},u)$. Geometrically this corresponds to maximal surfaces in a totally geodesic copy of $\Ads^3\cong\PSL(2,\R)$ in $\PSL(2,\C)$. 
\end{enumerate} 

The totally real embeddings in the second and third item then play a key role in our proof of Theorem~\ref{thm:main}. Here we only outline, at a heuristic level, two fundamental ingredients of our approach:
\begin{itemize}[leftmargin=2em]
    \item From the description of Section \ref{subsec:intro tensors} above, the tensors $(g,\mathbf i, \mathbb I, \mathbb J, \mathbb K)$ are, roughly speaking, the only possible holomorphic extensions of the tensors that are provided by the hyperKähler structure of $\mathcal{AF}(S)$ and the para-hyperKähler structure of $\mathcal{MGH}(S)$.
    \item To prove that the relations that define a complex hyperKähler structure hold, we use that, from the known hyperKähler  or  para-hyperKähler structures, certain relations hold on $\mathcal{AF}(S)$  of $\mathcal{MGH}(S)$,  both of which are half-dimensional, totally real submanifold. Therefore the relations hold everywhere by uniqueness of the holomorphic extension (see Appendix \ref{appendix:extension}). 
\end{itemize}

The second point justifies why we need to take, in order to pursue our approach, $\mathbb{M}_0(S)$ to be \emph{the connected component} of $\mathbb{M}(S)$ that contains the totally real copies of $\mathcal{AF}(S)$ and $\mathcal{MGH}(S)$.

However, the copy of $\QF(S)$ described in  item (1) is a \emph{complex} submanifold, and therefore no analytic continuation argument can be used from it. In other words, we are not able to use any of the properties of the complex bi-Kähler structure of Loustau-Sanders to infer that certain properties hold globally on $\mathbb{M}_0(S)$. Still, our complex hyperKähler structure recovers precisely the Loustau-Sanders structure on the holomorphic copy of quasi-Fuchsian space.  

\subsection{Complex area functional}\label{subsec:intro area}
In the last part of our paper, we provide additional constructions on the space $\mathbb M_0(S)$, which also provide a more geometric interpretation of some of the elements that constitute the complex hyperKähler structure of Theorem \ref{thm:main}.

First, we recall that the elements of the spaces $\mathbb M(S)$, and its connected component $\mathbb M_0(S)$, can be thought as of equivalence classes of 5-tuples $(J_1,q_1,-J_2,\overline{q_2},u)$, where $J_1$ and $-J_2$ are complex structures on $S$ compatible with opposite orientations, $q_1$ and $\overline{q_2}$ are holomorphic quadratic differentials with respect to each of these complex structures, and $u$ is a smooth, $\C$-valued function on $S$ satisfying a ``Gauss equation'' (see Equation \eqref{eq: Gauss PSL(2,C)}). 

This description allows to define a holomorphic $\C^*$-action on the space of such equivalence classes of 5-tuples, as follows:
$$\lambda\cdot [J_1,q_1,-J_2,\overline{q_2},u]:=
[J_1,\lambda q_1,-J_2,\lambda^{-1}\overline{q_2},u]~,$$
for $\lambda$ in $\C^*$. It is immediately seen that this action extends the well-known $S^1$-action on $\mathcal{AF}(S)$ and $\mathcal{MGH}(S)$ given by multiplication of the holomorphic quadratic differential by a unit complex number $e^{i\theta}\in S^1$. It is moreover known from \cite{trautwein2018infinite} and \cite{MazSepTamb} that such $S^1$-action is Hamiltonian on both $\mathcal{AF}(S)$ and $\mathcal{MGH}(S)$, an Hamiltonian function being the area of the minimal/maximal surface.

Unfortunately, unlike the $S^1$-action, the ``extended'' $\C^*$-action is not globally defined on our space $\mathbb M_0(S)$. The issue is that, by technical reasons, in our definition of $\mathbb M_0(S)$ we need to impose that the 5-tuple $[J_1,q_1,-J_2,\overline{q_2},u]$ has the property that the holonomy of the corresponding minimal immersion is a doubly irreducible  representation in $\PSL(2,\C)^2$ (see Section \ref{subsec: minimal immersions}), that is, both $\PSL(2,\C)$ factors are irreducible representations. However, this condition is open, and thus the $\C^*$-action is locally defined, that is, the orbit of every point is well-defined for an open subset of $\C^*$ containing $1$.

Having established these preliminary observations, we prove:

\begin{thm}\label{thm:area intro}
    The complex area functional $\mathcal A:\mathbb M_0(S)\to\C$ is a Hamiltonian function with respect to $\omega_{\mathbb I}$, that is:
    $$\mathrm d\mathcal A=\omega_{\mathbb I}(\xi,\cdot)$$
    where $\xi$ is the infinitesimal generator of the (local) $\C^*$-action.
\end{thm}

By the complex area functional, we mean the following. Given an equivariant minimal immersion of $\widetilde S$ in $\mathrm{PSL}(2,\C)$, the area form of the immersion is a complex-valued 2-form on $\widetilde S$, associated with the (complex) first fundamental form, and descends to $S$ by equivariance. Hence one can define the area functional as the integral over $S$ of such complex-valued 2-form, and we prove (Proposition \ref{prop:area functional well defined}) that it descends to a holomorphic function $\mathcal A$ on $\mathbb M_0(S)$. 

The $\C^*$-action is also crucially used to prove that $\QF(S)$ is totally geodesic in $\mathbb M_0(S)$, since it is the fixed point of a (local) isometric action. The totally geodesic property of $\mathcal{AF}(S)$ and $\mathcal{MGH}(S)$ is proved by a more general argument that relies on the fact that they are totally real and that the restriction of the holomorphic Riemannian metric $g$ is real-valued, although it could also be proved by a similar fixed point argument for certain anti-holomorphic involutions (Remark \ref{rmk:holo involution}). 

\subsection*{Organization of the paper} 
In Section \ref{sec:prel} we collect preliminary background on holomorphic Riemannian metrics, the geometry of $\PSL(2,\C)$, basics of quasiconformal Teichmüller theory and holomorphic quadratic differentials, complex and Bers metrics. In Section \ref{sec: minimal immersions} we introduce and study minimal immersions in $\PSL(2,\C)$, leading to the definition of the spaces $\mathbb M(S)$ and $\mathbb M_0(S)$. In Section \ref{sec:holonomy} we study the holonomy map, proving Theorem \ref{thm: local biholos intro}. In Section \ref{sec:complex_hyperkahler} we introduce the tensors as described in Section \ref{subsec:intro tensors} and prove that $(g,\mathbf i, \mathbb I, \mathbb J, \mathbb K)$ is a complex hyperKähler metric, together with the three items of Theorem \ref{thm:main}. In Section \ref{sec:C*} we define the local $\C^*$-action, prove Theorem \ref{thm:area intro} and prove the ``moreover'' part of Theorem \ref{thm:main}, that is, that $\mathcal{QF}(S)$, $\mathcal{AF}(S)$ and $\mathcal{MGH}(S)$ are totally geodesic. In Section \ref{sec:geometric interpretation of K} we construct the map $\Theta$ and prove the identity \eqref{eq:pullbackK} and the second identity in \eqref{eq:omegaIK intro}. In Section \ref{sec:bi-complex} we discuss the bi-complex perspective, we prove the identities \eqref{eq:bicomplex identities}, and we describe bi-holomorphic Darboux coordinates for $\omega_{G_\tau}$ that are, in a certain sense, Fenchel-Nielsen coordinates.

We include two appendices. In the first, for the sake of completeness, we provide the statement about the uniqueness of the holomorphic extension of a real analytic function defined on a totally real, half-dimensional submanifold, which is the key statement we use in all analytic continuation arguments. In the second, we provide a formula for the linearization of the mean curvature operator for minimal surfaces in $\PSL(2,\C)$, adapting in our setting some very classical computations in the Riemannian framework.

\subsection*{Acknowledgments} 
C.E., N.R. and A.S. 
were funded by the European Union (ERC, GENERATE, 101124349). F.M. was funded by the European Union (ERC, PosLieRep, 101018839). Views and opinions expressed
are however those of the author(s) only and do not necessarily reflect those of the European Union
or the European Research Council Executive Agency. Neither the European Union nor the granting
authority can be held responsible for them. 

A.T. acknowledges support from the Italian Ministry of University and Research (MUR) through the FIS3 project “Geometry and Topology of higher Teichm\"uller spaces”, grant no. FIS-2024-05142.

\subsection*{GenAI declaration} No GenAI tool has been used in any stage of this work. All the mathematical arguments have been developed without the use of AI, and the manuscript has been entirely human written.

\section{Preliminaries} \label{sec:prel}
\noindent In this preliminary section, we will outline the theory of holomorphic Riemannian manifolds, positive complex metrics on surfaces, and Bers metrics (see \cite{BEE, ElE}). 

\subsection{Holomorphic Riemannian manifolds}
\label{section: holo Riem mflds}
\noindent Let $\mathbb X$ be a complex manifold of complex dimension $n$ and with almost complex structure $\mathbf i:T\mathbb X\to T\mathbb X$. We will denote by $T^{1,0}\mathbb{X}$ the holomorphic tangent bundle of $\mathbb{X}$. With a slight abuse of notation, we will continue to denote by $\mathbf{i}$ its $\C$--linear extension to the complexified tangent bundle $\C T\mathbb{X}$, which thus acts on $T^{1,0}\mathbb{X}$ by multiplication by $i$.

    

\begin{defn}\label{def:holomorphic_riemannian_metric}
A \emph{holomorphic Riemannian metric} $\inners$ on $\mathbb X$ is a holomorphic section of $\mathrm{Sym}^2(T^{*1,0}\mathbb X)$ that at every point gives a non-degenerate bi-linear form on $T^{1,0}\mathbb{X}$. The pair $(\mathbb X, \inners)$ is called a \emph{holomorphic Riemannian manifold}.
\end{defn}

Such objects have been studied by several authors (see, for instance, \cite{dumitrescu2001, dumitrescu2009}). 
Since every smooth section of $T^{* 1,0}\mathbb{X}$ vanishes identically on $T^{0,1}\mathbb{X}$, the data of a holomorphic Riemannian metric $\langle \cdot, \cdot \rangle$ determines a degenerate bilinear form on $\mathbb{C} T \mathbb{X}$ with respect to which $(\mathbb{C} T \mathbb{X})^\perp = T^{0,1}\mathbb{X}$.
Moreover, for a given holomorphic Riemannian metric $\inners$ on $\mathbb X$, both its real part $\Re\inners$ and its imaginary part $\Im\inners$ are bilinear forms of signature $(n,n)$ on the (real or holomorphic) tangent space at every point of $\mathbb X$.

Holomorphic Riemannian manifolds $(\mathbb X, \inners)$ have a natural notion of \emph{holomorphic Levi-Civita connection $D$}:

\begin{defn}\label{def:homorphic levi-civita}
The \emph{holomorphic Levi-Civita connection} of a holomorphic Riemannian manifold $(\mathbb X, \inners)$ is the unique affine connection on $T\mathbb X$ whose natural $\C$--linear extension
\[
D : \Gamma(\mathbb{C} T\mathbb{X}) \times \Gamma(\mathbb{C} T\mathbb{X}) \longrightarrow \Gamma(\mathbb{C} T\mathbb{X}), \qquad D_U V := D(U,V)~,
\]
satisfies the following properties:
\begin{enumerate}
    \item $T^{1,0} \mathbb{X}$ is $D$--invariant, i.e. $D_U(\Gamma(T^{1,0} \mathbb{X})) \subset \Gamma(T^{1,0} \mathbb{X})$ for every $U \in \Gamma(\C T\mathbb{X})$;
    \item The restriction of $D$ over $\Gamma(T^{1,0} \mathbb{X})$ is torsion-free, i.e.
    \[
    D_{V}W - D_W V = [V,W] 
    \]
    for all $V, W \in \Gamma(T^{1,0}\mathbb{X})$.
    \item For every holomorphic local section $Z$ of $T^{1,0}\mathbb{X}$ and for any $U \in \Gamma(\C T\mathbb{X})$,
    \[
    D_{\mathbf{i}U} Z = \mathbf{i}(D_U Z)~.
    \]
    \item $D$ is compatible with the metric $\langle \cdot, \cdot \rangle$, i.e.
    \[
    W(\langle U, V \rangle) = \langle D_W U, V \rangle + \langle U, D_W V \rangle
    \]
    for any $U, V, W \in \Gamma(\C T \mathbb{X})$.
\end{enumerate}
\end{defn}

\begin{remark}
    Property (1) is equivalent to asking that the ($\C$--linear extension of the) almost-complex structure $\mathbf{i}$ is $D$--parallel. Properties (1) and (2) imply that, given any section $V \in \Gamma(T^{1,0} \mathbb{X})$, the $(0,1)$--part of $D_\bullet V$ coincides with $\bar{\partial} V$, where $\bar{\partial}$ denotes the del-bar operator of the holomorphic vector bundle $T^{1,0}\mathbb{X}$.
\end{remark}

For a proof of the existence and uniqueness of the holomorphic Levi-Civita connection, see \cite{lebrun1983}.

As in the Riemannian setting, starting from $D,$ one can define a $(0,4)$--type curvature tensor $R$, and, for all $2$--dimensional complex vector subspaces $V< T_p \mathbb X$ on which $\inners|_{V}$ is non-degenerate, one can define the sectional curvature of $V$. Moreover, since $D$ coincides with the Levi-Civita connection of both $\Re\inners$ and $\Im\inners$, the complex curvature tensor $R$ satisfies the same symmetries as in the real Riemannian case. 

\subsection{Fundamental examples}

The simplest example of a holomorphic Riemannian manifold is $\mathbb C^n$ equipped with $\langle\underline z, \underline w\rangle=\sum_{k=1}^n z_k w_k$, for $\underline z =(z_1,\dots, z_n),\underline w =(w_1,\dots, w_n)\in \C^n$. In fact, one can see that
$(\mathbb \C^n, \inners)$ is the unique geodesically-complete simply connected holomorphic Riemannian manifold with constant sectional curvature zero up to isometry (see \cite[Theorem 2.6]{BEE}). We now introduce two explicit examples of holomorphic Riemannian manifolds that will play a fundamental role in this article. 

\subsubsection{The space of oriented geodesics of \texorpdfstring{$\mathbb H^3$}{H3}}\label{sec:space_of_oriented_geodesics}
We will consider two holomorphic Riemannian manifolds of complex dimension two and constant sectional curvature $-1$, which are defined as follows:
\begin{itemize}
    \item $\widetilde{\mathbb G}$ is the space of \emph{oriented} (unparametrized) geodesics of $\mathbb H^3$;
    \item ${\mathbb G}$ is the space of \emph{unoriented} (unparametrized) geodesics of $\mathbb H^3$.
\end{itemize}
Clearly, $\widetilde{\mathbb G}$ (which is simply connected) is a 2-to-1 cover of ${\mathbb G}$ via the map that forgets the orientation. Now, $\widetilde{\mathbb G}$ is endowed with a complex structure via the identification with 
$\CP^1\times\CP^1\setminus \Delta$, where $\Delta$ is the diagonal in the product, through the map that sends an oriented geodesic to its endpoints in $\CP^1\cong \partial \Hyp^3$. Moreover, one can define  a holomorphic Riemannian metric on $\widetilde{\mathbb G}$ as follows: given an affine chart $(U,z)$ for $\CP^1$, we define the metric on $U\times U\setminus \Delta$ through the chart $z\times z=(z_1,z_2)$ as
\[
\inners_{\widetilde{\mathbb G}} := -\frac 4 {(z_1-z_2)^2} \mathrm{d}z_1\cdot \mathrm{d}z_2~.
\]
As shown in \cite[Proposition 2.11]{BEE}, this description is independent of the choice of the affine chart, and it defines a holomorphic Riemannian metric on $\widetilde{\mathbb G}$. Moreover, $(\widetilde{\mathbb G}, \inners_{\widetilde{\mathbb G}})$ is the unique simply connected geodesically-complete holomorphic Riemannian manifold of constant (sectional) curvature $-1$ and complex dimension 2. 

The action of $\mathbb Z_2$ on $\widetilde{\mathbb G}$ that switches the orientation  corresponds to the action on $\CP^1\times\CP^1\setminus \Delta$ generated by the  involution $(z_1,z_2)\mapsto(z_2, z_1)$. It is thus easily checked that this action preserves both the complex structure and the holomorphic Riemannian metric. Hence $\mathbb G=\widetilde{\mathbb G}/\mathbb Z_2$ inherits a holomorphic Riemannian metric of constant sectional curvature $-1$.

The identity component of the isometry group of both $\widetilde{\mathbb G}$ and ${\mathbb G}$ is the group $\PSL(2,\C)$, acting via the natural action on the (oriented or unoriented) geodesics of $\mathbb H^3$.

\subsubsection{The holomorphic Riemannian manifold \texorpdfstring{$\mathrm{PSL}(2,\mathbb{C})$}{PSL(2,C)}}\label{sec: PSL(2,C)} Recall that, given a complex Lie group $G$, its Lie algebra $\mathfrak g$ is defined as the vector space of left-invariant holomorphic vector fields, and it is naturally identified with $T^{1,0}_{\mathrm{id}}G$. When $G$ is semi-simple, the complex Killing form is a non-degenerate $\C$--bilinear form on $\mathfrak g$. 
Since it is also $Ad$--invariant, we can extend it to a bilinear form $Kill$ on the whole holomorphic tangent bundle $T^{1,0}\PSL(2,\C)$ by push-forward by left (or, equivalently, right) multiplication. This defines a holomorphic Riemannian metric on $G$.

Now, applying this construction to the complex Lie group $\PSL(2,\C)$, we obtain a complete holomorphic Riemannian manifold with constant sectional curvature. It will be more convenient to replace the Killing form by the $\C$--bilinear form $\inners:=\frac 1 8 Kill$, one of the advantages being that $(\PSL(2,\C), \inners)$ has constant sectional curvature equal to $-1$, as shown in \cite[Theorem 2.9]{BEE}. The group $\PSL(2,\C)\times\PSL(2,\C)$ acts by isometries by left and right multiplication and can be identified with the connected component containing the identity of the isometry group of $(\PSL(2,\C), \inners)$.  

As shown in \cite[Remark 2.7.1]{bonsante2020}, $\PSL(2,\C)$ admits a global holomorphic 3-form $\Omega$ such that $\|\Omega\|_{\inners}^2=1$, namely a metric volume form. Every isometry pushes $\Omega$ to either $\Omega$ or $-\Omega$. The space $\PSL(2,\C)\times \PSL(2,\C)$ corresponds in fact to the index-2 subgroup of isometries of $\PSL(2,\C)$ that preserve $\Omega$.

\subsubsection{Totally geodesic submanifolds}\label{subsubsec:tot geod}
There are three important totally geodesic submanifolds of $\PSL(2,\C)$, that we now describe.
\begin{enumerate}[leftmargin=2em]
\item The space $\mathbb G$ (introduced in Section~\ref{sec:space_of_oriented_geodesics}) is isometrically embedded in $\PSL(2,\C)$ via the map that associates to an unoriented geodesic $\gamma$ of $\mathbb H^3$ the unique elliptic orientation-preserving isometry of $\mathbb H^3$ that fixes $\gamma$ and has order two. The image of the embedding is  the space of projective classes of traceless matrices in $\PSL(2,\C)$, and is totally geodesic. 

The Lie group embedding of the isometry group of $\mathbb G$ into the isometry group of $\PSL(2,\C)$ associated to such embedding of $\mathbb G$ is the diagonal embedding:
\[
\Delta:\PSL(2,\C)\to\PSL(2,\C)\times \PSL(2,\C)~.
\]


\item The hyperbolic space $\mathbb H^3$ is itself embedded isometrically and in a totally geodesic way inside $\PSL(2,\C)$ by seeing it as its associated Riemannian symmetric space: taking a corresponding Cartan splitting $\mathfrak{sl}_2(\C)=\mathfrak K \oplus \mathfrak p$, the holomorphic metric $\inners$ restricts to a Riemannian metric on $\mathfrak p$, and $\mathbb H^{3}\cong \exp(\mathfrak p)=\{A\in \PSL(2,\C)\ |\ A=\ \overline{^tA^{-1}}\}$.

This copy of $\mathbb H^3$ can be interpreted geometrically as the subset of $\PSL(2,\C)$ whose elements are purely loxodromic isometries sending a fixed ``origin'' $o$ to a given point $x$ in $\mathbb H^3$. 

The corresponding embedding at the level of isometry groups is  
\[
\Delta^\dagger:\PSL(2,\C)\to\PSL(2,\C)\times \PSL(2,\C)~,
\]
where we recall that $\Delta^\dagger(A)=(A,\overline{{}^t A^{-1}})$.

    \item The inclusion $\PSL(2,\R)\subset \PSL(2,\C) $ is a totally geodesic isometric immersion of the real Lie group $\PSL(2,\R)$, endowed with the $\PSL(2,\R)$--invariant pseudo-Riemannian metric of signature $(2,1)$ (that is, a Lorentzian metric) induced by $(1/8)$ times the (real) Killing form on the (real) Lie algebra. The restriction of $\langle\cdot,\cdot\rangle_A$ to $T_A\PSL(2,\R)$ is real and coincides with this Lorentzian metric.
    
    It is known that this $3$--dimensional Lorentzian manifold is isometric to the Anti-de Sitter space $\Ads^3$  (see, for instance, \cite[\S~3.1]{bonsante2020}). The identity component of the isometry group of $\Ads^3$, via this $\PSL(2,\R)$ model, is isomorphic to $\PSL(2,\R)\times\PSL(2,\R)$, and the embedding of the isometry groups is  the obvious real embedding:
    \[\iota_\R:\PSL(2,\R)\times \PSL(2,\R)\to \PSL(2,\C)\times \PSL(2,\C)~.
    \]

\end{enumerate}

In the following, by an abuse of notation, when we mention the spaces $\mathbb G$, $\mathbb{H}^3$ and $\Ads^3$, we will always consider them as totally geodesic isometrically embedded submanifolds of $\PSL(2,\C)$. 


\subsection{The space of complex structures and the Teichmüller space}\label{sec: beltrami differentials}

Let $S$ be a connected, closed (compact with empty boundary), oriented surface with negative Euler characteristic. By the Newlander–Nirenberg theorem, we can view a complex structure on $S$ as an almost-complex structure $J$ compatible with the orientation
in the following sense: there exists a point $p \in S$ and a tangent vector $v \in T_p S$, with $v \neq 0$, such that $(v, J_p v)$ is a positive oriented basis of $T_p S$. (In fact, if such $p$ and $v$ exist, then the same property holds for any other point and non-zero tangent vector in $S$.) In this framework, a holomorphic chart on $(S,J)$ is a map $z : U \to \mathbb{C}$ that is compatible with the orientations of $S$ and $\mathbb{C}$ and whose differential is $\mathbb{C}$-linear with respect to $J$.

    We will denote:
\begin{itemize}
    \item The space of complex structures on $S$ by $\mathcal{C}(S)$.
    \item The \emph{Teichmüller space} of $S$, that is, the quotient of $\mathcal C(S)$ under the pull-back action of $\Diff_0(S)$, the group of diffeomorphisms isotopic to the identity, by $\mathcal{T}(S)$.
\end{itemize}

It is a classical fact that $\mathcal{T}(S)$ has a structure of real manifold of dimension $3|\chi(X)|$. It is also well-known that $\mathcal{T}(S)$ (resp. $\mathcal{C}(S)$) is endowed with a complex (resp. complex Fréchet) structure, that we now introduce.

\subsubsection{Beltrami differentials} For a fixed complex structure $J_0 \in \mathcal{C}(S)$, let $\mathrm{Belt}_{J_0}(S)$ be the space of Beltrami differentials on $(S,J_0)$, that is, the space of sections
\[
\mathrm{Belt}_{J_0}(S) : = \Gamma(T^{1,0}_{J_0} S \otimes T^{0,1*}_{J_0} S)~,
\]
where $T^{1,0}_{J_0} S$ (resp. $T^{0,1*}_{J_0} S$) denotes the holomorphic tangent bundle (resp. anti-holomorphic cotangent bundle) of $(S,J_0)$. A Beltrami differential $\eta \in \mathrm{Belt}_{J_0}(S)$ can be expressed in a holomorphic chart $(U,z)$ as $\eta = \eta_z \, \frac{\mathrm{d}\overline{z}}{\mathrm{d}z}$, and the locally defined function $\eta_{z}$ transforms on overlapping complex coordinate charts $(U,z)$ and $(V,w)$ as
$$
\eta_z = \eta_w \, \frac{\big(\overline{\frac{\partial w}{\partial z}}\big)}{\big(\frac{\partial w}{\partial z}\big)}~.
$$

Recall that there is a natural $L^\infty$--norm on the space of Beltrami differentials given by $\|\eta\|_{L^\infty} := \sup_{S} |\eta|$, where $|\eta| : S \to \mathbb{R}_{>0}$ denotes the non-negative function $|\eta|(p) := |\eta_z(p)|$, as we write $\eta = \eta_z \, \tfrac{\mathrm{d}\overline{z}}{\mathrm{d}z}$. (Here $|\eta|(p)$ does not depend on the choice of the chart at $p$ due to the aforementioned coordinate-change rule.) In what follows, we will denote by $\mathrm{Belt}^1_J(S)$ the open unit ball inside the normed space $(\mathrm{Belt}_J(S), \| \cdot \|_{L^\infty})$.

Given two complex structures $J_0, J \in \mathcal{C}(S)$ and an orientation-preserving diffeomorphism $f : (S,J_0) \to (S,J)$, the \emph{Beltrami differential of $f$}, denoted here by $\mu(f)$, measures how far $f$ is from being a biholomorphism between $(S,J_0)$ and $(S,J)$. The tensor $\mu(f)$ can be expressed locally as
\[
\mu(f) = \frac{\partial_{\overline{z}}(\zeta \circ f)}{\partial_z(\zeta \circ f)} \, \frac{\mathrm{d}z}{\mathrm{d}\overline{z}}~, 
\]
where $(U_0,z)$ and $(U,\zeta)$ are local charts for $(S,J_0)$ and $(S,J)$, respectively, such that $f(U_0) \subset U$. Notice that $\partial_z(\zeta \circ f)$ is no-where vanishing because $f$ preserves the orientation of $S$, and $\mu(f)$ is independent of the choice of local charts. Furthermore, $\mu(f)$ satisfies $\|\mu(f)\|_{L^\infty} < 1$.

The space $\mathcal C(S)$ admits a natural structure of complex Fréchet manifold that can be described as follows. Fixed $J_0\in \mathcal C(S)$, consider the Beltrami differential $\mu = \mu_J \in $ of the identity map $id\colon (S, J_0)\to (S,J)$, as $J$ varies in $\mathcal C(S)$. By the Generalized Riemann mapping theorem, the map $J \mapsto \mu_J$ defines a one-to-one correspondence
\[
\mathcal{C}(S) \longrightarrow \mathrm{Belt}_{J_0}^1(S)~. 
\]

By the chain rule for Beltrami differentials, the pull-back complex Fréchet structure does not depend on the choice of the fixed $J_0$.
This provides a model for the holomorphic tangent space at $J_0\in \mathcal C(S)$ as
$$
T_{J_0}^{1,0}\mathcal C(S)\cong \mathrm{Belt}_{J_0}(S)~.
$$

The Teichmüller space $\mathcal T(S)$ then has a natural complex (finite-dimensional) manifold structure for which the projection $\mathcal C(S)\to \mathcal T(S)$ is holomorphic (for more details see \cite[\S 4]{tromba1992}).

\subsubsection{Alternative description of the almost-complex structure}
Observe that, formally, a tangent vector $\dot{J} \in T_{J}\mathcal{C}(S)$ is naturally a section of $\mathrm{End}(TS)$ that anti-commutes with $J$.
The induced almost-complex structure on $T\mathcal C(S)$ is very easy to describe. The following result makes this description precise (see, for instance, \cite[\S 1.1]{tromba1992} and \cite[\S 2.1]{BMS2015II}). 

\begin{prop}
    The natural identification between $T_J\mathcal C(S)$ and $T_{J_0}^{1,0}\mathcal C(S)\cong \mathrm{Belt}_{J_0}(S)$ corresponds to the map
    \begin{equation}\label{eq:ident mu Jdot}
    \dot J\mapsto \dot\mu=\frac{1}{4i}\left(\dot{J} - i J \dot{J}\right)~.
    \end{equation}
    As a consequence, the almost-complex structure of $\mathcal{C}(S)$ acts on $T_J\mathcal{C}(S)$ by
        \begin{equation}\label{eq:J to Jdot}
    \dot J\mapsto J\dot{J}~.
    \end{equation}
\end{prop}
\begin{proof}
Any local complex coordinate $z$ for $(S,J)$ satisfies
\begin{equation}\label{eq:def_complex_chart}
    \mathrm{d}z\circ J= i\,\mathrm{d}z~, \qquad \mathrm{d}\overline z \circ J=-i\,\mathrm{d}\overline z~,
\end{equation}
or, equivalently,
\[
J(\partial_{\overline z})=-i\,\partial_{\overline z}~, \qquad J(\partial_z)=i\,\partial_z~.
\]
Now, a smooth path of complex structures on $S$ corresponds to a smooth $1$-parameter family of endomorphisms $(J_{t})_t$ and, by the Generalized Riemann mapping theorem, we can find a corresponding family of local coordinates $z_t$ for $(S,J_t)$ varying regularly in $t$. Differentiating Equation~\eqref{eq:def_complex_chart} at $t=0$, we obtain
\begin{align*}
    \mathrm{d}\dot z\circ J+ \mathrm{d}z\circ \dot{J} &= i \,\mathrm{d}\dot z~, \\
    \mathrm{d}\dot {\overline z}\circ J+ \mathrm{d}\overline z\circ \dot{J} &= -i\, \mathrm{d}\dot {\overline z}~,
\end{align*}
and, evaluating in $\partial_{\overline z}$ e $\partial_z$, we deduce
\begin{align}\label{eq: descrizione dot J}
   \mathrm{d}z(\dot{J}\partial_{\overline z})&= 2i\,\partial_{\overline z} \dot z~, & \mathrm{d}\overline z (\dot{J}\partial_{\overline z})& = 0~, \\
   \mathrm{d}z(\dot{J}\partial_z) & = 0~, & \mathrm{d}\overline z(\dot{J}\partial_z) & = -2i \,\partial_z\dot{\overline z}~. \notag
\end{align}
By construction, the tensor
\begin{equation}
\label{eq: diff di Beltrami}
    \dot{\mu} := \frac{d}{dt} \left( \frac{\partial_{\overline z} z_t}{\partial_z z_t} \frac {\mathrm{d}\overline z}{\mathrm{d}z} \right)= \partial_{\overline z}\dot z \frac {\mathrm{d}\overline z}{\mathrm{d}z} = \frac 1 {2i} \mathrm{d}z(\dot{J}\partial_{\overline z})  \frac {\mathrm{d}\overline z}{\mathrm{d}z}:=\dot{\mu}_{z}\,\frac {\mathrm{d}\overline z}{\mathrm{d}z}~,
\end{equation}
is the first-order variation at $t = 0$ of the Beltrami differential of the identity maps $id : (S,J_0) \to (S,J)$ and, in particular, it is an element of $\mathrm{Belt}_{J_0}(S)$. From Equation~\eqref{eq: descrizione dot J}, we deduce that $\dot{J}$ can be written in local coordinates as
\[
\dot{J} = \mathrm{d}\overline{z}(\dot{J} \partial_z) \, \partial_{\bar{z}} \otimes \mathrm{d}z + \mathrm{d}z(\dot{J} \partial_{\bar{z}}) \, \partial_z \otimes \mathrm{d}\overline{z}~.
\]
Similarly, 
\[
- i \, J \dot{J} = - \mathrm{d}\overline{z}(\dot{J} \partial_z) \,\partial_{\bar{z}} \otimes \mathrm{d}z + \mathrm{d}z(\dot{J} \partial_{\bar{z}}) \,\partial_z \otimes \mathrm{d}\overline{z}~.
\]
Therefore, combining the two expressions above with Equation~\eqref{eq: diff di Beltrami}, we have
\begin{align*}
    \frac{1}{2}\left(\dot{J} - i J \dot{J}\right) & = \mathrm{d}z(\dot{J} \partial_{\bar{z}}) \, \partial_z \otimes \mathrm{d}\overline{z}=2i\,\dot{\mu}~.
\end{align*}
    This proves the formula \eqref{eq:ident mu Jdot}. The formula \eqref{eq:J to Jdot} then follows from the immediate check that multiplication of $\dot\mu$ by $i$ corresponds to pre-composition of $\dot J$ by $J$.
\end{proof}

\subsubsection{Opposite orientation} Denote by $\overline S$ the same surface endowed with the opposite orientation. The space $\mathcal C(\overline S)$ of complex structures on $\overline S$ is endowed with the complex Fréchet structure for which the conjugation map
\begin{align*}
    \mathcal C(S)&\to \mathcal C(\overline S)\\
    J &\mapsto -J
\end{align*}
is anti-holomorphic. Notice that, in the identification $T_{-J}^{1,0}\mathcal C(\overline S)\cong \Gamma(T^{1,0}_{-J}S\otimes T^{0,1 *}_{-J} S)$, the differential map gives $\dot \mu\mapsto \dot{\overline \mu}$, hence the complex structure on $\mathcal C(\overline S)$ is just given by the multiplication times $i$, namely $\dot{\overline \mu}\mapsto i\dot{\overline \mu}$.

Similarly, $\mathcal T(\overline S)=\mathcal C(\overline S)/\Diff_0(S)$ is a complex manifold for which the map $[J]\mapsto [-J]$ is antiholomorphic.

\subsection{The space of holomorphic quadratic differentials}\label{sec: hq differentials} 
Another important object in this context is the bundle $\mathrm{QD}(S)$ over $\mathcal{C}(S)$ whose fiber at any point $J \in \mathcal{C}(S)$ coincides with the vector space $\mathrm{QD}_J(S)$ of holomorphic quadratic differentials with respect to $J$, that is, all sections of the holomorphic line bundle $\mathcal{K}_J^2\to(S,J)$, which can be written locally as $q = \phee(z)\,\mathrm{d}z \otimes \mathrm{d}z$ with $\partial_{\bar z} \phee = 0$. The space $\mathrm{QD}(S)$ comes with a natural complex Fréchet manifold structure for which the inclusion 
\[
\mathrm{QD}(S)=\{(J, q)\in \mathcal C(S)\times \Gamma(\mathrm{Sym}_2(TS, \C))\ |\ q \in H^0(\mathcal K_J^2) \}\hookrightarrow \mathcal C(S)\times \Gamma(\mathrm{Sym}_2(TS, \C))
\]
is holomorphic as a map between Fréchet manifolds. We will therefore denote elements in $T_{(J, q)}^{1,0}\mathrm{QD}(S)$ as pairs $(\dot \mu, \dot q)$, with the complex Fréchet structure corresponding to the multiplication times $i$ on each factor.

The pull-back action of $\mathrm{Diff}_0(S)$ on $\mathrm{QD}(S)$ determines a quotient $\mathbb{QD}(S)$, which carries a natural structure of vector bundle over Teichmüller space, whose fiber over $[J]\in\mathcal{T}(S)$ is isomorphic, as a vector space, to the fiber of $\mathrm{QD}(S)$ over any representative $J \in [J]$. In fact, $\mathbb{QD}(S)$ inherits a complex manifold structure for which the projection $\mathrm{QD}(S)\to \mathbb{QD}(S)$ is holomorphic and the bundle $\mathbb{QD}(S)\to \mathcal T(S)$ is a holomorphic vector bundle.

There is a natural pairing between quadratic and Beltrami differentials given by \begin{equation}\label{eq:pairing_QD_Belt}
\begin{aligned}
\mathrm{QD}_J(S)&\times\mathrm{Belt}_J(S)\longrightarrow\C \\ & (q,\mu)\longmapsto q\cdot\mu= \int_S \phee \mu_z\, \frac{i}{2}\mathrm d z\wedge \mathrm{d}\overline z~,
 \end{aligned}
 \end{equation}
where locally $q=\phee \,\mathrm{d}z\otimes\mathrm{d}z$ and $\mu=\mu_z\,\frac{\mathrm{d}\overline z}{\mathrm d z}$. The vector subspace
\[
\mathcal V:=\{\mu \in \mathrm{Belt}_J(S) \ | \ q\cdot \mu=0, \ \text{for all} \ q\in\mathrm{QD}_J(S)\}
\]
coincides precisely with the kernel of the differential of the quotient map $\pi:\mathrm{Belt}^1_J(S)\to\mathcal{T}(S)$ at $J\in\mathcal C(S)$. As a consequence, we can identify the holomorphic tangent space $T^{1,0}_{[J]}\mathcal T(S)$ with $\mathrm{Belt}_J(S)/\mathcal V =: \mathrm{Belt}([J])$, and, using the pairing above, we can also identify the holomorphic cotangent space $T_{[J]}^{* 1,0} \mathcal T(S)$ with $\mathrm{QD}_J(S)$.

Later in the paper, similar to what has been said for $\mathcal C(S)$ and \(\mathcal{T}(\overline S)\), we will also be interested in studying the spaces $\mathrm{QD}(S)$ and $\mathbb{QD}(\overline{S})$, equipped with the structures of complex Fréchet manifold and complex manifold respectively for which the maps
\begin{align*}
\mathrm{QD}(S)\to \mathrm{QD}(\overline S) \qquad\qquad (J,q)\mapsto (-J, \overline q)& \qquad \text{and}\\ 
\mathbb{QD}(S)\to \mathbb{QD}(\overline S) \qquad\qquad [J,q]\mapsto [-J, \overline q]&
\end{align*}
are antiholomorphic. 

Notice that the inclusion $\mathrm{QD}(\overline S)\hookrightarrow \mathcal C(\overline S)\times \mathrm{Sym}_2(TS, \C)$ is holomorphic, namely on $T^{1,0}\mathrm{QD}(\overline S)$ the complex structure writes as $(\dot{\overline \mu}, \dot{\overline q})\mapsto (i\dot{\overline \mu}, i \dot{\overline q})$. This description passes to the quotient on $T^{1,0}\mathbb{QD}(\overline S)$. Throughout the paper, we will often work with $\mathrm{QD}(S)\times\mathrm{QD}(\overline S)$ and $\mathbb{QD}(S)\times\mathbb{QD}(\overline S)$ with the standard product complex structures.


\subsection{Positive complex metrics}
In this section we describe some aspects of the theory of (positive) complex metrics on surfaces, as developed in \cite[\S~6.1]{BEE}. For the purposes of this section, it is sufficient to work with a smooth oriented surface even though in the subsequent sections we will work with a closed surface with negative Euler characteristic. For this reason, we will continue to denote the surface under consideration by $S$. Let us denote with $\C TS$ (resp. $\C T^*S$) the complexified tangent (resp. cotangent) bundle of $S$. 

\begin{defn}
A \emph{complex metric} $g$ on a smooth (real) surface $S$ is a smooth section of $\mathrm{Sym}^2(\mathbb C T^*S)$ which is pointwise non-degenerate. More precisely, $g$ is a non-degenerate symmetric $\C$--bilinear form defined on $\C T_pS$, for all $p\in S$, which varies smoothly with the point. 
\end{defn}
\noindent From the definition, it is clear that any Riemannian metric on $S$ defines a complex metric by taking its $\C$--bilinear extension. Despite these objects are clearly referring to a complex framework, they are very different from holomorphic Riemannian manifolds, introduced in Section~\ref{section: holo Riem mflds}, as we are not fixing a complex structure on $S$, and complex metrics are not required to be holomorphic in any sense. Nevertheless, as we will see in Section~\ref{subsec:immersions_into_PSL(2,C)}, these objects arise naturally from immersions of manifolds into holomorphic symmetric spaces. By applying the Gram-Schmidt algorithm, one can easily see that for every complex metric $g$ there exist local orthonormal frames $(e_1,e_2)$, namely $e_1$ and $e_2$ are local complex vector fields such that $g(e_j, e_k)=\delta_{jk}$.
Every complex metric $g$ has a natural \emph{complex Levi-Civita connection}, which is defined as the unique affine connection $\nabla^g: \Gamma( \C TM) \to \Gamma(\mathrm{End}(\C TM))$ that is torsion-free and compatible with the metric, i.e., $\nabla^g g=0$. Using $g$ and $\nabla^g$, one can define in a very natural way most of the classic tensors that are considered in Riemannian geometry, such as:
\begin{itemize}
     \item the complex Riemann curvature tensors
     \begin{align}
         \quad  R^g(X,Y)Z&=\nabla^g_X\nabla^g_Y Z-\nabla^g_Y\nabla^g_X Z- \nabla^g_{[X,Y]}Z~,\\
     \quad  R^g(X,Y,Z,T)&= g(R(X,Y)Z,T),\quad \text{for all} \ X,Y,Z,T\in \Gamma(\C TS)~;
     \end{align}
    \item the complex Gauss curvature 
    $$\K^g=\frac{R^g(X,Y,Y,X)}{g(X,X)g(Y,Y)- (g(X,Y))^2}~,$$ 
    equivalently defined for any linearly independent $X,Y \in \Gamma(\C TS) $;
    \item the gradient $\nabla^g f$ of a function $f$, defined by $\mathrm{d}f=g(\nabla^gf,\cdot)$;
    \item the Hessian of a function $f$: $$\Hess^g(f) (X,Y)= g(\nabla^g_X \nabla^gf, Y)~;$$
    \item the Laplacian of a function $f$: $$\Delta_g f= \tr_g(\Hess^g(f))= \Hess^g(f)(e_1,e_1)+\Hess^g(f)(e_2,e_2)~,$$
    for $(e_1,e_2)$ a local $g$--orthonormal basis.
\end{itemize}

\noindent Many of the properties of these tensors that hold for Riemannian metrics hold for complex metrics as well, e.g. the Riemann curvature tensor and its derivative have the same symmetries as the corresponding notions in Riemannian geometry. 

In particular, we point out that some classic formulas for conformal metrics hold, as in the following. 

\begin{defn}\label{def conformal metrics}
    Two complex metrics $g$ and $\widehat g$ on a smooth surface $S$ are \emph{conformal} if there exists a smooth function $\varrho\colon S\to \mathbb C^*$ such that $\hat g=\varrho \, g$.
\end{defn}

Given two conformal metrics $g$ and $\widehat g$ as in Definition~\ref{def conformal metrics}, it turns out that their Laplacian operators satisfy $\Delta_{\widehat g}= \frac 1 {\varrho} \Delta_g$, and that their Gauss curvatures are related by the following identity:
    \begin{equation}\label{eq: conf curvature} 
    \mathrm K^{\widehat g}= \frac 1 \varrho \big(\mathrm K^g-\frac 1 2\Delta_g \log(\varrho)\big)~.
   \end{equation}
Notice that, while $\log(\varrho)$ is only defined locally, $\Delta_g \log(\varrho)$ is well-defined globally (see e.g. \cite[Corollary 3.18]{rungi2024complex} and \cite[\S~5.1]{emam2024bers}).

We now impose an additional condition to complex metrics on $S$ which will determine an interesting subfamily. By standard linear algebra, a complex metric $g$ on $S$ has exactly two isotropic directions at each point: namely, for each point $p\in S$, the set $\{v\in \C T_pS\ |\ g_p(v,v)=0 \}$ consists of the union of two complex lines of $\mathbb C T_p S$, i.e. to two points in $\mathbb P_{\C}\big(\C T_pS\big)$. Also notice that $\mathbb P_\R\big(T_pS\big)$ embeds into the sphere $\mathbb P_{\C}\big(\C T_pS\big)$ as a round circle, cutting $\mathbb P_{\C}\big(\C T_pS\big)$ into two disks.
\begin{defn}\label{def: positive}
    A complex metric $g$ is \emph{positive} if for each $p\in S$, the two isotropic directions are contained in distinct components of $\mathbb P_{\C}(\C T_pS)\setminus \mathbb P_\R(T_pS)$.
\end{defn}
\noindent Observe that the orientation on $S$ induces an orientation to the circle 
$\mathbb P_\R(T_pS)\subset \mathbb P_{\C}(\C T_pS)$, allowing to distinguish a ``positive'' and a ``negative'' connected component of $\mathbb P_{\C}(\C T_pS)\setminus \mathbb P_\R(T_pS)$. Again, Riemannian metrics can be seen as positive complex metrics just by extending them by $\C$--linearity: its isotropic directions are the eigenlines of the induced almost complex structure, whose $(+i)$ and $(-i)$--eigenlines lie in the positive and in the negative component of $\mathbb P_{\C}(\C T_pS)\setminus \mathbb P_\R(T_pS)$, respectively. We can actually generalize this picture to positive complex metrics as follows.

Let $g$ be a positive complex metric and let $\ell_+(g_p),\ell_-(g_p)$ be the two isotropic lines in $\mathbb P_\R\big(T_pS\big)\subset \mathbb P_{\C}(\C T_pS)$ that lie respectively in the positive and in the negative connected component. As already observed, since the two directions are distinct and vary smoothly with $p$, they define complex line subbundles such that $\C TS=\ell_+(g)\oplus\ell_-(g).$ They also determine unique complex structures $J_1(g)\in\mathcal C(S), -J_2(g)\in\mathcal C(\overline S),$ characterized by \begin{equation}\label{eq: complex-structures-isotropic-lines}
T^{0,1}_{J_1(g)}S=\ell_-(g)~,\quad T^{0,1}_{-J_2(g)}S=\ell_+(g)~.
\end{equation}
Notice that this construction only uses the orientation of \(S\), and does not require \(S\) to be closed. The same data can equivalently be encoded by a bi-complex structure. \begin{defn}\label{def: bicomplex structure} A {bi-complex structure} on $\C TS$ is a $\C$--linear bundle endomorphism $\mathcal J:\C TS\to\C TS$ such that:
\begin{itemize}
    \item $\mathcal J^2=-\mathrm{Id}$;
    \item for all $p\in S$, the eigenspaces $V_i(\mathcal J_p)$ and $V_{-i}(\mathcal J_p)$ of $\mathcal J$ relative to the eigenvalues $\pm i$ lie respectively in the positive and in the negative connected component of $\mathbb P_{\C}(\C T_pS)\setminus \mathbb P_\R(T_pS)$.
\end{itemize}
\end{defn}From the definition, it is clear that every complex structure on \( S \) induces a bi-complex structure by taking its \( \mathbb C \)-linear extension. The converse is not true, since not all bi-complex structures \(  \mathcal J \) arise from complex structures on \( S \), that is, they do not preserve the real subbundle \( TS \). 

One can easily show 
that the space of bi-complex structures on \(S\) 
is in bijection with the product \( \mathcal C(S)\times\mathcal C(\overline S) \): this can be seen by assigning to each $\mathcal J$ the unique pair $(J_1, -J_2)\in \CSCS$ such that 
\[
V_{-i}(\mathcal J) =V_{-i} (J_1) \qquad \text{and} \qquad V_i(\mathcal J)= V_{-i}(-J_2)~.
\]
Under the correspondence just described, the bi-complex structures that are $\C$--linear extensions of complex structures on $TS$ are mapped into the anti-diagonal inside $\mathcal{C}(S)\times\mathcal C(\overline S)$, namely the subset $\{(J,-J) \ | \ J\in\mathcal{C}(S)\}$. The connection between bi-complex structures and positive complex metrics is given by the following result.

\begin{prop}[{\cite[Proposition 6.4]{BEE}}]\label{prop:unique_bicomplex_structure}
Given a positive complex metric $g$ on an oriented surface $S$ there exists a unique bi-complex structure $ \mathcal J$ compatible with $g$, namely $g( \mathcal J\cdot, \mathcal J\cdot)=g$. 
\end{prop}
\noindent If $z$ and $\bar{w}$ denote local coordinates for the complex structures $J_1$ and $-J_2$, respectively, (where $J_1$ and $-J_2$ are the complex structures associated with a bi-complex structure $\mathcal J$ as in Proposition \ref{prop:unique_bicomplex_structure}) then the positive metric $g$ can be written as
$$
g=\lambda\,\mathrm{d}z\mathrm{d}\overline w~,
$$
where $\lambda$ is a local nowhere-vanishing smooth function with complex values. Moreover, the two tangent vectors $\partial_{\bar{z}}$ and $\partial_{w}$ generate the two isotropic directions of $g$ and satisfy $ \mathcal J(\partial_{\bar z})=-i\partial_{\bar z}$ and $ \mathcal J(\partial_w)=i\partial_w$. 

Let us now consider $S$ \emph{closed} oriented surfaces. From the discussion above, one can define a map
\[
\left\{\text{positive complex metrics}\right\}/C^\infty(S,\C^*)\rightarrow \mathcal C(S)\times \mathcal C(\overline S)~,
\]
that associates to any representative in a conformal class of positive complex metrics the two complex structures we introduced above.
This map is in fact a bijection, and an inverse is constructed through \emph{Bers metrics}, which we define in the next section.


\subsection{Bers' simultaneous uniformization and Bers metrics}\label{sec: Bers metrics}

The purpose of this section is to introduce, when $S$ is a closed surface, a canonical way of assigning a positive complex metric to each pair $\ccpair\in \CSCS$. This canonical positive complex metric will be called \emph{the Bers metric} of $(J_1,-J_2)$ and will be denoted by $\hpair$. The essential result is the following:

\begin{prop}[{\cite[Theorem B]{ESholodependence}}]
    \label{prop: Bers metrics holo dependence}
    Let $S$ be a closed oriented surface. The map 
    \begin{align*}
        \mathcal C(S)\times\mathcal C(\overline S) &\to \mathrm{Sym}_2(\C TS, \C)\\
        \ccpair &\mapsto  \hpair~,
    \end{align*}
    has the following properties:
    \begin{itemize}
        \item The tensor $\hpair$ is a positive complex metric of constant curvature $-1$.
        \item The pair of complex structures associated to $\hpair$ are precisely $J_1$ and $-J_2$.
        \item $\vb*h(J, -J)$ is the unique hyperbolic Riemannian metric whose conformal class is $J$.
    \end{itemize}
    Moreover, when considered as a map
  between complex Fréchet manifolds, it is holomorphic.
\end{prop}

Although we will not need the explicit construction of $\hpair$, we outline it here for the sake of completeness. 
We first recall that a representation \(\rho:\pi_1(S)\to \PSL(2,\mathbb{C})\) is called \emph{quasi-Fuchsian} if its limit set inside \(\partial_\infty\mathbb{H}^3 \cong \mathbb{C}\mathrm{P}^1\) is a Jordan curve. The set of all quasi-Fuchsian representations, defined up to conjugation, will be denoted by \(\mathcal{QF}(S)\). It is classically known that \(\mathcal{QF}(S)\) forms a smooth open subset of the character variety \(\chi\big(\pi_1(S),\PSL(2,\C)\big)\) (see e.g.~\cite{marden1974geometry}).
The construction of $\hpair$ relies on the well-known Bers' Simultaneous Uniformization Theorem \cite{SU}, which we briefly recall. Any quasi--Fuchsian representation $\rho\colon \pi_1(S)\to \PSL(2,\C)$ acts freely and properly discontinuously on the complement of its limit set in $\mathbb{C}\mathrm{P}^1$ and preserves its connected components $\Omega_\rho^\pm$. (We distinguish the connected components by requiring that the limit set is positively oriented with respect to $\Omega_\rho^+$ and negatively oriented for $\Omega_\rho^-$.) As such, $\rho$ comes associated with a pair $([J_1], [-J_2]) \in \mathcal{T}(S) \times \mathcal{T}(\overline{S})$ encoding the complex structures of $\Omega_\rho^+/\rho(\pi_1(S))$ and $\Omega_\rho^-/\rho(\pi_1(S))$. The Bers' Simultaneous Uniformization theorem then can be rephrased by saying that, for all $(J_1, -J_2)\in \mathcal{C}(S)\times \mathcal{C}(\overline{S})$, there exists a quasi-Fuchsian representation $\rho:\pi_1(S)\to \PSL(2,\C)$ that admits $\rho$--equivariant holomorphic diffeomorphisms
$$f_1: (\widetilde S, \tilde J_1) \to \Omega_\rho^+`,\qquad \qquad  \overline{f_2}: (\widetilde S, -\tilde{J_2}) \to \Omega_\rho^-~.
$$
Moreover, the representation $\rho$ is unique up to conjugation and, for a fixed $\rho$, the $\rho$--equivariant holomorphic diffeomorphisms $f_1$ and $f_2$ are uniquely determined.

Now, it turns out that, given an immersion of a surface in $\widetilde{\mathbb G}\cong\CP^1\times\CP^1\setminus\Delta$, if the projection on each $\CP^1$--factor is a local diffeomorphism, then the pull-back metric is a positive complex metric, and the two complex structures associated to it are precisely the complex structures pulled back from each $\CP^1$--factor. So, the pull-back of $\inners_{
\widetilde{\mathbb G}}$ via the immersion $(f_1, \overline{f_2}): \widetilde S\to \Omega^+_\rho\times \Omega^-_\rho\subset\CP^1\times\CP^1\setminus\Delta$, is the desired positive complex metric $\hpair$, which is thus canonically defined. It can be expressed as:
\begin{equation*}
\hpair:= -\frac 4 {(f_1-\overline{f_2})^2}\mathrm{d}f_1\,\mathrm{d}\overline{f_2}~.
\end{equation*}

When $J_1=J_2$, the positive complex metric $\hpair$ is the complex-linear extension of the Riemannian metric of constant curvature $-1$ in the conformal class of $J_1$. In the general case, as a result of we will discuss in Section~\ref{sec: minimal immersions}, the metric $\hpair$ has constant curvature $-1$. In fact, as proved in \cite[Theorem 6.15]{BEE}, Bers metrics form a connected component of the subset of positive complex metrics of constant curvature $-1$.

\section{Minimal immersions in \texorpdfstring{$\mathrm{PSL}(2,\mathbb{C})$}{PSL(2,C)}}
\label{sec: minimal immersions}

\subsection{Immersions of surfaces into \texorpdfstring{$\mathrm{PSL}(2,\mathbb{C})$}{PSL(2,C)}}\label{subsec:immersions_into_PSL(2,C)}

Complex metrics arise naturally from the theory of immersions inside holomorphic Riemannian manifolds. In this section, we study immersions in $\PSL(2,\C)$. More specifically, we will introduce a subclass of immersions in $\PSL(2,\C)$ that simultaneously generalizes almost-Fuchsian minimal surfaces in $\mathbb{H}^3$ and spacelike maximal surfaces in $\Ads^3$, as well as immersions in $\mathbb G$. 

While we are 
 interested in studying immersions $\sigma:\widetilde{S} \to \PSL(2,\C)$ that are equivariant with respect to a representation $\rho:\pi_1(S) \to \PSL(2,\C) \times \PSL(2,\C)$, the first definitions and properties we introduce are of a local nature. Therefore, for the moment we will work with a generic orientable surface $\Sigma$. From Section~\ref{subsec:regular min imm}, $\Sigma$ will be replaced by the universal cover $\widetilde{S}$ of a \emph{closed} orientable surface $S$.

In the following, given an immersion $\sigma:\Sigma\to\PSL(2,\C)$, with a slight abuse of notation, we will denote by $\mathrm{d}\sigma_{p}$ the $\C$--bilinear extension of the differential of $\sigma$, so that $\mathrm{d}\sigma_{p}:\C T_{p}\Sigma \rightarrow \C T_{\sigma(p)}\PSL(2,\C)$. 

\begin{defn}\label{def: admissible and positive immersion}
    Let $\sigma:\Sigma\to\PSL(2,\C)$ be an immersion. The pull-back bilinear form $\Is:=\sigma^*{\inners}$ is called the (complex) \emph{first fundamental form} of $\sigma$. Moreover, we say that:
    \begin{itemize}
        \item  $\sigma$ is \emph{admissible} if $\Is$ is a complex metric on $\C T\Sigma$.
        \item $\sigma$ is \emph{positive} if $\Is$ is a positive complex metric on $\C T\Sigma$.
    \end{itemize}
\end{defn} 

Recall that positive complex metrics were defined in Definition~\ref{def: positive}. In particular, we note that by our definition, a positive immersion is automatically admissible.

Observe that, if $\sigma$ is admissible, then $\sigma$ is an immersion, and $\pi^{1,0} \circ \mathrm{d}\sigma_p : \C T_{p}\Sigma \rightarrow T^{1,0}_{\sigma(p)}\PSL(2,\C)$ is injective. Moreover, along an admissible immersion $\sigma$ there exist exactly two \emph{local} orthogonal vector fields on which the complex quadratic form $\langle\cdot,\cdot\rangle$ is identically equal to $1$, say $N$ and $-N$, along $\sigma$. 

 As shown in \cite[Theorem 6.5]{BEE}, if $\sigma$ is a positive immersion, then there exists a global volume form $\mathrm{d}A_{\Is}$ for $\Is$ on $\Sigma$, canonically consistent with the orientation on $\Sigma$. Hence, there exists a \emph{global} equivariant choice of the normal vector field $N$ along $\sigma$ that satisfies $\langle N,N \rangle = 1$, uniquely defined so that, for all $v, w\in \C T_p\Sigma$, we have 
\begin{equation}
    \Omega(v,w, N)= \mathrm{d}A_{\Is}(v,w)~.
\end{equation}
where $\Omega$ is the holomorphic volume form introduced in Section~\ref{sec: PSL(2,C)}.

Let $D$ be the holomorphic Levi-Civita connection associated with the holomorphic metric on $\PSL(2,\C)$. Similarly to what happens in the classical theory of hypersurfaces in Riemannian geometry, one can see by metric compatibility that $\inner{\mathrm DN, N}=0$. Hence we can give the following definitions.

\begin{defn}
    Let $\sigma:\Sigma\to\PSL(2,\C)$ be a positive immersion. Then:
    \begin{itemize}
        \item The \emph{complex shape operator} of $\sigma$ is the endomorphism $B$ of $\C T\Sigma$ defined by $\mathrm{d}\sigma\circ B=-(DN)\circ \mathrm{d}\sigma$.
        \item The \emph{complex second fundamental form} of $\sigma$ is the symmetric bilinear form on $\C T\Sigma$ defined by $\II=\Is(B\cdot, \cdot)$.
    \end{itemize}
\end{defn}

The symmetry of $\II$ or, equivalently, the fact that $B$ is $\Is$--self-adjoint, follows by standard computations (see \cite[Proposition~3.5]{BEE}), where it is also shown that
$$
\sigma^*D= \nabla^\Is + \II \cdot N~,
$$
where $\nabla^\Is$ denotes the complex Levi-Civita connection of $\Is$. An analog of the Gauss and Codazzi equations hold:
\begin{align}
\label{eq: Gauss PSL(2,C)}
    \text{Gauss equation:} \qquad & \mathrm K_{\Is} =-1 +\mathrm{det}_{\Is}{\II}\\
\label{eq: Codazzi PSL(2,C)}
    \text{Codazzi equation:} \qquad & {\nabla^\Is}\II \text{ is tri-symmetric} .
\end{align}
In fact, for $\Sigma$ simply connected, given two tensors $\Is$ and $\II$ solving the Gauss and Codazzi equations, there exists an immersion into $\PSL(2,\C)$ with those tensors as immersion data, as explained in the following theorem.

\begin{thm}[{\cite[Theorem A]{BEE}}]
\label{thm: integration GC}
Let $\Sigma$ be a simply connected surface, let $\Is$ be a positive complex metric with complex Levi-Civita connection $\nabla^\Is$, and let $\II\in \mathrm{Sym}_2(\mathbb{C} T\Sigma)$. Then, there exists a positive immersion $\sigma:\Sigma\to \PSL(2,\C)$ with $\Is$ and $\II$ as first and second fundamental form if and only if they satisfy Equations~\eqref{eq: Gauss PSL(2,C)}-\eqref{eq: Codazzi PSL(2,C)}. 
Moreover, 
$\sigma$ is unique up to post-composing with elements in $\PSL(2,\C)\times \PSL(2,\C)$. 
\end{thm}

\begin{remark}\label{rmk:uniqueness post composition}
Observe that, in the notation of Theorem~\ref{thm: integration GC}, if $\sigma_1,\sigma_2\colon \Sigma\to \PSL(2,\C)$ have the same first and second fundamental forms, then there exists exactly one isometry $\phi\in \PSL(2,\C)\times\PSL(2,\C)$ such that $\sigma_2 = \phi \circ \sigma_1$. Indeed, at any given point $p\in \Sigma$ and linearly independent $v,w\in T_p\Sigma$, the map $\phi$ is such that $(\mathrm{d}\phi)_{\sigma_1(p)}((\mathrm{d}\sigma_1)_p(v))=(\mathrm{d}\sigma_2)_p(v)$, $(\mathrm{d}\phi)_{\sigma_1(p)}((\mathrm{d}\sigma_1)_p(w))=(\mathrm{d}\sigma_2)_p(w)$ and, as a consequence of $\phi^*\Omega=\Omega$, also $(\mathrm{d}\phi)_{\sigma_1(p)}(N_{\sigma_1(p)}) = N_{\sigma_2(p)}$: these conditions completely characterize $\mathrm{d}\phi$ at $\sigma_1(p)$, hence they characterizes $\phi$ globally. 
\end{remark}



Recall from Section \ref{subsubsec:tot geod} that $\mathbb G$, $\mathbb H^3$, and $\Ads^3$ all embed in $\PSL(2,\C)$ in an isometric totally geodesic way. As a result, the study of immersions of surfaces into $\PSL(2,\C)$ includes the analogue theories for $\mathbb G$, $\mathbb H^3$, and $\Ads^3$, as explained in the following proposition.

\begin{prop}[{\cite[Theorem C]{BEE}}] \label{prop:three cases immersions}
    Let $\Sigma$ be a simply connected surface, let $\sigma\colon \Sigma\to \PSL(2,\C)$ be a positive immersion with fundamental forms $\Is$ and $\II$. 
    Then, up to composing $\sigma$ with elements in $\PSL(2,\C)^2$:
    \begin{enumerate}[leftmargin=2em]
        \item $\sigma(\Sigma)\subset \mathbb G$ if and only if $\II=0$.
        \item $\sigma(\Sigma)\subset\mathbb H^3$ if and only if $\Is$ is a (real-valued) Riemannian metric and $\II$ is real-valued. In this case, $\Is$ and $\II$ coincide with the usual fundamental forms in hyperbolic geometry. 
        \item $\sigma(\Sigma)\subset\Ads^3$ if and only if $\Is$ is a (real-valued) Riemannian metric and $i\II$ is real-valued. In this case, $\Is$ and $i\II$ coincide with the usual fundamental forms in Anti-de Sitter geometry. 
    \end{enumerate}
\end{prop}

\begin{remark}
\label{rmk: af e mgh 1}
Let us briefly comment on these three items in terms of the Gauss-Codazzi equations~\eqref{eq: Gauss PSL(2,C)}-\eqref{eq: Codazzi PSL(2,C)}. 

In the case of item (1), we have $\II=0$, so by Equation~\eqref{eq: Gauss PSL(2,C)} the first fundamental form $\Is$ has constant curvature -1. We thus recover the fact, already mentioned in Section~\ref{sec: Bers metrics}, that Bers metrics have constant curvature $-1$: in fact, they are obtained by pull-back of the metric of $\widetilde{\mathbb G}$ from an immersion of $\Sigma$ into $\widetilde{\mathbb G}$, so, after composing with the projection $\widetilde{\mathbb G}\to\mathbb G$ (which is a local isometry), the first fundamental form is realized as in item (1). 

For item (2), all the objects involved in Equation~\eqref{eq: Gauss PSL(2,C)} become real, and one recovers precisely the Gauss-Codazzi equations defining immersions in $\mathbb H^3$. 


Finally, let us consider item (3). Since the pull-back metric is assumed to be complex positive (see Definition~\ref{def: positive}), the isotropic directions do not belong to the real tangent space to the surface. Therefore the first fundamental form can never have signature $(1,1)$, i.e. it is necessarily Riemannian. Hence, the unit normal vector field to $\Sigma$ inside $T \Ads^3$ is timelike, and coincides with $i \, N$, where $N$ is the unit normal in $T \PSL(2,\mathbb{C})$. This is why the second fundamental form of $\Sigma$ as an immersed surface in $\Ads^3$ is equal to $i\, \II$. As we write Equation~\eqref{eq: Gauss PSL(2,C)} in terms of $i \, \II$, we get $\mathrm K_{\Is} =-1 -\mathrm{det}_{\Is}(i \, \II)=-1+\mathrm{det}_{\Is}(\II)$, yielding the usual equations for the Anti-de Sitter case (see \cite[\S~6.1]{bonsante2020}). 

\end{remark}


\subsection{Minimal immersions}\label{subsec: minimal immersions}

Here we introduce minimal immersions in $\PSL(2,\C)$ in a way similar to what is done in the Riemannian case for $\mathbb H^3$. Remarkably, while in that setting a minimal immersion determines one holomorphic quadratic differential via its (real) second fundamental form, for $\PSL(2,\C)$  the complex second fundamental form $\II$ determines \emph{two} quadratic differentials $q_1,\overline{q_2}$. Those are holomorphic with respect to $J_1$ and $-J_2$, that is, with respect to the two complex structures on the surface induced by the unique bi-complex structure $\mathcal J$ compatible with the  complex first fundamental form $\Is$ (see Proposition~\ref{prop:unique_bicomplex_structure}).

\begin{defn}
    Let $\sigma\colon \Sigma\to \PSL(2,\C)$ be a positive  immersion. Then $\sigma$ is \emph{minimal}  if $\trace_{\Is}\II=0$.
\end{defn}

The second fundamental form depends on the choice of the  vector field $N$ such that $\langle N,N \rangle$, but only up to a sign, since the two possible such vector fields are $N$ and $-N$. So the definition of minimal immersion is well-posed. 

 Recall that, by Theorem~\ref{thm: integration GC}, positive  immersions are completely characterized by the immersion data $(\Is,\II)$ satisfying Equations~\eqref{eq: Gauss PSL(2,C)}-\eqref{eq: Codazzi PSL(2,C)}. For minimal ones, $\II$ has an interesting alternative description.

\begin{thm}
\label{thm: II sum of quadratic diffs}
    Let $\Is$ be a positive complex metric on an oriented surface $\Sigma$ and let $(J_1,-J_2)$ be the associated pair of complex structures. Then $\trace_{\Is}\II=0$ and $\nabla^\Is\II$ is tri-symmetric if and only if $\II$ can be expressed as
    \begin{equation}\label{eq: II sum hqd}
        \II=\frac 1 2 q_1+\frac 1 2\overline{q_2}~,
    \end{equation}
    for  $q_1\in \QD_{J_1}(\Sigma)$ and $\overline{q_2}\in \QD_{-J_2}(\Sigma)$.
    Moreover,  $q_1$ and $\overline{q_2}$ are unique.
\end{thm}
\begin{proof}
The proof relies on ideas from \cite{emam2024bers, rungi2024complex}. 
Let $z$ and $\overline w$ be local holomorphic coordinates for the pull-back complex structures of $J_1$ and $-J_2$ on $\Sigma$, respectively. Clearly, $\mathrm{d}z$ and $\mathrm{d}\overline w$ are linearly independent at every point of $\Sigma$, as they have different kernels, corresponding to the two isotropic directions of $\Is$. In particular, $\mathrm{d}z$ and $\mathrm{d}\overline w$ span $\C T^*\Sigma$ wherever they are simultaneously defined.

Being $\II$ a $(2,0)$--tensor it can locally be written as 
    \[
    \II= \alpha(z,\overline w)\,\mathrm{d}z^2+ \beta(z,\overline w) \,\mathrm{d}z\mathrm{d}\overline w +\gamma(z,\overline w) \,\mathrm{d}\overline w^2~,
    \]
    for smooth local functions $\alpha,\beta, \gamma$. To prove the statement, we have to show that $\tr_{\Is}\II=0$ and $\nabla^{\Is}\II$ is tri-symmetric if and only if $\alpha $ is $J_1$-holomorphic, $\gamma$ is $(-J_2)$-holomorphic, and $\beta=0$. Let $\Is=\varrho \, \mathrm{d}z\mathrm{d}\overline w$ be the positive metric written in the isotropic local frame $\{\partial_{\overline z},\partial_w\}$, where $\varrho=\varrho(z,\overline w)$ is a smooth, nowhere-vanishing function. Let us define $s^2:=\Is(\partial_{\overline z},\partial_w)=\frac{1}{2}\varrho(\partial_wz)(\partial_{\overline z}\overline w)$ and consider the $\Is$--orthonormal local frame
    $$
    \Big\{e_1:=\frac{1}{s\sqrt 2}(\partial_{\overline z}+\partial_w), e_2:=i\frac{1}{s\sqrt 2}(\partial_{\overline z}-\partial_w)\Big\}~,
    $$
    for a choice of a square root of $\varrho$, possibly in a smaller contractible open subset. Using this frame we can compute
    \begin{align*}
    \trace_{\Is}\II &=\II(e_1,e_1)+\II(e_2,e_2) \\ &=\frac{1}{2s^2}\Big(\II(\partial_{\overline z}+\partial_w,\partial_{\overline z}+\partial_w)-\II(\partial_{\overline z}-\partial_w,\partial_{\overline z}-\partial_w)\Big)\\&=\frac 2 {\Is(\partial_{\overline z}, \partial_w)} \II(\partial_{\overline z}, \partial_w) \\ &=2\varrho^{-1}\beta(z,\overline w)~,
    \end{align*}
    hence $\trace_{\Is}\II=0$ if and only if $\beta=0$. By differentiating $\II$ we get
    \begin{equation}\label{eq: nabla II try simmetric proof}
    \nabla^\Is \II= \mathrm{d}\alpha \otimes \mathrm{d}z^2 + 2\alpha\,(\nabla^\Is \mathrm{d}z)\cdot \mathrm{d}z + \mathrm{d}\gamma \otimes \mathrm{d}\overline w^2+ 2\gamma\,(\nabla^\Is \mathrm{d}\overline w)\cdot \mathrm{d}\overline w~.
    \end{equation}
   We show that $\nabla^{\Is}\mathrm{d}z$ is symmetric. By metric compatibility, we have $\Is(\nabla^{\Is} \partial_{\overline z}, \partial_{\overline z})=\frac 12 \mathrm{d}(\Is(\partial_{\overline z}, \partial_{\overline z}))=0$, hence $\nabla^{\Is} \partial_{\overline z}\in \Gamma(\mathbb{C}T^*S \otimes \mathrm{Span}_{\C}(\partial_{\overline z}))$ since $\Is$ is nondegenerate. Moreover, the identities
   \[
   (\nabla^{\Is} \mathrm{d}z)(\partial_{\overline z})= \mathrm{d}(\mathrm{d}z(\partial_{\overline z})) - \mathrm{d}z(\nabla^{\Is} \partial_{\overline z}) =0
   \]
    and
    \[
   (\nabla_{\partial_{\overline z}}^\Is \mathrm{d}z)(\partial_{z})= \partial_{\overline z}(1) - \mathrm{d}z(\nabla_{ \partial_{\overline z}}^\Is \partial_z)= -\mathrm{d}z(\nabla_{\partial_z}^\Is { \partial_{\overline z}} ) = 0
   \]
   imply that $\nabla^\Is \mathrm{d}z \in\mathrm{Span}_{\C}(\mathrm{d}z^2)$. Similarly we observe that $\nabla^\Is \mathrm{d}\overline w\in \mathrm{Span}_{\C}(\mathrm{d}\overline w^2)$. By Equation~\eqref{eq: nabla II try simmetric proof}, we directly obtain that $\nabla^\Is \II$ is tri-symmetric if and only if $\tau:= \mathrm{d}\alpha \otimes \mathrm{d}z^2 + \mathrm{d}\gamma \otimes \mathrm{d}\overline w^2$ is tri-symmetric. In order to conclude the proof, notice that if $\alpha$ and $\gamma$ are holomorphic for $J_1$ and $-J_2$ respectively, then $\tau= \partial_z\alpha\, \mathrm{d}z^3 + \partial_{\overline w} \gamma \, \mathrm{d}\overline w^3$, which is clearly tri-symmetric. Conversely, if $\tau$ is tri-symmetric, imposing $\tau( \partial_{\overline z}, \partial_z, \partial_w)= \tau(  \partial_z, \partial_{\overline z},\partial_w)$ and $\tau(\partial_{w}, \partial_{\overline w}, \partial_{\overline z})= \tau(\partial_{\overline w},\partial_{w}, \partial_{\overline z})$, we get respectively that $\partial_{\overline z} \alpha=0$ and $\partial_w \gamma=0$, hence $\alpha$ is $J_1$--holomorphic and $\gamma$ is $(-J_2)$--holomorphic.
\end{proof}

Let us now specify the three particular cases of immersions characterized in Proposition~\ref{prop:three cases immersions} to the particular context of minimal immersions. 

\begin{enumerate}[leftmargin=2em]
    \item When $\II=0$, which is the case of immersions taking values in $\mathbb G$ (up to isometries), we have $q_1=\overline{q_2}=0$, by the uniqueness part of Theorem~\ref{thm: II sum of quadratic diffs}. 
    \item When $\Is$ is a Riemannian metric and $\II$ is real, which is the case for immersions into $\mathbb H^3$ (up to isometries), we have $q_1=q_2$.
    \item When $\Is$ is a Riemannian metric and $\II$ is  purely imaginary, which is the case for immersions into $\Ads^3$ (up to isometries), we have $q_1=-q_2$.
\end{enumerate}

We conclude with a straightforward consequence of Theorem~\ref{thm: II sum of quadratic diffs} and Equation~\eqref{eq: Gauss PSL(2,C)}, that will allow us to recover minimal immersions from a 4-tuple $(J_1,-J_2,q_1,\overline{q_2})$ together with a complex-valued function that solves a certain PDE (Equation~\eqref{eq: Gauss varrho} below). 

\begin{cor}
\label{cor: GC semplificato}
    Let $h$ be a positive complex metric of curvature $-1$ on an oriented surface $\Sigma$, whose associated complex structures are $J_1$ and $-J_2$, let $\II$ be a complex $(0,2)$--symmetric tensor of the form $\II=(q_1+\overline{q_2})/2$ for $q_1\in \QD_{J_1}(\Sigma)$ and $\overline{q_2}\in \QD_{-J_2}(\Sigma)$, and let $\varrho\colon \Sigma\to \C^*$ be a smooth function. Then the pair $(\Is:=\varrho h, \II)$ satisfies the Equation~\eqref{eq: Gauss PSL(2,C)} if and only if $\varrho$ solves
        \begin{equation}
            \label{eq: Gauss varrho}
            \Delta_{h} \log(\varrho)= 2\varrho + \frac{2q_1\overline{q_2}}{\varrho h^2} -2~.
        \end{equation}

\end{cor}
In Equation~\eqref{eq: Gauss varrho},  $\log(\varrho)$ can be locally defined as an arbitrary determination of the logarithm of $\varrho$, and we recall that $\Delta_{h} \log(\varrho)$ is actually globally well-defined.
\begin{proof}
Let $z$ and $\overline w$ be local coordinates for $J_1$ and $-J_2$, respectively. Denoting $q_1=\phee \mathrm{d}z^2$, $\overline{q_2}= \overline \psi \mathrm{d}\overline w^2$, $h=\lambda \mathrm{d}z \mathrm{d}\overline w$, the tensors $\Is$ and $\II$ are given in the local frame $(\partial_w, \partial_{\overline z})$ by the matrices
\[
\Is \longleftrightarrow \begin{pmatrix}
    0 & \frac 12 \lambda \rho \partial_{\overline z}\overline w \partial_w z \\
    \frac 12 \lambda \rho \partial_{\overline z}\overline w \partial_w z & 0 
\end{pmatrix} \qquad \II \longleftrightarrow
\begin{pmatrix}
    \phee (\partial_w z)^2 & 0\\
    0 & \overline{\psi} (\partial_{\overline z} \overline w)^2
\end{pmatrix}.
\]

    Hence, we compute $${\det}_{\Is}\II= -\frac {\phee \overline \psi}{\varrho^2 \lambda^2 }=-\frac{q_1\overline{q_2}}{\varrho^2 h^2}.$$ As a result, through Equation~\eqref{eq: conf curvature}, we deduce that the Gauss equation~\eqref{eq: Gauss PSL(2,C)} is equivalent to Equation~\eqref{eq: Gauss varrho}.
\end{proof}

\subsection{Equivariant and regular minimal immersions}\label{subsec:regular min imm}


We now get back to the setting of a \emph{closed} oriented surface $S$. We will apply the results from Sections~\ref{subsec:immersions_into_PSL(2,C)} and~\ref{subsec: minimal immersions} to $\Sigma=\widetilde{S}$, the universal cover of $S$. The fundamental observations that we use are the following. First, given a positive (minimal) immersion $\sigma:\widetilde{S}\to \PSL(2,\C)$ which is \emph{equivariant} with respect to some representation of $\pi_1(S)$ to $\PSL(2,\C)\times \PSL(2,\C)$ (that we will call the \emph{holonomy} of $\sigma$), all the tensors that we constructed on $\widetilde{S}$, namely $\Is,\II,J_1,q_1,-J_2,\overline{q_2}$, descend to well-defined tensors on $S$. Second, once $J_1$ and $-J_2$ are fixed, Proposition~\ref{prop: Bers metrics holo dependence} provides us with a canonical choice of a positive complex metric on $S$ (that we denoted by $\hpair$) whose associated complex structures are $J_1$ and $-J_2$.

Combining Theorem~\ref{thm: integration GC}, Theorem~\ref{thm: II sum of quadratic diffs} and Corollary~\ref{cor: GC semplificato}, we obtain the following converse result, which will be the starting point of our construction of the deformation space. 

\begin{thm}\label{thm: GC semplification 2}
    Let $S$ be a closed oriented surface. Given:
    \begin{itemize}
        \item two complex structures $J_1$ and $J_2$ compatible with the orientation of $S$,
        \item two holomorphic quadratic differentials $q_1\in \QD_{J_1}(S)$ and $\overline{q_2}\in \QD_{-J_2}(S)$, and
        \item a smooth function $\varrho:S\to\C^ *$.
    \end{itemize} 
    there exists an equivariant positive minimal immersion $\sigma\colon\widetilde S\to \PSL(2,\C)$ with $\Is=\varrho \, \hpair$ and $\II=(q_1+\overline{q_2})/2$ as first and second fundamental form if and only if $\varrho$ solves
    \begin{equation}
            \label{eq: Gauss varrho on S}
            \Delta_{\hpair} \log(\varrho)= 2\varrho + \frac{2q_1\overline{q_2}}{\varrho \, \hpair^2} -2~.
        \end{equation}
\end{thm}
\begin{proof}
    First, a direct consequence of Theorem~\ref{thm: integration GC} and Remark~\ref{rmk:uniqueness post composition} is that equivariant immersions of $\widetilde{S}$ are uniquely determined (up to post-composition with elements of $\PSL(2,\C)\times \PSL(2,\C)$) by pairs $(\Is,\II)$ on $S$, solving Equations~\eqref{eq: Gauss PSL(2,C)} and~\eqref{eq: Codazzi PSL(2,C)}. By Theorem~\ref{thm: II sum of quadratic diffs} and Corollary~\ref{cor: GC semplificato}, this is equivalent to pairs of the form $(\Is=\varrho \, \hpair,\II=(q_1+\overline{q_2})/2)$, together with a solution of Equation~\eqref{eq: Gauss varrho} for $h=\hpair$.
\end{proof}

For simplicity, in the following, we will restrict to solutions of Equation~\eqref{eq: Gauss varrho} of the form $\varrho=e^{2u}$ with $u\colon S\to \C$. Clearly, the subset of positive metrics of the form $e^{2u}\hpair$ still determines an open subset of $\mathrm{Sym}_2(\C TS)$ that contains Riemannian metrics (i.e. for $J_1=J_2$ and $u$ real-valued). 

This motivates the following definition.

\begin{defn}
    A positive  minimal immersion $\sigma\colon\widetilde S\to \PSL(2,\C)$ is called a \emph{regular minimal immersion} if its complex first fundamental form $\Is$ is of the form $\Is=e^{2u}\hpair$.
\end{defn}

For $\varrho=e^{2u}$, Equation~\eqref{eq: Gauss varrho} is equivalent to 
\begin{equation}
\label{eq: Gauss semplif in u}
 \Delta_{\hpair} u= e^{2u} + {e^{-2u}} \frac{q_1\overline{q_2}}{\hpair^2} -1~.
\end{equation}
As a consequence of this equation and Theorem~\ref{thm: GC semplification 2}, we introduce the following notation.

\begin{defn}
    Given a regular minimal immersion $\sigma\colon\widetilde S\to \PSL(2,\C)$, we denote by $\gamma_\sigma$ the 5-tuple $(J_1,q_1, -J_2, \overline{q_2}, u)$ such that $\Is=e^{2u}\hpair$, $\II=(q_1+\overline{q_2})/2$ and $u$ solves Equation~\eqref{eq: Gauss semplif in u}.
\end{defn}

Theorem~\ref{thm: GC semplification 2} tells us that $\sigma$ is entirely determined by the 5-tuple $\gamma_\sigma$ up to post-composition in $\PSL(2,\C)\times \PSL(2,\C)$.
So, recalling that $\mathrm{QD}(S)$ denotes the bundle of holomorphic quadratic differentials over $\mathcal C(S)$, with its natural complex Fréchet manifold structure (see Section~\ref{sec: beltrami differentials}), we define the functional
\begin{gather*}
    \mathscr G\colon \mathrm{QD}(S)\times \mathrm{QD}(\overline S) \times C^{\infty}(S,\C)\to C^{\infty}(S,\C)~, \\
    \mathscr G(J_1,q_1,-J_2, \overline{q_2}, u)=\Delta_{\hpair} u - {e^{-2u}} \frac{q_1\overline{q_2}}{\hpair^2} - e^{2u} +1~.
\end{gather*}
The map $\mathscr G$ is holomorphic in the complex Fréchet sense (see \cite[\S~6.1]{ESrank2} and \cite[Proposition 6.8]{ESholodependence}). By Corollary~\ref{cor: GC semplificato}, regular minimal immersions are therefore parametrized by the space
\[
\mathrm{M}(S)
=
\left\{
(J_1,q_1,-J_2,\overline{q_2},u)
\;\middle|\;
\begin{gathered}
\text{$q_1\in\QD_{J_1}(S)$, $\overline{q_2}\in\QD_{-J_2}(S)$, and $u\in C^\infty(S,\mathbb C)$} \\
\text{satisfy $\mathscr G(J_1,q_1,-J_2,\overline{q_2},u)=0$}
\end{gathered}
\right\}~.
\]

For a given $\gamma=(J_1,q_1,-J_2, \overline{q_2}, u)\in \mathrm M(S)$, let $\partial_u\mathscr G\colon C^{\infty}(S,\C)\to C^{\infty}(S,\C)$ be the linearization of $\mathscr G$ at $\sigma$ in the last entry, namely 
\begin{equation}
\label{eq: linearized Gauss}
\partial_u \mathscr G (v)= \Delta_{\hpair} v +  \frac{2q_1\overline{q_2}}{e^{2u}\hpair^2}v  -2e^{2u}v~.
\end{equation}
The operator $\partial_u\mathscr G|_{\gamma}$ is Fredholm of index zero for all $\gamma\in \mathrm M(S)$ (see \cite[Proposition 6.8]{ESrank2}), hence it is invertible if and only if it is injective or surjective. 

\begin{defn}\label{def:regular immersion}
    A regular minimal immersion $\sigma$ is called \emph{infinitesimally rigid} if $\partial_u\mathscr G$ is bijective at $\gamma_{\sigma}$.
\end{defn} 

 For our configuration space of minimal immersions, we will need to require an additional condition on their holonomy. We fix some notation concerning the representation spaces. Let $\Hom^{irr}(\pi_1(S), \PSL(2,\C))$ be the open dense subset of homomorphisms $\pi_1(S)\to \PSL(2,\C)$ that are \emph{irreducible}, namely the corresponding action on $\C^2$ does not preserve a complex line. In the identification $\Hom(\pi_1(S), \PSL(2,\C)^2)\cong (\Hom(\pi_1(S), \PSL(2,\C)))^2$, we also denote the space
 \[
 \Hom^*(\pi_1(S), \PSL(2,\C)^2):=  (\Hom^{irr}(\pi_1(S), \PSL(2,\C)))^2
 \]
 and call its elements \emph{doubly irreducible} homomorphisms. 
 
 The quotient of $\Hom^{irr}(\pi_1(S), \PSL(2,\C))$ under the conjugacy action gives an open dense subset $\chi^{irr}(\pi_1(S), \PSL(2,\C))$ (in fact, a Zariski open subset) of the entire character variety $\chi(\pi_1(S), \PSL(2,\C))$ that is contained in its smooth locus (\cite[Theorem 4.7]{cooper2015non}, \cite{goldman1984}). In the identification $\chi(\pi_1(S), \PSL(2,\C)^2)\cong \chi(\pi_1(S), \PSL(2,\C))^2$ (see \cite[Remark 6.3]{florentino2025mixed}), we define
 \begin{equation}
 \begin{split}
 \chi^*(\pi_1(S), \PSL(2,\C)^2)&:= \Hom^*(\pi_1(S), \PSL(2,\C)^2)  /(\PSL(2,\C)^2) \\
 &=\chi^{irr}(\pi_1(S), \PSL(2,\C))\times \chi^{irr}(\pi_1(S), \PSL(2,\C))
 \end{split}
 \end{equation}
 
We define
\begin{equation}\label{eq: M*(S)}
\mathrm{M}^*(S)
=
\left\{
\gamma_{\sigma}\in\mathrm{M}(S)
\;\middle|\;
\begin{gathered}
\sigma \text{ is infinitesimally rigid and} \\
\text{its holonomy is in }
\chi^*\bigl(\pi_1(S),\PSL(2,\C)^2\bigr)
\end{gathered}
\right\}~.
\end{equation}

Observe that $\mathrm M^*(S)$ is an open subset of $\mathrm M(S)$. We now describe three particular cases of regular minimal immersions: isometric immersions of Bers metrics in $\mathbb G$, almost-Fuchsian minimal immersions in $\mathbb H^3$ and maximal immersions in $\Ads^3$.

\begin{enumerate}[leftmargin=2em]
\item Equivariant minimal immersions in $\mathbb G$
 have vanishing complex second fundamental form, and, as a consequence of Equation~\eqref{eq: Gauss PSL(2,C)}, their complex first fundamental form $\Is$ has constant curvature $-1$. Such minimal immersions are regular (according to Definition~\ref{def:regular immersion}) if additionally $\Is$ is of the form $\hpair$ for some $\ccpair\in \mathcal C(S)\times\mathcal C(\overline S)$. (For the explicit construction of the immersion, see Section~\ref{sec: Bers metrics} and Remark~\ref{rmk: af e mgh 1}.) In that case, as a consequence of the construction, they correspond to the 5-tuple $(J_1,0,-J_2,0,0)$. As such, we can define the natural inclusion
 \begin{equation}\label{eq:inclusion CC in M}
    \begin{split}
        \mathcal C(S)\times\mathcal C(\overline S)&\longrightarrow \mathrm{M}(S) \\
        \ccpair &\longmapsto (J_1,0,-J_2,0,0)
    \end{split}
  \end{equation}
whose image corresponds to the space of regular totally geodesic immersions in $\PSL(2,\C)$ of constant curvature $-1$. 
 
\item An equivariant minimal immersion $\sigma:\widetilde S\to\mathbb H^3$ is called \emph{almost-Fuchsian} if its second fundamental form in $\mathbb H^3$ (that we still denote by $\II$) satisfies the condition $\|\II\|^2\leq 2$. Almost-Fuchsian immersions are uniquely determined, up to orientation-preserving isometries of $\mathbb H^3$, by a triple $(J,q,u)$, where: $J$ is a complex structure on $S$, $q$ is a holomorphic quadratic differential for $J$, and $u\in C^\infty(S,\R)$ is a solution of the equation
 \begin{equation}
\label{eq: gauss in H3}
 \Delta_{h_J} u = {e^{-2u}} \frac{|q|^2}{h_J^2} + e^{2u} -1 
\end{equation} 
that satisfies $e^{-4u}|q|^2/h_J^2<1$, where $h_J=\vb*h(J,-J)$ denotes the unique (real) hyperbolic metric in the conformal class of $J$. (As we write $\Is = e^{2 u} h_J$ and $\II = \Re(q)$, Equation~\eqref{eq: gauss in H3} is a reformulation of the Gauss equation in $\mathbb H^3$, and the condition $e^{-4u}|q|^2/h_J^2<1$ is equivalent to the bound $\|\II\|^2\leq 2$.) Thus, we define the following space:
\begin{align*}
    \mathrm{AF}(S)&=\left\{(J, q, u) \bigg|\  (J,q)\in\QD(S),u\in C^\infty(S,\R), \frac{|q|^2}{e^{4u} h_J^2} <1,  \mathscr G(J, q, -J, \overline q, u)=0 \right\}~.
\end{align*}
     We then have a natural inclusion 
    \begin{equation}\label{eq:inclusion AF in M}
    \begin{split}
        \mathrm{AF}(S)&\longrightarrow \mathrm{M}(S) \\
        (J,q,u)&\longmapsto (J, q, -J, \overline q, u)
    \end{split}
  \end{equation}
    which corresponds geometrically to composing $\sigma$ with the totally geodesic embedding of $\mathbb H^3$ in $\PSL(2,\C)$, described in the second item from Section~\ref{subsubsec:tot geod}.

\item Likewise, equivariant maximal immersions in $\Ads^3$ (now without any assumption on the second fundamental form) are uniquely determined, up to isometries of $\Ads^3$, by a triple $(J,q,u)$, where: $J$ is a complex structure on $S$, $q$ is a holomorphic quadratic differential for $J$, and $u\in C^\infty(S,\R)$ satisfies
 \begin{equation}
\label{eq: gauss in ads3}
 \Delta_{h_J} u=  - {e^{-2u}} \frac{|q|^2}{h_J^2} + e^{2u} -1~,
\end{equation}
which is equivalent to $\mathscr G(J, q, -J, -\overline q, u)=0$. (Similarly to the almost-Fuchsian setting, Equation~\eqref{eq: gauss in ads3} is equivalent to the Gauss equation in $\Ads^3$ as we write $\Is = e^{2 u} h_J$ and $i\, \II = \Re(q)$.) So, we define 
    \begin{align*}
    \mathrm{MGH}(S)&=\left\{(J, q, u) \bigg|\  (J,q)\in\QD(S),u\in C^\infty(S,\R),  \mathscr G(J, q, -J, -\overline q, u)=0 \right\}~.
\end{align*}
We have a natural inclusion 
    \begin{equation}\label{eq:inclusion MGH in M}
    \begin{split}
        \mathrm{MGH}(S)&\longrightarrow \mathrm{M}(S) \\
        (J,q,u)&\longmapsto (J, q, -J, -\overline q, u)~,
    \end{split}
  \end{equation}
which corresponds to composing $\sigma$ with the totally geodesic embedding of $\Ads^3$ in $\PSL(2,\C)$ as described in the third item from Section~\ref{subsubsec:tot geod}.
\end{enumerate}

In view of realizing the spaces of immersions in $\mathbb G$, of almost-Fuchsian minimal immersions in $\Hyp^3$, and of maximal immersions in $\Ads^3$ as subsets of the space $\mathbb M_0(S)$ (that will be defined in the next section), it would be natural at this point to expect that the images of the inclusions~\eqref{eq:inclusion CC in M},~\eqref{eq:inclusion AF in M}, and~\eqref{eq:inclusion MGH in M} are  contained in $\mathrm M^*(S)$. While for inclusions \eqref{eq:inclusion AF in M} and~\eqref{eq:inclusion MGH in M} this is not hard to  prove (see Proposition \ref{prop: AF e MGH in M*} below), for inclusion \eqref{eq:inclusion CC in M} this has been proved in \cite{emam2024bers} but only up to a closed subset of empty interior. We provide here the precise statement, which will still be sufficient for our purposes.

\begin{prop}[{\cite[Theorem~5.16]{emam2024bers}}]\label{prop: dense open subset quasi-fuchsian}
    The subset $$\mathcal B_{reg}:=\{\ccpair\in \mathcal C(S)\times\mathcal C(\overline S)\  |\ (J_1, 0, -J_2, 0, 0)\in \mathrm M^*(S)\}$$ is open and dense in $\CSCS$. Moreover, the intersection of $\mathcal B_{reg}$ with every fiber $F$ of the projection $\mathcal C(S)\times\mathcal C(\overline S)\to\TSTS$ is connected and dense in $F$.
\end{prop}

We now prove the desired inclusion property for almost-Fuchsian minimal immersions in $\Hyp^3$ and maximal immersions in $\Ads^3$.

\begin{prop}
\label{prop: AF e MGH in M*}
  The images of   $\mathrm{AF}(S)$ and $\mathrm{MGH}(S)$ under the inclusions~\eqref{eq:inclusion AF in M} and~\eqref{eq:inclusion MGH in M}, respectively, are contained inside $\mathrm M^*(S)$.
\end{prop}

\begin{proof}
We need to check that, given a 5-tuple of the form $\gamma=(J, q, -J, \pm\overline q, u)$ for $(J,q,u)$ in $\mathrm{AF}(S)$ (resp. $\mathrm{MGH}(S)$), the linearized operator $\partial_u\mathscr G$ is injective.
Observe that, at $\gamma$, $\partial_u\mathscr G$ restricts to a linear operator $L\colon C^{\infty}(S,\R)\to C^{\infty}(S,\R)$, and, by $\C$--linearity, $\partial_u\mathscr G$ is injective if and only if $L$ is injective.

By Equation~\eqref{eq: linearized Gauss}, we have
\[
L(v)=\partial_u \mathscr G (v)= \Delta_{h_J} v -2e^{2u}\left(1\pm \frac{|q|^2}{e^{4u}h_J^2} \right)  v~,
\]
for $v\in C^{\infty}(S,\R)$, where $h_J=\vb*h(J,-J)$.
The sign is negative for the case of $\mathrm{AF}(S)$ and positive for the case of $\mathrm{MGH}(S)$. In either case (recall that in the former case we assume $e^{-4u}|q|^2/h_J^2<1$), we deduce:
\[
\int_S v L(v) \, \mathrm{d}A_h \leq  \int_S  v (\Delta_{h_J} v)  \,\mathrm{d}A_h = - \int_S \|\nabla^{h_J} v\|^2 \,\mathrm{d}A_h \le 0~,
\]
with equality if and only if $v$ is identically zero. This concludes the argument for the injectivity of $L$.
\end{proof}

\subsection{The structure of \texorpdfstring{$\mathrm{M}^*(S)$}{M*(S)} and its quotient by the \texorpdfstring{$\Diff_0$}{Diff0}--action}\label{sec: the deformation spaces}

In this section, we first analyze the structure of $\mathrm M^*(S)$ as an infinite-dimensional space and then introduce the finite-dimensional deformation space on which we will construct a complex hyperK\"ahler metric.\\
In the sequel, we consider $\mathrm{QD}(S),\mathrm{QD}(\overline S), \mathbb {QD}(S)$, and $\mathbb {QD}(\overline S)$ with the complex structures defined as in Section \ref{sec: hq differentials}. 
An application of the implicit function theorem yields the following result.
\begin{prop}[{\cite[Proposition 6.1]{ESrank2}}] 
\label{prop: Frechet structure}
$\mathrm{M}^*(S)$ is a complex Fréchet submanifold of $\mathrm{QD}(S)\times \mathrm{QD}(\overline S) \times C^{\infty}(S,\C)$ and the projection $\pi_0:\mathrm{M}^*(S)\to \QD(S)\times\QD(\overline S)$ is a local biholomorphism.

\end{prop}
We will therefore denote a generic element of $T_{\gamma} ^{1,0}\mathrm M^*(S)$ as $(\dot \mu_1, \dot q_1, \dot {\overline {\mu_2}}, \dot {\overline {q_2}}, \dot u)$.

\begin{remark}
The projection map $\pi_0\colon \mathrm M^*(S)\to \mathrm{QD}(S)\times \mathrm{QD}(\overline S)$, according to Proposition~\ref{prop: Frechet structure}, is a local biholomorphism and is clearly equivariant for the (diagonal) action of $\mathrm{Diff}_0(S)$. 
Observe that this implies that the space $\mathrm M^*(S)/\mathrm{Diff}_0(S)$ is not finite dimensional. As a consequence, one can find uncountably many regular minimal surfaces in $\PSL(2,\C)$ with the same holonomy. 
Nevertheless, we will now explain how to construct (using a similar process as in \cite[\S~6.1]{ESrank2}) a meaningful finite dimensional quotient of $\mathrm M^*(S)$. 
\end{remark}

Compose the map $\pi_0$ with the double isotopy projection $\mathrm{QD}(S)\times \mathrm{QD}(\overline S)\to \mathbb{QD}(S) \times\mathbb{QD}(\overline S)$, to get a holomorphic map 
$[\pi_0]\colon\mathrm M^*(S)\to \mathbb{QD}(S) \times\mathbb{QD}(\overline S)$.

\begin{defn}\label{def:M(S)}
    We define $\mathbb{M}(S)$ as the topological quotient of $\mathrm M^*(S)$ obtained by contracting to points the connected components of the fibers of $[\pi_0]$. We will denote by $\mathbb{M}_0(S)$ the connected component of $\mathbb M(S)$ containing the elements $[(J,0,-J,0 ,0)]$, which parametrize an isomorphic copy of Teichm\"uller space inside $\mathbb M(S)$.
\end{defn}
\begin{thm} [{\cite[Proposition 6.3]{ESrank2}}]\label{thm: Kernel projection equal to Lie derivatives}
    The space $\mathbb M(S)$ has a structure of complex manifold for which the quotient map $\mathrm M^*(S)\to \mathbb M(S)$ is holomorphic, and the natural projection $\pi\colon \mathbb M(S)\to \mathbb{QD}(S) \times\mathbb{QD}(\overline S) $ is a local biholomorphism. 
\end{thm}
The connected component $\mathbb{M}_0(S)$ just introduced (of complex dimension $12g-12$) is precisely the space on which we will later define the complex hyperK\"ahler structure. 

The last goal of this section is to show that the space $\mathbb{M}_0(S)$ contains a copy of 
$\mathcal{T}(S)\times\mathcal{T}(\overline S)$ and of the spaces of almost-Fuchsian immersions in $\mathbb{H}^3$ and spacelike minimal immersions in $\Ads^3$, where the inclusions are essentially induced by~\eqref{eq:inclusion CC in M},~\eqref{eq:inclusion AF in M} and~\eqref{eq:inclusion MGH in M}. These copies will be a complex submanifold in the former case, and totally real submanifolds in the latter two.

We begin with the case of $\mathcal{T}(S)\times\mathcal{T}(\overline S)$.

\begin{prop}
    \label{prop: TSTS dentro mathbb M}
   Consider the space $\mathcal B_{reg}$ defined in Proposition~\ref{prop: dense open subset quasi-fuchsian}. The inclusion 
    \begin{equation}\label{eq: Breg dentro M*}
        \begin{split}
                \mathcal B_{reg} &\hookrightarrow \mathrm M^*(S)\\
        (J_1, -{J_2})&\mapsto (J_1, 0, -{J_2}, 0, 0)
        \end{split}
    \end{equation}
    descends to a well-defined proper holomorphic embedding 
    \begin{equation}
    \label{eq: TSTS in mathbb M}
            \TSTS\hookrightarrow \mathbb M_0(S)
    \end{equation}
    given by $[(J_1, -{J_2})]\mapsto[(J_1, 0, -{J_2}, 0, 0)]$. 
\end{prop}
\begin{proof}
    Fix $F$ a fiber of the projection $\CSCS\to \TSTS$. By Proposition~\ref{prop: dense open subset quasi-fuchsian}, the intersection $F\cap \mathcal B_{reg}$ is non-empty and connected, hence 
    \[
    \{  (J_1, 0, -{J_2}, 0, 0) \in \mathrm M^*(S) \ |\ (J_1, -J_2)\in F\cap \mathcal B_{reg}\} 
    \]
    is connected. By the construction of $\mathbb M(S)$, the image of this subset via the projection $\mathrm M^*(S)\to \mathbb M(S)$ is just the point $\{ [ (J_1, 0, -{J_2}, 0, 0) ]\} $, which is therefore the same for all $(J_1, -J_2)\in F\cap \mathcal B_{reg}$. This proves that the map~\eqref{eq: TSTS in mathbb M} is well-defined.

To see that it is holomorphic, observe that the map in Equation~\eqref{eq: Breg dentro M*} is trivially holomorphic, hence Equation~\eqref{eq: TSTS in mathbb M} is holomorphic as it can be locally seen as the composition of holomorphic local inverses of $\mathcal B_{reg}\to \TSTS$ with the map~\eqref{eq: Breg dentro M*}, and the projection $\mathrm M^*(S)\to \mathbb M(S)$.
    
    Finally, the map~\eqref{eq: TSTS in mathbb M} is an inverse of the submersion $\mathbb M(S)\to \TSTS$ on the first and third factor, hence it is a proper embedding.
\end{proof}

Let us now explain how the isomorphism classes of points $[\gamma]\in\mathbb{M}_0(S)$ actually generalize those of almost-Fuchsian immersions in $\mathbb{H}^3$ and spacelike minimal immersions in $\Ads^3$, whose deformation spaces are respectively given by
\begin{equation}\label{eq: definition of almost fuchsian and MGH after quotient by DIff}
\mathcal{AF}(S):=\mathrm{AF}(S)/{\mathrm{Diff}_0(S)}~, \quad \text{and} \quad \mathcal{MGH}(S):=\mathrm{MGH}(S)/{\mathrm{Diff}_0(S)}~.
\end{equation}
Both of these deformation spaces carry a smooth manifold structure of real dimension $12g-12$ (for $\mathcal{AF}(S)$, this follows from the original work of Uhlenbeck \cite{uhlenbeck1983}, see also \cite{taubes,huang-lucia-tarantello2023}; for $\mathcal{MGH}(S)$, this was first proved in the work of Mess \cite{mess07}, see also \cite[Section 5]{bonsante2020}). From the inclusions of Equations~\eqref{eq:inclusion AF in M} and~\eqref{eq:inclusion MGH in M}, we obtain the following result.
\begin{prop}
\label{prop: AF and MGH are totally real submanifolds}
  The inclusions $\mathrm{AF}(S)\subset \mathrm M^*(S)$ and $\mathrm{MGH}(S)\subset \mathrm M^*(S)$ descend to totally real full dimensional embeddings of $\mathcal{AF}(S)$ and $\mathcal{MGH}(S)$ into $\mathbb M_0(S)$.
\end{prop}

\begin{proof}
    We show it for $\mathcal{AF}(S)$, and the proof for $\mathcal{MGH}(S)$ follows in the same fashion. As we discussed in Proposition~\ref{prop: AF e MGH in M*}, the space of almost-Fuchsian minimal immersions corresponds to the subset of $\mathrm M^*(S)$ given by $\mathrm{AF}(S)$. On this subspace, the fibers of the projection $[\pi_0]$ to $\mathbb{QD}(S)\times\mathbb{QD}(\overline S)$ correspond to the orbits of the (diagonal) action of $\mathrm{Diff}_0(S)$, hence they are connected. As a result, this inclusion descends to a continuous injective map $\mathcal{AF}(S)\hookrightarrow\mathbb M_0(S)$. To check that it is a smooth totally real immersion (see Definition~\ref{def: totally real submanifold}), it is sufficient to observe that, on local charts on which $\pi\colon \mathbb M_0(S)\to \mathbb{QD}(S)\times\mathbb{QD}(\overline S)$ restricts to a biholomorphism, the inclusion of $\mathcal{AF}(S)$ into  $\mathbb M_0(S)$ corresponds to the inclusion of $\mathbb{QD}(S)$ given by $[J,q]\mapsto ([J,q],[-J, \overline q])$, whose image is easily checked to be totally real. Finally, to see that this inclusion is an embedding, observe that the map $[J,q]\mapsto [(J,q,-J, \overline q, u)]$ allows to see $\mathcal{AF}(S)$ as the image of a section, defined on an open subset of $\mathbb{QD}(S)$, of the submersion 
    $$
    \mathbb M_0(S)\xrightarrow{\pi} \mathbb{QD}(S)\times\mathbb{QD}(\overline S)\xrightarrow{pr_1} \mathbb{QD}(S)~,
    $$
    where $pr_1$ is the projection on the first factor: as a consequence, $\mathcal{AF}(S)$ is an embedded submanifold of $\mathbb M_0(S)$.
\end{proof}

\begin{remark}\label{rmk: restriction of pi on AF and MGH}
Consider the natural projections $\pi_{\mathcal{AF}}\colon \mathcal {AF}(S)\to \mathbb{QD}(S)$ and $\pi_{\mathcal{MGH}}\colon \mathcal{MGH}(S)\to \mathbb{QD}(S)$ given by $[J, q,u]\mapsto [J,q]$. The first one is an embedding onto a \emph{proper} open subset of $\mathbb{QD}(S)$ (see again \cite{uhlenbeck1983}, or the more recent work \cite{bronstein-smith2024} which also provides finer information on the image) while the second one is a diffeomorphism (\cite[Theorem 3.8]{kirill07}). The restrictions of $\pi$ to $\mathcal{AF}(S)$ and $\mathcal{MGH}(S)$ essentially coincides with $\pi_{\mathcal{AF}}$ and $\pi_{\mathcal{MGH}}$, in the sense that \begin{equation}
    \begin{split}
        \pi([J,q, -J ,\overline q, u]) &= ([J,q], [-J, \overline q])= \Big(\pi_{\mathcal{AF}} ([J,q,u]),\overline{\pi_{\mathcal{AF}} ([J,q,u])}
        \Big)~, \\
         \pi([J,q, -J ,-\overline q, u]) &= ([J,q], [-J, -\overline q])= \Big(\pi_{\mathcal{MGH}} ([J,q,u]),\overline{\pi_{\mathcal{MGH}} ([J,-q,u])}\Big)~. 
\end{split}
   \end{equation}
\end{remark}

\section{The holonomy map}\label{sec:holonomy}
The aim of this section is to prove Theorem~\ref{thm: local biholos intro}. Along the way, we will introduce some technical tools, namely \emph{convenient $C^{\infty}$ maps} (see Section~\ref{subsec:holo holo}) and the \emph{complex Lie derivative} (see Section~\ref{subsec: complex lie derivatives}), which will allow us to carry out an infinitesimal analysis of the holonomy map. We will prove Theorem~\ref{thm: local biholos intro} in three steps: first, by proving that the holonomy map is a holomorphic map defined on the complex Fr\'echet space \(\mathrm{M}^*(S)\) (see Theorem~\ref{thm: hol su M^* is holomorphic}), then by showing that it descends to a holomorphic map defined on its finite-dimensional quotient \(\mathbb{M}(S)\) (see Theorem~\ref{thm: hol from mathbb M is holo}), and finally by showing that such map is a local biholomorphism (see Theorem~\ref{thm: hol local biholo}).

\subsection{The holonomy map is holomorphic}\label{subsec:holo holo}

In this section, we prove that the holonomy map ${hol}:\mathrm M^*(S)\to\chi^*\big(\pi_1(S),\PSL(2,\C)^2\big)$ is holomorphic (see Theorem~\ref{thm: hol su M^* is holomorphic}). To discuss the regularity of this map, one should first recall what differentiability means for maps defined on a Fréchet space. This subject is generally subtle, but various equivalent characterizations of the notion of Fréchet $C^\infty$ have been described (see e.g. \cite{keller1974differential}). For our purposes, we will focus on two notions, which we now recall. 

Given two Fréchet spaces $E$, $F$ and an open subset $U \subset E$, we say that $f\colon U\to F$ is \emph{Gâteaux $C^1$} (or simply $C^1$) if it admits all directional derivatives at any point, and the induced differential map $\mathrm{d}f\colon U\times E\to F$ (which, a priori, is not linear with respect to $E$) is continuous as a map between open subsets of Fréchet spaces. By iterating this notion, we say that $f$ is smooth (or Gâteaux $C^{\infty}$, or simply $C^{\infty}$) if the iterated differential maps $\mathrm{d}^k f$ are Gâteaux $C^1$ for all $k$. (Observe that, for maps between Banach spaces, $\mathrm{d}^k f$ being $C^1$ implies that $\mathrm{d}^{k-1} f$ is Fréchet--differentiable, so the notions of Fréchet $C^\infty$ and Gâteaux $C^\infty$ are equivalent in such setting). 

For a distinguished class of locally convex vector spaces that includes Fréchet spaces, called \emph{convenient} vector spaces, there exists a notion of \emph{convenient $C^{\infty}$--maps}. We say that $f$ is convenient $C^{\infty}$ if it maps $C^{\infty}$--curves to $C^{\infty}$--curves (see e.g. \cite{kriegl1997convenient} for a detailed exposition on the subject). Clearly, if $f$ is smooth then $f$ is convenient $C^{\infty}$. Surprisingly, the converse holds too.

\begin{thm}[Generalized Boman's Theorem]
\label{thm: boman}
    A map between Fréchet spaces is smooth if and only if it maps $C^{\infty}$--curves to $C^{\infty}$--curves.
\end{thm}

This statement is usually referred to as the ``Generalized Boman's Theorem'', as it extends the original work of Boman~\cite{boman1967differentiability} for finite-dimensional vector spaces. For a proof of Theorem~\ref{thm: boman}, we refer to Egeileh-Wurzbacher~\cite[Proposition~3.13]{egeileh2017infinite}, or Neeb~\cite{neeb2005infinite}. (Notice that the $c_{\infty}$--topology considered in the references is equivalent to the standard topology on $E$ by Kriegl-Michor~\cite[Theorem~4.11]{kriegl1997convenient}.)

Observe that, by restricting to local charts, Theorem~\ref{thm: boman} applies also to maps between real Fréchet--manifolds, where a map between real Fréchet--manifolds is said to be smooth if it is smooth in local charts or, equivalently, if it maps smooth curves to smooth curves. 

Finally, whenever $E$ and $F$ are complex Fréchet--vector spaces, we say that a smooth function $f\colon U \subset E \to F$ is holomorphic if $\mathrm{d}f_u \colon E\to F$ is $\C$--linear for all $u\in U$.

To show that $\mathrm{hol}$ is holomorphic, we will first prove that it is smooth, and then that it is holomorphic by applying the following result.

\begin{prop}[{\cite[Theorem D]{BEE}}]
\label{prop: holo smooth dependence}
    Let $\Lambda$ be a smooth (resp. complex) manifold, and let  ${(\Is_{\lambda}, \II_{\lambda})}_{\lambda \in \Lambda}$ be a family of equivariant immersion data of $\widetilde S$ into $\PSL(2,\C)$ whose holonomy representations are in $\chi^{*}\big(\pi_1(S),\PSL(2,\C)^2\big)$.
    
    Assume that, at each point $p\in S$, the maps $\Is_{\cdot}|_p,\II_{\cdot}|_p\colon \Lambda\to \mathrm{Sym}_2(\C T_p S)$ are smooth (resp. holomorphic) in $\lambda \in \Lambda$. Then there exists a map
   \[
        \sigma: \Lambda \times \widetilde S\longrightarrow \PSL(2,\C)
    \]
    such that $\sigma(\lambda, \cdot) \colon \widetilde S\to \PSL(2,\C)$ is an equivariant immersion with immersion data $(\Is_{\lambda}, \II_{\lambda})$, and such that $\sigma(\cdot, p)\colon \Lambda\to \PSL(2,\C)$ is smooth (resp. holomorphic) for every $p\in \widetilde S$.
   
    If the metrics $\Is_{\lambda}$ are positive, and the holonomy map
    \begin{align*} &\Lambda\longrightarrow\chi\big(\pi_1(S),\PSL(2,\C)^2\big) \\ &\lambda \longmapsto {hol}(\sigma_\lambda)
    \end{align*}
has image in $\chi^*\big(\pi_1(S),\PSL(2,\C)^2\big)$, then this map is smooth (resp. holomorphic).
\end{prop}

\begin{proof}
    The original result from Bonsante-El Emam \cite{BEE} is stated only in the holomorphic category. However, the proof argument from \cite{BEE} can also be applied to deduce the result in the smooth setting. Indeed, if we assume that the family ${(\Is_{\lambda}, \II_{\lambda})}_{\lambda \in \Lambda}$ depends smoothly on $\lambda$, Proposition~5.5 in \cite{BEE} implies that there exists a smooth map $\Phi\colon \Lambda\times \widetilde S\to \PSL(2,\mathbb{C})^2$ such that the evaluation map of $\Phi$ at any point of $\PSL(2,\C)$ gives a smooth map $\sigma\colon \Lambda\times \widetilde S\to \PSL(2,\C)$ as in Proposition~\ref{prop: holo smooth dependence}. 
    To see the part of the statement involving the holonomy, observe first of all that $\sigma_{\lambda}$ is invariant for at most one homomorphism $\rho_{\lambda}\colon \pi_1(S)\to \PSL(2,\C)^2$ by Remark \ref{rmk:uniqueness post composition}. More precisely, for all $\gamma\in \pi_1(S)$, $\rho_{\lambda}$ is characterized by the conditions $\sigma_{\lambda}\circ \gamma=  \rho_{\lambda}(\gamma)\circ \sigma_{\lambda}$ and $N_{\sigma_{\lambda}\circ \lambda} = d( \rho_{\lambda}(\gamma)) (N_{\sigma_{\lambda}})$, where $N_{\sigma_{\lambda}}$ is the unit normal of $\sigma_{\lambda}$. These conditions allow to see that the map $\lambda\mapsto \rho_{\lambda}(\gamma)$ is smooth and compute its derivatives explicitly.
\end{proof}

\begin{remark}
\label{rmk: pointwise holo and global holo}
The pointwise smoothness/holomorphicity condition assumed in Proposition~\ref{prop: holo smooth dependence} is equivalent to requiring that the maps $\Lambda\to \Gamma(\mathrm{Sym}_2(\C TS))$ given by $\lambda \mapsto \Is_\lambda$ and $\lambda \mapsto \II_{\lambda}$ are smooth/holomorphic as maps between Fréchet manifolds. Indeed, as shown in \cite[Theorem~I.3.12]{kriegl1997convenient}, 
we have the identification 
\begin{align*}
C^{\infty}(\Lambda, C^{\infty}(S,  \mathrm{Sym}_2(\C TS)))& \cong C^{\infty}(\Lambda\times S, \mathrm{Sym}_2(\C TS))\\
\left(\lambda\mapsto \alpha_{\lambda}\right) &\mapsto   \left((\lambda,p)\mapsto \alpha_{\lambda}(p)\right)~.
\end{align*}
Also, the holomorphicity conditions are trivially equivalent by the identification between the differential maps of the two formats.
\end{remark}

We adapt Proposition~\ref{prop: holo smooth dependence} to our setting. We first show that passing from elements in $\mathrm M^*(S)$ to the data $(\Is, \II)$ gives a holomorphic map.
\begin{lem}
    \label{lem: from M* to I II}
    The map
    \begin{equation}
    \label{eq: from M* to I II}
    \begin{split}
        \mathrm M^*(S)&\to \Gamma(\mathrm{Sym}_2(\C TS))^2\\
        \gamma= (J_1, q_1, -J_2, \overline{q_2}, u) &\mapsto (\Is_{\gamma}, \II_{\gamma}):= (e^{2u}\hpair, q_1+\overline{q_2})~,
    \end{split}
    \end{equation}
    is holomorphic.
\end{lem}
\begin{proof}
     By Proposition~\ref{prop: Bers metrics holo dependence}, the map $\ccpair \mapsto \hpair$ is holomorphic. Therefore, the function $\gamma\mapsto (\Is_{\gamma}, \II_{\gamma})$ is holomorphic as it is the restriction of the following holomorphic map:
     \begin{gather*}
     \mathcal C(S)\times \Gamma(\mathrm{Sym}_2(\C TS))\times \mathcal C(\overline S) \times \Gamma(\mathrm{Sym}_2(\C TS))\times C^{\infty}(S,\mathbb C)\to (\Gamma(\mathrm{Sym}_2(\C TS))^2 , \\
     (J_1, \alpha, -J_2, \beta, u) \mapsto (e^{2u}\hpair, \alpha+\beta) .
     \end{gather*}
\end{proof}



We now deduce the following adaptation of Proposition~\ref{prop: holo smooth dependence} to the setting of this section, which is the essential ingredient to prove Theorem~\ref{thm: hol su M^* is holomorphic}.

\begin{lem}
\label{lem: holo dependence}
        Let $\Lambda$ be a smooth (resp. complex) manifold and let $F:\Lambda\to \mathrm M^*(S)$ be a smooth (resp. holomorphic) map. Then the composition $\mathrm{hol}\circ F:\Lambda \to  \chi^{*}\big(\pi_1(S),\PSL(2,\C)^2\big)$ is smooth (resp. holomorphic).
\end{lem}
\begin{proof}
The composition of $F$ with the map from Equation~\eqref{eq: from M* to I II} is a smooth/holomorphic map $\Lambda\to (\Gamma(\mathrm{Sym}_2(\C TS)))^2$, by Lemma~\ref{lem: from M* to I II}. By Remark~\ref{rmk: pointwise holo and global holo}, the composed map is smooth/holomorphic pointwise, hence it gives a family $\lambda\mapsto (\Is_{\lambda}, \II_{\lambda})$ as in the hypothesis of Proposition~\ref{prop: holo smooth dependence}. The lemma then follows by Proposition~\ref{prop: holo smooth dependence}.
\end{proof}

We now have all the ingredients to prove the following.

\begin{thm}
 \label{thm: hol su M^* is holomorphic}
 The holonomy map ${hol}\colon \mathrm M^*(S)\to \chi^*(\pi_1(S),\PSL(2,\C)^2)$ is holomorphic.
\end{thm}
\begin{proof}
    By Lemma~\ref{lem: holo dependence}, $\mathrm{hol}$ maps $C^{\infty}$--curves into $C^{\infty}$--curves and, hence, it is smooth by Theorem~\ref{thm: boman}. 
    
    To see that the map is holomorphic, we fix $\sigma\in \mathrm M^*(S)$ and $\dot \sigma \in T_{\sigma}\mathrm M^*(S)$. Then, by taking a local chart, we can find a (small) complex submanifold $L$ of $\mathrm M^*(S)$ of complex dimension $1$ passing through $\sigma$ and such that $T_{\sigma}L=\mathrm{Span}_{\C}(\dot \sigma)$. By Lemma~\ref{lem: holo dependence}, the restriction of $\mathrm{hol}$ to $L$ is holomorphic and, hence, $\mathrm{d}(\mathrm{hol})_{\sigma}(i \, \dot \sigma) = i \, \mathrm{d}(\mathrm{hol})_{\sigma}(\dot \sigma)$. We conclude that $\mathrm{d}(\mathrm{hol})_{\sigma}$ is $\C$--linear for all $\sigma\in \mathrm M^*(S)$. Thus, $\mathrm{hol}$ is holomorphic. 
\end{proof}

\subsection{The complex Lie derivative and the holonomy map on \texorpdfstring{$\mathbb{M}(S)$}{M(S)}} \label{subsec: complex lie derivatives}
 
We shall now explain how the holonomy map studied above descends to a well-defined holomorphic map $\mathrm{hol}: \mathbb{M}(S) \rightarrow \chi^*(\pi_1(S),\PSL(2,\C)^2)$. This will follow from a general principle (Lemma~\ref{lem: invariance of map under Diff-action}) which makes use of complex Lie derivatives originally defined in \cite[\S~4.1]{emam2024bers}, which we now recall.


Let $\upalpha$ be a (complex) tensor on $S$ and let $X\in \Gamma(\C TS)$. We define the \emph{complex Lie derivative} of $\upalpha$ in the direction of $X$ as the complex tensor
\begin{equation}\label{eq: complex_lie_derivative}
\mathscr L_X \upalpha= \mathscr L_{\Re(X)}\upalpha+ i\,\mathscr L_{\Im(X)}\upalpha~.
\end{equation} 
One motivation for this definition is that, for all holomorphic $d$--differentials
$\upalpha\in H^0(S,\mathcal K^d_J)$ and $X\in \Gamma( \C TS)$, we have $\mathscr L_X \upalpha= \mathscr L_{X_J}\upalpha= \mathscr L_{X_J+ \overline{X_J}}\upalpha$. We explain in more details how to compute the complex Lie derivative of a complex structure. Let $J \in \Gamma(\mathrm{End}(TS))$ be a complex structure on $S$ and extend it by $\C$--linearity to an endomorphism of $\C TS$. With a slight abuse of notation, we will still denote by $J$ this extension. Then, the complex Lie derivative of $J$ is
\[
    (\mathscr{L}_{X}J)Y:= [X,JY]-J[X,Y] 
\]
for all complex vector fields $X,Y \in \Gamma(\C TS)$. Here the Lie bracket is $\C$--bilinear, hence Equation~\eqref{eq: complex_lie_derivative} holds. In this way, $\mathscr{L}_{X}J$ is naturally an element in $\C T_{J}\mathcal{C}(S)$. As explained in Section~\ref{sec: beltrami differentials}, the space $\mathcal{C}(S)$ has a natural complex structure and we can thus split $\C T_{J}\mathcal{C}(S)$ into its holomorphic and anti-holomorphic part. More concretely, a complex tangent vector $\dot{J}$ is an endomorphism of $\C TS$ that anticommutes with $J$. This implies that $\dot{J}(T^{1,0}S)=T^{0,1}S$ and vice versa. Therefore,
 \[
     \C T_{J}\mathcal{C}(S) =  \Omega^{0,1}(S, T^{1,0}S) \oplus \Omega^{1,0}(S, T^{0,1}S)~,
\]
and this is exactly the eigenspace decomposition with respect to the complex structure defined in Section $\ref{sec: beltrami differentials}$ with the first factor being the holomorphic tangent bundle: indeed, $\Omega^{0,1}(S, T^{1,0}S)$ identifies with the space of Beltrami differentials of $(S,J)$.
In particular, the complex Lie derivative $\mathscr{L}_{X}J$ can be written as sum of a holomorphic and anti-holomorphic tangent vector. A simple computation shows that the $(1,0)$-part can be written as
\[
    (\mathscr{L}^{1,0}_{X}J)Y=((\mathscr{L}_{X^{1,0}+\overline{X^{1,0}}}J)Y^{1,0})^{1,0} \ ,
\]
where, for any $Z\in \Gamma(\C TS)$, we have set $Z^{1,0}=\frac{1}{2}(Z-iJZ)$.

\begin{remark} The complex Lie derivative of a complex structure was originally defined in \cite{emam2024bers} in a slightly different way. The discussion above shows that the definition in \cite{emam2024bers} of complex Lie derivative of $J$ is equal to the $(1,0)$--part of the one considered here. In particular, since in this paper we will mostly work with the holomorphic tangent bundle of the space of complex structures (and of holomorphic quadratic differentials) the two notions actually coincide and we will remove the superscript $(1,0)$ from the notation.
\end{remark}

The complex Lie derivative allows to give a nice description of the vectors that are tangent to the fibers of the projection $\mathrm{M}^*(S)\to \mathbb M(S)$, as explained in the following Proposition.

 \begin{prop}[{\cite[Proposition~6.6]{ESrank2}}]
     \label{prop: tangent fiber M*}
     For all $\gamma\in \mathrm M^*(S)$, an element $X\in T^{1,0}_{\gamma}\mathrm M^*(S)$ lies in the kernel of the differential of $\mathrm M^*(S)\to \mathbb M(S)$ at $\gamma$ if and only if it is of the form 
     $\mathscr L_X\gamma$ for some $X\in \Gamma(\C TS)$, where $\mathscr L_X\gamma$ denotes the complex Lie derivative $\mathscr L_X$ on each entry of the 5-tuple $\gamma$.
 \end{prop}

As explained in \cite{ESrank2}, this model for the tangent space to the fibers gives a nice and simple criterion to establish when a holomorphic map from $\mathrm M^*(S)$ descends to a map from $\mathbb M(S)$. This is explained in the  following Lemma.

\begin{lem}[{\cite[Proposition 6.7]{ESrank2}}]\label{lem: invariance of map under Diff-action}
Let $\mathbb X$ be a complex manifold and let $\hat F:\mathrm M^*(S)\to\mathbb X$ be a holomorphic map. Then $\hat F$ descends to a well-defined holomorphic map $F:\mathbb M(S)\to\mathbb X$ if and only if $\hat F$ is invariant under the action of $\operatorname{Diff}_0(S)$ on $\mathrm M^*(S)$.  Equivalently, $F$ descends to $\hat F$ if and only if $\mathrm d\hat F\big(\mathscr L_{\mathsf X}\gamma)=0$ for all $\mathsf X\in \Gamma(TS)$ and $\gamma\in\mathrm{M}^*(S)$.

\end{lem}
By Lemma~\ref{lem: invariance of map under Diff-action}, we can immediately deduce that the holonomy map from $\mathrm M^*(S)$ descends to $\mathbb M(S)$.

\begin{thm}
\label{thm: hol from mathbb M is holo}
    The holonomy map ${hol}: \mathrm M^*(S)\to \chi^*\big(\pi_1(S),\PSL(2,\C)^2\big)$ descends to a well-defined holomorphic function $\mathrm{hol}: \mathbb M(S)\to \chi^*\big(\pi_1(S),\PSL(2,\C)^2\big)$.
\end{thm}
\begin{proof}
    The map $hol$ is holomorphic by Theorem~\ref{thm: hol su M^* is holomorphic}, and it is clearly $\mathrm{Diff}_0(S)$-invariant. The statement now follows from Lemma~\ref{lem: invariance of map under Diff-action}.
\end{proof}

In short, we have the following commutative diagram 

\begin{center}
\begin{equation}
\label{eq: diagramma proiezioni}
   \begin{tikzcd}
\mathrm{QD}(S)\times \mathrm{QD}(\overline S)  \arrow[d] & \mathrm M^*(S) \arrow[l, "\pi_0"] \arrow[d] \arrow[dr, "{hol}" ] \\
\mathbb{QD}(S)\times\mathbb{QD}(\overline S) & \mathbb M(S) \arrow[l, "\pi"] \arrow[r, "\mathrm{hol}"'] & \chi^{*}\big(\pi_1(S), \PSL(2,\C)^2\big)
\end{tikzcd}
\end{equation}
\end{center}

with $\pi$ and $\pi_0$ local biholomorphisms.

\begin{remark}
    Most of the statements in this section also appear in \cite{ESrank2} in the setting of complex harmonic $ G$--bundles for $ G=\PSL(2,\R)\times \PSL(2,\R)$, corresponding to complex harmonic maps into $\mathbb G\times \mathbb G$. In fact, the two approaches produce the same space $\mathrm M^*(S)$ of immersion data, hence the same quotient $\mathbb M(S)$, despite they parametrize immersions in different spaces. Since the data have a different geometric meaning, the holonomy maps from $\mathrm M^*(S)$ and $\mathbb M(S)$ defined in the two papers could be different a priori, therefore we had to prove independently that the holonomy maps defined in this paper are holomorphic. Observe that, a posteriori, the holonomy maps in the two papers coincide on $\mathbb M_0(S)$, as the two holomorphic maps coincide on the totally real locus $\mathcal{MGH}(S)$: indeed, the correspondence between minimal immersions into $\Ads^3$ and $G$-Hitchin Higgs bundles corresponding to a pair of quadratic differentials $(q,-q)$ is obtained explicitly via the the non-abelian Hodge correspondence and the Gauss map on $\Ads^3$. %
    
    In fact, we believe that the two holonomy maps coincide on $\mathrm M^*(S)$ via an analogue approach using the Gauss map, even if attention should be paid to the relation with the complex harmonic Higgs bundles setting treated in \cite{ESrank2}. Since these aspects are beyond the scope of this paper, we defer further clarification to some future work of the first author and Sagman.
\end{remark}

\subsection{The holonomy map is a local biholomorphism}\label{sec: holonomy map local biholo}
The aim of this section is to prove the following statement. 
\begin{thm}
\label{thm: hol local biholo}
    The holonomy map $\mathrm{hol}:\mathbb M(S)\to \chi^{*}\big(\pi_1(S), \PSL(2,\C)^2\big)$ is a local biholomorphism.
\end{thm}
To prove this, we will determine the kernel of the differential map of the holonomy map $\mathrm{hol}\colon \mathrm M^*(S)\to \chi^* (\pi_1(S),\PSL(2,\C)^2)$, or, equivalently, the infinitesimal deformations of a regular minimal immersion that do not change the holonomy at first order. We introduce some framework for this purpose. 

\subsubsection{The tangent space to the character variety}\label{subsubsec:H^1}
First, both in the proof and later in the paper we will use the following well-known description of the tangent space of the character variety. Suppose that $[\rho]$ is a smooth point of $\chi(\pi_1(S),G)$, for $G$ a real semisimple Lie group.  Then one has the following canonical identification (see \cite[\S 1.5]{goldman1984} for details):

$$T_{[\rho]}\chi(\pi_1(S),G)\cong H^1_{\mathrm{Ad} \rho}(\pi_1(S),\mathfrak g)~.$$

We recall that the first cohomology group $H^1_{\mathrm{Ad} \rho}(\pi_1(S),\mathfrak g)$ is the quotient real vector space $Z^1_{\mathrm{Ad} \rho}(\pi_1(S),\mathfrak g)/B^1_{\mathrm{Ad} \rho}(\pi_1(S),\mathfrak g)$, where :
\begin{itemize}
    \item $Z^1_{\mathrm{Ad} \rho}(\pi_1(S),\mathfrak g)$ is the space of cocyles, namely the functions $\tau:\pi_1(S)\to \mathfrak g$ satisfying $\tau(\gamma\eta)=\rho(\gamma)\tau(\eta)+\tau(\gamma)$.
    \item $B^1_{\mathrm{Ad} \rho}(\pi_1(S),\mathfrak g)$ is the space of coboundaries, namely the functions $\tau:\pi_1(S)\to \mathfrak g$ satisfying $\tau(\gamma)=V-\mathrm{Ad} \rho(\gamma)V$ for some $V\in\mathfrak g$.
\end{itemize}
The identification is induced by the map sending an infinitesimal variation $\dot \rho$  to $\dot \rho \cdot \rho^{-1}$. It is easily checked that the cocycle condition is verified as a consequence of the group morphism property, while the infinitesimal variations arising from the action by conjugaction correspond to coboundaries. Hence the map $\dot \rho\mapsto \dot\rho\cdot\rho^{-1}$ descends to
 $T_{[\rho]}\chi(\pi_1(S),G)$, and is a vector space isomorphism with $H^1_{\mathrm{Ad}\rho}(\pi_1(S),\mathfrak g)$.

When $G$ is a \emph{complex} semisimple Lie group, $\chi(\pi_1(S),G)$ is a \emph{complex} affine algebraic variety. If $[\rho]$ is again a smooth point,  there is a completely analogous \emph{complex} vector space isomorphism
$$T^{1,0}_{[\rho]}\chi(\pi_1(S),G)\cong H^1_{\mathrm{Ad} \rho}(\pi_1(S),\mathfrak g)~.$$
See \cite[\S 13]{sikora2012} for details.
In the proof of Theorem \ref{thm: hol local biholo}, and later in the next sections, we will use this descriptions of the holomorphic tangent space of the character variety for $G=\PSL(2,\C)^2$.

For future use, we also recall that, by the de Rham theorem, we have the following isomorphism of (real or complex) vector spaces:
\begin{equation}\label{eq: iso H1 group derham}
    H^1_{\mathrm{Ad} \rho}(\pi_1(S),\mathfrak g)\cong H^1_{\mathrm{dR}}(\pi_1(S),\mathfrak g_{\mathrm{Ad}\rho})~,
\end{equation}
where $H^1_{\mathrm{dR} \rho}(\pi_1(S),\mathfrak g_{\mathrm{Ad}\rho})$ denotes the first de Rham coholomology with values in the flat vector bundle $$\widetilde S\times_{\mathrm{Ad}\rho} \mathfrak g:=(\widetilde S\times \mathfrak g)/\sim$$ and $\sim$ is the equivalent relation defined by $$(p,V)\sim (\gamma(p),\mathrm{Ad}(\rho(\gamma)V)$$ for $\gamma$ in $\pi_1(S)$. 
This allows  to see the elements of  $T^{1,0}_{[\rho]}\chi(\pi_1(S),G)$ as $(\widetilde S\times_{\mathrm{Ad}\rho} \mathfrak g)$-valued 1-forms, on which (when $G$ is complex) the almost complex structure of $\chi(\pi_1(S),G)$ acts, trivially, by multiplication by $i$. 

Moreover, making use of the isomorphism \eqref{eq: iso H1 group derham}, one can define the Atiyah-Bott-Goldman symplectic form $\omega_G$ on (the smooth part of) $\chi(\pi_1(S),G)$, by taking the wedge of two 1-forms with values in $\widetilde S\times_{\mathrm{Ad}\rho} \mathfrak g$, paired with a non-degenerate bilinear form on $\mathfrak g$, typically (a multiple of) the Killing form. In the cases of interest in this work, $G$ is the group of projective 2-by-2 matrices with values in $\R$, $\C$ or, in Section \ref{subsec: bi-complex Goldman}, the algebra $\C_\tau$. We will define the Goldman form using the trace as the pairing in $\mathfrak g$, as follows:

\begin{equation}\label{eq: Goldman symplectic form}\omega_G([\alpha],[\beta])=\int_S \tr(\alpha\wedge\beta)~.
\end{equation}



\subsubsection{Space of equivariant immersions} We define
$$\widehat{\mathscr U} := \{(\sigma,\rho)\ | \ \sigma\colon \widetilde S\to \PSL(2,\C)\ \text{ is a $\rho$--equivariant immersion for some $\rho\in \mathrm{Hom}^{*}\big(\pi_1(S), \PSL(2,\mathbb{C})^2\big)$}  \}.$$
The set $\widehat{\mathscr U}$ can be equipped with a complex Fréchet structure, as we now describe. Let $E$ be the set 
$$
E=\widetilde S\times \PSL(2,\C)\times \mathrm{Hom}^{*}\big(\pi_1(S), \PSL(2,\mathbb{C})^2\big) /\sim ,
$$ 
where $\sim$ denotes the equivalence relation
$$(x, A, \rho) \sim (\gamma(x), \rho(\gamma)A, \rho), \ \forall\gamma\in \pi_1(S) .$$ 

The space $E$ admits a natural structure of fiber bundle over $S$, whose bundle map $E \to S$ is given by the composition of the projection onto the first component with the universal covering map onto $S$. Also, observe that the fibers are complex manifolds, as they are identified with copies of $\PSL(2,\C)\times \mathrm{Hom}^{*}\big(\pi_1(S), \PSL(2,\mathbb{C})^2\big)$. Every element $(\sigma,\rho)\in  \widehat{\mathscr{U}}$ induces a unique section of the bundle $E\to S$ given by $[x]\mapsto [(x,\sigma(x),\rho)]$. 
Vice versa, any section $\delta \colon S\to E$ lifts to a unique map
$$\widetilde\delta=(\id,\delta_2, \delta_3)\colon \widetilde S \to \widetilde S\times \PSL(2,\C)\times \mathrm{Hom}^{*}\big(\pi_1(S), \PSL(2,\mathbb{C})^2\big) . $$
If we require the component $\delta_3$ to be constant, i.e. $\delta_3\equiv \rho$, then $\sigma:=\delta_2\colon \widetilde S\to \PSL(2,\C)$ is a $\rho$--equivariant map. For any such $\delta$, the assignment $\delta \mapsto (\delta_2, \delta_3)$ identifies $\widehat{\mathscr U}$ with a complex Fréchet submanifold of the complex Fréchet manifold $\Gamma(E)$ of smooth sections of the bundle $E\to S$. In other words, $\widehat{\mathscr U}$ corresponds to the subset of sections in $\Gamma(E)$ whose lifts have third component that is constant.

Now, recall that an immersion 
\(
\sigma : S \to \mathrm{PSL}(2,\mathbb{C})
\)
is called  positive if 
\(
\sigma^*\langle\cdot,\cdot\rangle|_{\C TS}
\)
is a complex positive metric on $\mathbb{C}TS$ (see Definition~\ref{def: positive} and Definition~\ref{def: admissible and positive immersion}).
We consider the space
\[
\mathscr U:=\left\{ \sigma\colon \widetilde S\to \PSL(2,\C)\ |\ \sigma \text{ is positive and equivariant for some $\rho\in \mathrm{Hom}^{*}\big(\pi_1(S), \PSL(2,\mathbb{C})^2\big)$} \right\}.
\]
Since, as observed right after Theorem~\ref{thm: integration GC}, every element $\sigma\in \mathscr U$ is equivariant for exactly one element $\rho\in {\mathrm{Hom}^{*}\big(\pi_1(S), \PSL(2,\mathbb{C})^2\big)}$ (and $\rho$ depends smoothly on $\sigma$, see Proposition~\ref{prop: holo smooth dependence}), we can think of $\mathscr U$ as an open subset of $\widehat{\mathscr U}$ by identifying $\sigma$ with $(\sigma,\rho)$ and noticing that positivity is an open condition on the map $\sigma$. As such, $\mathscr U$ inherits a complex Fréchet manifold structure. 

Each $\sigma\in\mathscr{U}$ determines a canonical unit normal vector field $N_{\sigma}$, together with the following tensors (which descend to the closed surface $S$): the (positive) first fundamental form $\Is_{\sigma}$, its associated bicomplex structure $J_{\sigma}$ (see Proposition~\ref{prop:unique_bicomplex_structure}), and the second fundamental form $\II_{\sigma}$.

If we consider $\mathscr U$ as a subset of the space of sections $\Gamma(E)$, we can think of tangent vectors in $T_\sigma \mathscr U$ as vector fields $\dot \sigma\in \Gamma(\widetilde S, \sigma^*T\PSL(2,\C))$ that satisfy a suitable relation that involves the holonomy $\rho$ of $\sigma$ and its corresponding infinitesimal variation $\dot \rho$. Specifically, the map $\dot \sigma$ satisfies
\begin{equation}
    \label{eq: dot sigma e dot rho}
    \dot \sigma |_{\gamma (p)}=   \dot \rho(\gamma) \cdot \sigma_0 (p) + \mathrm{d}(\rho_0(\gamma))_{\sigma(p)}(\dot \sigma |_p) ,
\end{equation}
for all $p\in \widetilde S$. Since the pull-back bundle $\sigma^*T\PSL(2,\C)$ splits as a direct sum \(\C TS \oplus \mathrm{Span}_{\C}(N_{\sigma})\), we can split every element \(\dot\sigma\) into the sum of a tangent and a normal vector field along $\widetilde S$.

\subsubsection{Structure of the proof} Having established this background, we now focus on the proof of Theorem~\ref{thm: hol local biholo}. More precisely, it only remains to prove that the differential of the holonomy map is injective. The proof will consist of two main steps. Suppose that we have a path $\sigma_t:\widetilde S\to \PSL(2,\C)$ of equivariant regular minimal immersions whose holonomy is constant at first order. One step (see Proposition~\ref{prop: minimal deformations with trivial holonomy } below) will show, applying infinitesimal rigidity, that the first variation of $\sigma_t$ is ``trivial'' in a complex sense, that is, $\dot\sigma=\mathrm{d}\sigma(\widetilde X)$ for some $X\in \Gamma(\C TS)$, where $\widetilde X$ denotes the lift to $\widetilde S$. The other step (see Proposition~\ref{prop: trivial sigma variation gives trivial gamma variation}) then shows that if $\dot\sigma=\mathrm{d}\sigma(\widetilde X)$, then the first derivative of the 5-tuple $\gamma_t:=\gamma_{\sigma_t}$ is of the form $\dot\gamma=\mathscr L_X\gamma$. Since the variations of this form are precisely the kernel of the projection from $\mathrm M^*(S)$  to $\mathbb M(S)$, this will conclude the proof.

It will be convenient, however, to treat the two steps in the reverse order. So, we begin with the following statement.

\begin{prop}\label{prop: trivial sigma variation gives trivial gamma variation}
    Let $\sigma_t\colon \widetilde S\to \PSL(2,\C)$ be a smooth path of regular minimal immersions with immersion data $\gamma_t$. If $\dot \sigma=\mathrm{d}\sigma(\widetilde X)$ for some $X\in \Gamma(\C TS)$, then $\dot\gamma=\mathscr L_X\gamma$.
\end{prop}

We split the proof of Proposition~\ref{prop: trivial sigma variation gives trivial gamma variation} into a series of lemmas.

\begin{lem}
\label{lem: mappa I II olomorfa}
The map 
\begin{equation}
\label{eq: mappa I, II di sigma}
\begin{split}
(\Is,\II)\colon \mathscr U &\to \Gamma(\mathrm{Sym}_2(\C TS))^2 \\
\sigma&\mapsto (\Is_{\sigma}, \II_{\sigma})
\end{split}
\end{equation}
is holomorphic.
\end{lem}
\begin{proof} 
Denote by $ \mathtt b $ the Killing form on $\mathfrak{sl}(2,\C)$. 
For $\sigma\in \mathscr U$, denote by $N_{\sigma}$ a unit normal vector field. Let $v\in \C T_xS$ be any vector such that $\Is_{\sigma}(v,v)\ne 0$, so that the triple $(v, \mathcal Jv, N_{\sigma})$ is an orthogonal basis of $T_{\sigma(x)}\PSL(2,\C)$. Let $\dot \sigma\in T_{\sigma}\mathscr U$ and denote by $\dot \Is_\sigma$ and $\dot \II_\sigma$ the infinitesimal variations of the first and the second fundamental form, respectively, as defined in the statement. 
By definition of the metric on $\PSL(2,\C)$, we have that, for all $w\in \C T\widetilde S$
\begin{equation}\label{eq: metrica complessa positiva in Lemma}
\Is_{\sigma}(w,w)=\mathtt b \left(L_{\sigma(x)}^* \left(\mathrm{d}\sigma(w)\right), L_{\sigma(x)}^* \left(\mathrm{d}\sigma(w)\right)\right).
\end{equation}

It follows that a smooth curve $\sigma_t$ in $\mathscr U$ determines a family of  tensors $\Is_{\sigma_t}$ that depend smoothly on $t$.
Similarly, $\mathcal J_{\sigma}(v)$ is characterized by the following equations:
    \begin{equation}
    \label{eq: bicomplessa in Lemma}
    \begin{aligned}
        &\Is_{\sigma}(\mathcal J_{\sigma}(v),\mathcal J_{\sigma}(v))=\Is_{\sigma}(v,v),\\
        & \Is_{\sigma}(\mathcal J_{\sigma}(v),v)=0.\end{aligned}\end{equation}
        
We observe that, since $v$ and   $\mathcal J_{\sigma}(v)$      form a basis of $\C TS$ whenever $v\neq 0$, if equations~\eqref{eq: metrica complessa positiva in Lemma} hold for some $v\neq 0$, then they hold for every $v$. 
        As before, a smooth curve $\sigma_t$ in $\mathscr U$ determines smoothly varying tensors $\mathcal J_{\sigma_t}$. 
 We may also repeat the same argument to show that $N_{\sigma_t}$ and, finally, $\II_{\sigma_t}$ depend smoothly on $t$ if $\sigma_t$ does, by observing that $N_{\sigma_t}$ is uniquely determined by:
    
  \begin{equation}
  \label{eq: normale in Lemma}
      \begin{aligned}
    &\mathtt b\left(L_{\sigma}^* \left(\mathrm{d}\sigma(v)\right), L_{\sigma}^* N_{\sigma}\right)=0\\
    &\mathtt b \left(L_{\sigma}^*\left( \mathrm{d}\sigma(J_{\sigma} v)\right),  L_{\sigma}^*N_{\sigma}\right) =0\\
    &\mathtt b \left(L_{\sigma}^*N_{\sigma} ,  L_{\sigma}^*N_{\sigma}\right) =1, 
\end{aligned}
 \end{equation}
while $\II_{\sigma_t}$ is defined as:

     \begin{equation}\label{eq: II in Lemma}\begin{aligned} 
        \II_{\sigma}(w,v) &=  -\Is_{\sigma} \left(w, \sigma^*\left(D_{\mathrm{d}\sigma(v)} N_{\sigma}\right)\right)  \\ 
        &=\mathtt b \left( L_{\sigma(x)}^*\left(\mathrm{d}\sigma(v) \right),L_{\sigma(x)}^*\left(D_{\mathrm{d}\sigma(v)} N_{\sigma}\right)\right).
    \end{aligned}\end{equation}
       
So, by the Generalized Boman's Theorem (see Theorem~\ref{thm: boman}), the map~\eqref{eq: mappa I, II di sigma} is $C^{\infty}$.

       Moreover, a direct computation by differentiating the equations~\eqref{eq: metrica complessa positiva in Lemma},~\eqref{eq: bicomplessa in Lemma},~\eqref{eq: normale in Lemma} and~\eqref{eq: II in Lemma}, an infinitesimal variation in the direction $i\dot \sigma$ produces the infinitesimal variations $i\dot \Is$, $i\dot{\mathcal J}_{\sigma}$, $i \dot N_{\sigma}$, $i\dot \II$.  Hence the differential of the map~\eqref{eq: mappa I, II di sigma} is $\C$--linear, and we therefore conclude that it is also holomorphic. 
\end{proof}

Now, it is an evident fact that a variation of an immersion along a vector field that is in the (real) tangent space corresponds to a variation of all the geometric quantities, such as the first and second fundamental form, by a Lie derivative. The next result describes in a similar manner the infinitesimal variation of $(\Is,\II)$ along directions in the \emph{complexified} tangent space.
 
\begin{lem}
\label{lem: complex tangential deformations}
Let $\sigma\in \mathscr U$, and assume that $\dot \sigma\in T_{\sigma}\mathscr U$ is such that $\dot \sigma= \mathrm{d}\sigma(\widetilde X)$, for some $X\in \Gamma(\C TS)$. Then the corresponding infinitesimal variations of the first and second fundamental forms are given by $$\dot \Is_\sigma= \mathscr L_X \Is_\sigma , \quad \quad \dot \II_\sigma= \mathscr L_{X} \II_\sigma.$$
\end{lem}
  \begin{proof}
  Let $\mathrm X\in \Gamma(TS)$ and let $\phi_{t\widetilde{\mathrm X}}$ denote the lift to the universal cover of the corresponding flow. Then, $\frac{\mathrm{d}}{\mathrm{d}t}|_{t=0} ( \sigma\circ \phi_{t\widetilde{\mathrm X}}) = \mathrm{d}\sigma(\widetilde{\mathrm X})\in \Gamma(\sigma^* T\PSL(2,\C))$. In particular, we can think of $\mathrm{d}\sigma(\widetilde{\mathrm X})$ as an element of $T_{\sigma}\mathscr U$. By the properties of the (real) Lie derivative, this infinitesimal deformation determines an infinitesimal variation of $\Is_\sigma$ and $\II_\sigma$ according to the Lie derivative through $\mathrm X$. By Lemma~\ref{lem: mappa I II olomorfa}, for all $X\in \Gamma(\C TS)$
  \begin{align*}
  \mathrm{d}(\Is, \II)(\mathrm{d}\sigma(X))&=  \mathrm{d}(\Is, \II)(\mathrm{d}\sigma(Re(X)))+i \, \mathrm{d}(\Is, \II)(\mathrm{d}\sigma(Im(X))) \\
  &=(\mathscr L_{Re(X)}\Is, \mathscr L_{Re(X)}\II) +i \, (\mathscr L_{Im(X)}\Is, \mathscr L_{Im(X)}\II)\\& =(\mathscr L_{X}\Is, \mathscr L_{X}\II)
  \end{align*}
as claimed. \end{proof}


\begin{lem}
\label{lem: map from M* to I and II}
    The map introduced in Equation~\eqref{eq: from M* to I II} has injective differential. Moreover, its differential at $\gamma$ maps $\mathscr L_X \gamma$ to $(\mathscr L_X\Is_{\gamma}, \mathscr L_X{\II}_\gamma)$ for all $X\in \Gamma(\C TS)$.
\end{lem}

\begin{proof}
We start by proving that the differential map is injective. Consider a path $\gamma^t=(J_1^t, q_1^t, -J_2^t, \overline q_2^t, u^t)\in \mathrm M^*(S)$, let $\Is_t$ and $\II_t$ be the image path and assume that both derivatives $\dot \Is, \dot \II$ vanish at $t = 0$. Let $z_t$ and $\overline {w_t}$ be a path of local coordinates for $J_1^t$ and $-J_2^t$, respectively. Also denote by $\mathcal J_t$ the path of bi-complex structures induced by $(J_1^t,-J_2^t)$. 
 
 We show that the derivative $\dot{\mathcal J}$ of $\mathcal J_t$ at time zero is such that $\dot{\mathcal J}=0$. The condition $\mathcal J_t^2=- id$ implies $\mathcal J\dot{\mathcal J}+\dot{\mathcal J}\mathcal J=0$, from which we see that $\dot {\mathcal J}$ maps the $i$--eigenspace to the $(-i)$ eigenspace and vice versa: equivalently, $\dot{\mathcal J}(\partial_{\overline z})
 \in \mathrm{Span}(\partial_{w})$ and $\dot{\mathcal J}(\partial_{w})\in \mathrm{Span}(\partial_{ \overline z})$. On the other hand, from $\Is_t= \varrho_t \mathrm{d}z_t\mathrm{d}\overline {w_t}$, 
 we get $\Is_t(\mathcal J_t\cdot, \cdot)+\Is_t(\cdot, \mathcal J_t\cdot)=0$, hence $\dot \Is=0$ implies:  $$\Is( \dot {\mathcal J}\cdot, \cdot)+\Is(\cdot, \dot{\mathcal J})=0.$$ By evaluating the left-hand side in $(\partial_{\overline z}, \partial_{\overline z})$, we have that $\Is(\dot {\mathcal J}(\partial_{\overline z}), 
 \partial_{\overline z})=0$: using that $\dot {\mathcal J}(\partial_{\overline z})\in Span(\partial_ w)$, $\Is(\partial_{\overline z}, \partial_{\overline z})=0$, and that $\Is$ is non-degenerate, we 
 conclude that $\dot{\mathcal J}(\partial_{\overline z})=0$. In the same fashion, $\dot{\mathcal J}(\partial_{ w})=0$, hence $\dot{\mathcal J}=0$.

Now, observe that, for any function $v$ we can write
 \[
 \mathrm{d}v=\frac{\partial_w v}{\partial_w z} \mathrm{d}z +\frac{\partial_{\overline z} v}{\partial_{\overline z} \overline w} \mathrm{d}\overline w, \] 
 as one can easily see by evaluating both sides in the basis $\{\partial_{\overline z}, \partial_w\}$. Hence,
$\mathrm{d}v\circ \mathcal J= i \mathrm{d}v$ (resp. $\mathrm{d}v\circ \mathcal J= -i \mathrm{d}v$) if and only if $v$ is holomorphic for  $J_1$ (resp. $-J_2$).
 
 By differentiating the expressions $\mathrm{d}z_t\circ \mathcal J_t= i \mathrm{d}z_t$ and $\mathrm{d}\overline w_t\circ \mathcal J_t= -i\mathrm{d}\overline w_t$ at $t=0$, we get that
 \begin{equation}
     \label{eq: infinitesimal dz and dw}
     \mathrm{d}\dot z\circ \mathcal J= i\mathrm{d}\dot z\qquad \mathrm{d}\dot{\overline w}\circ \mathcal J=-i \mathrm{d}\dot{\overline w}.
 \end{equation}
 We therefore conclude that $\dot z$ is $J_1$--holomorphic, hence the corresponding Beltrami differential is zero, implying that $\dot {c}_1=0$. In the same spirit, we get that $\dot{\overline w}$ is $(-J_2)$--holomorphic and $\dot{\overline {J_2}}=0$. 

 Since $\dot z$ and $\dot{\overline w}$ are holomorphic for $J_1$ and $-J_2$ respectively, $\dot q_1$ and $\dot{ \overline{q_2}}$ are multiples of $\mathrm{d}z^2$ and $\mathrm{d}\overline w^2$ respectively, hence $0=\dot \II=\dot q_1+\dot{ \overline{q_2}}$ implies that $\dot q_1=0=\dot{ \overline{q_2}}$.

Finally, since the projection map $\mathrm M^*(S)\to \mathrm{QD}(S)\times\mathrm{QD}(\overline S)$ has injective differential, we conclude that $\dot \gamma=0$.

To prove the second part of the statement, observe that the map defined in Equation~\eqref{eq: from M* to I II} is clearly equivariant for the action of $\mathrm{Diff}_0(S)$, hence its derivative maps $\mathscr L_{\mathsf X}\gamma$ to $(\mathscr L_{\mathsf X} \Is_{\gamma}, \mathscr L_{\mathsf X} \II_{\gamma} )$ for all $\mathsf X\in \Gamma(TS)$. For $X\in \Gamma(\C TS)$, we get the analogue result by writing $\mathscr L_X= \mathscr L_{\Re(X)}+ i\mathscr L_{\Im(X)}$ and using that the differential map is $\C$--linear.
\end{proof}

We have reached the conclusion of Proposition~\ref{prop: trivial sigma variation gives trivial gamma variation}.

\begin{proof}[Proof of Proposition~\ref{prop: trivial sigma variation gives trivial gamma variation}]
    By Lemma~\ref{lem: complex tangential deformations}, the fundamental forms vary infinitesimally according to $\dot \Is=\mathscr L_X \Is$ and $\dot \II=\mathscr L_X\II$. 
Since the differential map of $\gamma\mapsto (\Is_{\gamma}, \II_{\gamma})$ is injective and sends $\dot \gamma$ into $(\mathscr L_X \Is, \mathscr L_X\II)$, we conclude that $\dot \gamma=\mathscr L_X\gamma$ thanks to Lemma~\ref{lem: map from M* to I and II}. 
\end{proof}

Let us move on to the other step, namely the following statement.

\begin{prop}
\label{prop: minimal deformations with trivial holonomy }
  Consider the map
 \begin{align*}   
 \mathcal R: \ &\mathscr U \to \mathrm{Hom}^*\big(\pi_1(S),\PSL(2,\C)^2\big)\times C^{\infty}(S,\C) \\ & \ \ \sigma\longmapsto \ \ \big(\rho, H_\sigma\big) \end{align*}
 where $\rho=\mathrm{hol}(\sigma)$ is the holonomy of the immersion and $H_\sigma$ is the function on $S$ whose lift to $\widetilde S$ is the mean curvature $\mathrm{\tr}_{\Is_\sigma}(\II_\sigma)$ of $\sigma$. 
 
 If $\sigma$ is an equivariant regular minimal immersion that is infinitesimally rigid (or, equivalently, if $\gamma_\sigma \in \mathrm{M}^*(S)$), then $\mathrm{d}\mathcal R_{\sigma} (\dot \sigma)=(0,0)$ if and only if $\dot \sigma=\mathrm{d}\sigma(\widetilde X)$ for some $X\in \Gamma(\C TS)$.
\end{prop}
\begin{proof}
Let $\sigma: \widetilde{S} \to \PSL(2,\C)$ be an equivariant, regular minimal immersion. We first prove that $\mathrm{d}\sigma(\widetilde X)$ is tangent to $\mathscr U$ and that $\mathrm{d}\mathcal R(\mathrm{d}\sigma(\widetilde  X))=(0,0)$ for every $X\in \Gamma(\C TS)$. Observe that $\mathrm{d}\sigma(\widetilde X)$ is the derivative at time $t = 0$ of the path of immersions $\sigma_t\in \mathscr U$ given by $\sigma_t (p)= \exp_{\sigma(p)} (t\,\mathrm{d}\sigma(\widetilde X))$ (for sufficiently small values of $|t|$). Since $\mathrm{d}\sigma(\widetilde X)$ is $\rho$--invariant, the maps $\sigma_t : \widetilde{S} \to \PSL(2,\C)$ are equivariant with respect to the same representation $\rho$, independently of $t$. Moreover, by Lemma~\ref{lem: complex tangential deformations} and the properties of the Lie derivative, the infinitesimal deformation of the mean curvature along $\dot \sigma$ is given by 
$$
\left. \frac{\mathrm{d}}{\mathrm{d}t} \right|_{t=0} \widetilde H_\sigma =\mathscr L_{X}(\mathrm{tr}_{\Is}\II)= \mathscr L_X (H_\sigma)~,
$$
which equals zero since $\sigma$ is minimal. Hence, we conclude that $\mathrm{d}\mathcal R_{\sigma} \big( \mathrm{d}\sigma(\widetilde X)\big)=(0,0)$ for all $X\in \Gamma(\C TS)$.

We now suppose that $\sigma$ is infinitesimally rigid and show that the kernel of $\mathrm{d}\mathcal R_{\sigma}$ consists only of elements of the form $\mathrm{d}\sigma(X)$. 
So, assume $\dot \sigma\in T_{\sigma}\mathscr U $ is such that $\mathrm{d}\mathcal R_{\sigma}(\dot\sigma)=(0,0)$. By Equation~\eqref{eq: dot sigma e dot rho}, we see that $\dot \sigma$ is a $\rho$--invariant section of $\sigma^*T\PSL(2,\C)$, or equivalently a $\rho$--invariant vector field on $\PSL(2,\C)$ along $\sigma$. Denoting by $N$ the normal vector field along $\sigma$ (which is also $\rho$--equivariant), we can write $\dot \sigma= \mathrm{d}\sigma(\widetilde X)+ \widetilde fN$ for some $X\in \Gamma(\C TS)$ and some  complex function $f\colon  S\to \C$. 

Clearly $\widetilde fN\in T_{\sigma}\mathscr U$, since it is the derivative of the path of immersions $\sigma_t=\exp_{\sigma(p)} (t \widetilde fN)$ for $|t|$ small enough. Moreover, by linearity, we get that
\[
\mathrm{d}\mathcal R_{\sigma}(\widetilde f N)= - \mathrm{d}\mathcal R_{\sigma}(\mathrm{d}\sigma(\widetilde X)) = (0,0)~,
\] 
where the last equality holds by the first part of the argument. It remains to show that $\mathrm{d}\mathcal{R}_{\sigma} (\widetilde fN)=(0,0)$ if and only if $f=0$. 

In Appendix~\ref{appendix:inf mean curv} (see in particular Equation~\eqref{eq:linearized H on S}), we show that the $t$-derivative of $H_{\sigma_t}$ at $t=0$ has the following expression, that holds on the closed surface $S$ since all quantities of interest are $\pi_1(S)$-invariant:
$$
\left. \frac{\mathrm{d}}{\mathrm{d}t} \right|_{t=0} H_\sigma = e^{-2u}\partial_u \mathscr G (f)~,
$$
where we recall that the expression for $\partial_u \mathscr G$ has been given in~\eqref{eq: linearized Gauss}. By definition of infinitesimal rigidity, $\partial_u\mathscr G$ is bijective as a linear map from $C^{\infty}(S,\C)$ to itself, and $\dot H_\sigma=0$ by assumption, so we conclude that $f=0$. This concludes the proof. 
\end{proof}

We can finally prove that $\mathrm{hol}:\mathbb M(S)\to\chi^{*}\big(\pi_1(S),\PSL(2,\C)^2\big)$ is a local biholomorphism.

\begin{proof}[Proof of Theorem~\ref{thm: hol local biholo}] We already know that $\mathrm{hol}$ is a holomorphic map by Theorem~\ref{thm: hol from mathbb M is holo}. Therefore, $\mathrm{hol}$ is a local biholomorphism if and only its differential is injective (and, hence, bijective, as domain and codomain have the same dimension) at every point in $\mathbb{M}(S)$. By Theorem~\ref{thm: Kernel projection equal to Lie derivatives} and \cite[Proposition 6.6]{ESrank2}, this is equivalent to proving that, for all $\gamma\in \mathrm M^*(S)$, the differential of $\mathrm{hol}\colon \mathrm M^*(S)\to \chi^{*}\big(\pi_1(S),\PSL(2,\C)^2\big)$ at $\gamma$ has kernel given uniquely by elements of the form $\mathscr L_X\gamma$ for $X\in \Gamma(\C TS)$. 

Consider $\dot \gamma\in T_{\gamma}\mathrm M^*(S)$ and assume that $\mathrm{d} \mathrm{hol}_{\gamma}(\dot\gamma)=0$. For simplicity, fix a path $\gamma_t\in \mathrm M^*(S)$ starting at $\gamma$ with tangent vector $\dot \gamma$. Let $\sigma_t\colon \widetilde S\to \PSL(2,\C)$ be a smooth path of regular minimal immersions with immersion data $\gamma_t$ (constructed, for instance, using Proposition~\ref{prop: holo smooth dependence}). Denoting by $\rho_t\in \mathrm{Hom}\big(\pi_1(S), \PSL(2,\C)^2\big)$ the holonomy of $\sigma_t$, and recalling the description of the tangent space to the character variety from \ref{subsubsec:H^1}, the fact that $\mathrm{d} \mathrm{hol}_{\gamma}(\dot \gamma)=0$ implies that there exists $V\in \mathfrak g:=\mathrm{Lie}\big(\PSL(2,\C)^2\big)$ such that 
\[
\dot \rho \cdot {\rho}^{-1} = V- \mathrm{Ad}(\rho) (V)  \colon \pi_1(S)\to \mathfrak g~. 
\]
Hence, by replacing $\sigma_t$ with $\sigma'_t= \exp(-tV)\cdot  \sigma_t$, we get that $\sigma_t\in \mathscr U$ has holonomy $\rho_t'= \exp(-tV) \rho_t \exp(tV)$ and 
\[
\dot \rho' \cdot (\rho')^{-1}=  \mathrm{Ad}(\rho)(V)-V +\dot \rho =0~.
\]
Since, by assumption, $\sigma'_t$ is a path of minimal immersions and $\gamma$ is infinitesimally rigid, Proposition~\ref{prop: minimal deformations with trivial holonomy } implies that $\dot \sigma'=\mathrm{d}\sigma(\widetilde X)$ for some $X\in \Gamma(\C TS)$. By Proposition~\ref{prop: trivial sigma variation gives trivial gamma variation}, we conclude that $\dot\gamma=\mathscr L_X\gamma$.
\end{proof}

\section{The complex hyperKähler structure}\label{sec:complex_hyperkahler}

In this section, we construct the complex hyperK\"ahler structure on the space $\mathbb{M}_0(S)$, which, as we shall see, simultaneously extends the complex bi-Kähler structure of Loustau-Sanders \cite{loustau-sanders}, the hyperK\"ahler structure introduced by Donaldson~\cite{donaldson2003} on the almost-Fuchsian locus and the para-hyperK\"ahler introduced by the second, fourth and fifth author~\cite{MazSepTamb} on the Anti-de Sitter locus.

\subsection{Complex hyperKähler structures and properties}
We start by recalling the definition of complex hyperKähler metric which,
to the best of our knowledge, was first introduced in \cite[\S~2.1]{Bridgeland2021}.
\begin{defn}
\label{def: complex hyperkahler}
A \emph{complex hyperKähler structure} on a complex manifold $\mathbb X$ with almost complex structure $\mathbf i$ is the data of a $\mathbf i$--holomorphic Riemannian metric $g$ on $\mathbb X$ (see Definition~\ref{def:holomorphic_riemannian_metric}) and of a triplet of tensors $\mathbb I, \mathbb J,\mathbb K\in \mathrm{End}_{\C}(T^{1,0}\mathbb X)$ that satisfy the following properties:
\begin{itemize}
    \item $\mathbb I, \mathbb J,\mathbb K$ are $\mathbf i$--holomorphic;
    \item $\mathbb I^2= \mathbb J^2=\mathbb K^2=\mathbb{IJK}=-\mathrm{Id}$;
    \item $g(\mathbb I\cdot, \mathbb I\cdot)= g(\mathbb J\cdot, \mathbb J\cdot)= g(\mathbb K\cdot, \mathbb K\cdot)= g(\cdot,\cdot)$;
    \item $D \mathbb I=D \mathbb J= D \mathbb K=0$, where $D$ is the holomorphic Levi-Civita connection of $g$.
\end{itemize}
\begin{remark}\label{rmk: hyperkahler symplectic forms}

We observe that, as in the case of hyperKähler structures, the pairing
$
\omega_{\mathbb L} :=  g(\mathbb L \cdot, \cdot)
$
defines three complex symplectic forms that are now holomorphic with respect to $\mathbf i$, as $\mathbb L$ varies in $\{\Icpx,\Jcpx,\Kcpx\}$. 
\end{remark} 
\end{defn}

Similarly to what happens in the hyperK\"ahler setting, the structures $\mathbb{I}, \mathbb{J}, \mathbb{K}$ satisfy the quaternionic relations on the $\mathbf{i}$--holomorphic vector bundle $T^{1,0}\mathbb{X}$. Therefore, the complex dimension of a manifold endowed with a complex hyperK\"ahler structure is necessarily equal to $4n$. In local $\mathbf{i}$--holomorphic coordinates $(z_1,\dots,z_n)$ on an open set $U \subset \mathbb{X}$, the holomorphic metric can be written as
\(
g_{|_U} = \sum_{i,j} g_{ij}(z)\, \mathrm{d}z_i \mathrm{d}z_j,
\)
where $g_{ij}$ are holomorphic functions on $U$ for all $i,j=1,\dots,n$. Similarly, each of the tensors $\mathbb{I}, \mathbb{J}, \mathbb{K}$ in Definition~\ref{def: complex hyperkahler}, being $\mathbf{i}$--holomorphic, can be written locally as
\(
\sum_{i,j} a_{ij}(z)\, \mathrm{d}z_i \otimes \frac{\partial}{\partial z_j},
\)
where $a_{ij}$ are holomorphic functions for all $i,j=1,\dots,n$. 

\begin{remark}\label{rmk:extending_tensors}
Using the splitting
\(
\C T\mathbb{X} = T^{1,0}\mathbb{X} \oplus T^{0,1}\mathbb{X} \cong T^{1,0}\mathbb{X} \oplus \overline{T^{1,0}\mathbb{X}},
\)
we can define a real tensor $g_{\mathbb{R}}$ on $T\mathbb{X}$ and endomorphisms $\mathbb{I}_{\mathbb{R}}, \mathbb{J}_{\mathbb{R}}, \mathbb{K}_{\mathbb{R}}$ on $T\mathbb{X}$ by setting
\begin{align*}
    &g_\R(\mathsf X,\mathsf Y):=2\Re\big(g(X,Y)\big),\quad \text{for} \ \mathsf X=X+\overline X\  \text{and} \ \mathsf Y=Y+\overline Y, \text{with} \ X,Y\in \Gamma(T^{1,0}\mathbb X) \\ & \mathbb L_\R(\mathsf X):=\mathbb{L}(X)+\overline{\mathbb{L}(X)},\quad \text{for} \ \mathsf X=X+\overline X\  \text{and} \ X\in\Gamma(T^{1,0}\mathbb X), 
\end{align*}
where $\mathbb L$ (resp. $\mathbb L_\R$) denotes one of the endomorphisms $\mathbb I,\mathbb J,\mathbb K$ (resp. $\mathbb I_\R,\mathbb J_\R,\mathbb K_\R$). As can be easily verified, the following properties hold:
\begin{enumerate}
    \item The endomorphisms $\mathbb{I}_{\mathbb{R}}, \mathbb{J}_{\mathbb{R}}, \mathbb{K}_{\mathbb{R}}$ commute with $\mathbf{i}$ and they satisfy the quaternionic relations. Moreover, the $\C$--linear extension of $\mathbb{L}_{\mathbb{R}}$ to $\C T\mathbb{X}$ coincides with $\mathbb{L}$ when restricted over $T^{1,0}\mathbb{X}$, for every $\mathbb{L} \in \{\mathbb{I},\mathbb{J},\mathbb{K}\}$;
    \item $g_{\mathbb{R}}$ is a pseudo-Riemannian metric of neutral signature such that $g_{\mathbb{R}}(\mathbf{i} \cdot, \mathbf{i} \cdot) = - g_{\mathbb{R}}$ and it is compatible with each of the three endomorphisms $\mathbb{I}_{\mathbb{R}}, \mathbb{J}_{\mathbb{R}}, \mathbb{K}_{\mathbb{R}}$ on $T\mathbb{X}$, meaning that $g_{\mathbb{R}}(\mathbb{L}_{\mathbb{R}}\cdot,\mathbb{L}_{\mathbb{R}}\cdot) = g_{\mathbb{R}}$ for all $\mathbb{L} \in \{\mathbb{I},\mathbb{J},\mathbb{K}\}$; observe that the $\mathbb{C}$-bilinear extension of $g_{\mathbb R}$ to $\mathbb C T\mathbb X$ restricts to $g$ on $T^{1,0}\mathbb X$ and to zero on $T^{0,1}\mathbb X$.
    \item the non-degenerate $2$--forms $\omega_{\mathbb{L}}^{\R}:=g_{\R}(\mathbb{L}_{\R}\cdot, \cdot)$, for $\mathbb{L}_{\R} \in \{\mathbb{I}_{\R},\mathbb{J}_{\R},\mathbb{K}_{\R}\}$, are closed and, hence, determine three real symplectic forms;
    
\end{enumerate}

Although not strictly necessary for this paper, we find worthwile to explain that the almost-complex structures $\mathbb{I}_{\R}, \mathbb{J}_{\mathbb{R}}$ and $\mathbb{K}_{\mathbb{R}}$ are integrable in the classical (complex) sense.

\begin{lem}\label{lemma integrable almost complex for complex hyper}
    Given a complex hyperKähler structure on a complex manifold $(\mathbb X,\mathbf i)$, the endomorphisms $\mathbb{I}_{\mathbb{R}}, \mathbb{J}_{\mathbb{R}}, \mathbb{K}_{\mathbb{R}}$ are integrable almost-complex structures.
\end{lem}
\begin{proof}
    We shall show that the Nijenhuis $N_{\mathbb L_{\mathbb R}}$ vanishes identically, for $\mathbb L=\mathbb I,\mathbb J,\mathbb K$. 

For this purpose, just as for almost complex structures defined on the real tangent bundle, let us introduce a Nijenhuis tensor associated with the endomorphism \(\mathbb L\) as follows:
$$N_{\mathbb L}(X,Y):=[X,Y]+\mathbb L\big([\mathbb LX,Y]+[X,\mathbb LY]\big)-[\mathbb LX,\mathbb LY]\in T^{1,0}\mathbb X, \quad \text{for} \ X,Y\in T^{1,0}\mathbb X.$$
Indeed, since the almost complex structure \(\mathbf{i}\) is integrable, the distribution \(T^{1,0}\mathbb X\) is involutive, and therefore \(N_{\mathbb L}\) defines a well-defined \((1,2)\)-tensor on the holomorphic tangent bundle.  As in the classical setting, this condition  follows from the fact that the endomorphisms \(\mathbb L\) are parallel with respect to the holomorphic Levi--Civita connection (Definition~\ref{def: complex hyperkahler}).

Now, as a consequence of the identity 
\begin{equation}\label{eq: Nijen real and complex}
    N_{\mathbb L_{\mathbb R}}(\mathsf X,\mathsf Y)=2\operatorname{Re}\bigl(N_{\mathbb L}(X,Y)\bigr)
\end{equation}
where \(\mathsf X=X+\overline{X}\) and \(\mathsf Y=Y+\overline{Y}\), we deduce that $N_{\mathbb L_{\mathbb R}}=0$. The relation~\eqref{eq: Nijen real and complex} between the Nijenhuis tensors follows directly from the definition of \(\mathbb L_{\mathbb R}\) and from the fact that \([\overline{X},Y]=0\) for \(X,Y\in T^{1,0}\mathbb X\), after choosing local holomorphic extensions.
\end{proof}

In summary, from Remark~\ref{rmk:extending_tensors} and Lemma~\ref{lemma integrable almost complex for complex hyper} the quadruple $(g_{\R}, \mathbb{I}_{\R}, \mathbb{J}_{\mathbb{R}}, \mathbb{K}_{\mathbb{R}})$ has the same algebraic and differential geometric properties as a classical hyperK\"ahler structure on $\mathbb{X}$, except for the fact that $g_{\R}$ is indefinite (compare with \cite{Bridgeland2021}).
\end{remark}

Hitchin described in \cite[Lemma 6.8]{hitchin1987self} a sufficient condition for a Riemannian manifold endowed with three compatible almost-complex structures satisfying the quaternionic identities to define a hyperKähler manifold. A similar statement holds in this setting. Before proving the analogous result, it is worth observing that the eigenspaces of \(\mathbb L\) corresponding to the eigenvalues \(\pm i\), which we denote by \(\operatorname{Eig}_{\pm i}(\mathbb L)\), are by definition contained in the holomorphic tangent bundle \(T^{1,0}\mathbb X\). Their intersection is clearly trivial. Moreover, from the decomposition 
\[
X=\frac{1}{2}\bigl(X-i\,\mathbb L(X)\bigr)+\frac{1}{2}\bigl(X+i\,\mathbb L(X)\bigr) ,
\]
it follows that the eigenspaces corresponding to the eigenvalues \(\pm i\) generate the holomorphic tangent bundle. Hence \(T^{1,0}\mathbb X=\operatorname{Eig}_i(\mathbb L)\oplus\operatorname{Eig}_{-i}(\mathbb L)\). 

\begin{lem}\label{lem:integrability_complex_hyperkahler}
Let $\mathbb X$ be a complex manifold and let $\mathbf i$ be the associated integrable almost-complex structure on $T\mathbb X$. Suppose that there exist three $\mathbf{i}$--holomorphic sections $\Icpx,\Jcpx,\Kcpx$ of $\mathrm{End}_{\mathbb{C}}(T^{1,0}\mathbb X)$ such that $\mathbb I^2= \mathbb J^2=\mathbb K^2=\Icpx\Jcpx\Kcpx=-\mathrm{Id}$ and an $\mathbf i$--holomorphic Riemannian metric $g$ compatible with $\Icpx$, $\Jcpx$, and $\Kcpx$. 

If the 2-forms $\omega_\mathbb L=g(\mathbb L\cdot,\cdot)$ are closed for each $\mathbb L\in\{\Icpx,\Jcpx,\Kcpx\}$, then the quadruple $(\Icpx,\Jcpx,\Kcpx,g)$ defines a complex hyperK\"ahler structure on $\mathbb X$.
\end{lem}
\begin{proof}
We first observe that it is enough to prove that \(N_{\mathbb L}\equiv 0\) for every \(\mathbb L\in\{\mathbb I,\mathbb J,\mathbb K\}\). Indeed, once this is established, the classical formula in \cite[Chapter IX, Proposition 4.2]{kobayashi1996foundations}, extended to the present holomorphic setting, implies directly that if the symplectic forms \(\omega_{\mathbb L}\) are closed, then \(D\mathbb L=0\). In the classical setting of almost complex structures defined on the smooth tangent bundle of a manifold, the vanishing of the Nijenhuis tensor is equivalent to the involutivity of the holomorphic, or equivalently anti-holomorphic, tangent distribution. In our setting, as we shall explain shortly, in order to prove that \(N_{\mathbb L}\equiv 0\), we need to show that both eigendistributions \(\operatorname{Eig}_{\pm i}(\mathbb L)\) are involutive. Indeed, imposing the involutivity conditions
\[
\mathbb L\bigl([X\mp i\,\mathbb L(X),Y\mp i\,\mathbb L(Y)]\bigr)=\pm i\,[X\mp i\,\mathbb L(X),Y\mp i\,\mathbb L(Y)]
\]
and expanding the brackets, one finds that the holomorphic vector field \(N_{\mathbb L}(X,Y)\) belongs to both \(\operatorname{Eig}_i(\mathbb L)\) and  \(\operatorname{Eig}_{-i}(\mathbb L)\). Since these two eigenspaces have trivial intersection, this is possible if and only if \(N_{\mathbb L}(X,Y)=0\).

In the following, we will only show that $\operatorname{Eig}_i(\mathbb K)\subset T^{1,0}\mathbb X$ is involutive, since the computations for the \((-i)\)-eigenspace and the other endomorphisms are analogous. Given $Y,Z\in\Gamma(T^{1,0}\mathbb X)$, we have $$
\omega_{\mathbb I}(Y,Z)=g(\mathbb I(Y),Z)=-g(\Kcpx(Y),\Jcpx(Z))=-\omega_{\mathbb J}(Z,\Kcpx(Y))=\omega_{\mathbb J}(\Kcpx(Y),Z),
$$
hence $\iota_{Y}\omega_{\Icpx}=\iota_{\Kcpx(Y)}\omega_{\Jcpx}$. In particular, $\Kcpx(Y)=i \, Y$ if and only if $\iota_Y\omega_{\Icpx}=i\, \iota_{Y}\omega_{\Jcpx}$. In what follows, we denote by $\mathscr L$ the complex Lie derivative introduced in Section~\ref{subsec: complex lie derivatives}, which acts naturally on any space of tensors defined on $\mathbb X$ with complex values, and differentiates them along a complex vector field in $\mathbb{C}T \mathbb X$. Suppose that $\Kcpx(Y)=i\, Y$ and $\Kcpx(Z)=i\, Z$, then \begin{align*}
    \iota_{[Y,Z]}\omega_{\Icpx}&=\mathscr{L}_Y\big(\iota_{Z}\omega_\Icpx\big)-\iota_Z\big(\mathscr{L}_Y\omega_\Icpx\big) \\ &=\mathscr{L}_Y\big(i\,\iota_{Z}\omega_\Jcpx\big)-\iota_Z\big(\mathscr{L}_Y\omega_\Icpx\big) \\ &=\mathscr{L}_Y\big(i\,\iota_{Z}\omega_\Jcpx\big)-\iota_Z\mathrm d(\iota_Y\omega_\Icpx) \tag{Cartan's magic formula and $\omega_\Icpx$ closed}\\ &=\mathscr{L}_Y\big(i\,\iota_{Z}\omega_\Jcpx\big)-\iota_Z\mathrm d\big(i\,\iota_Y\omega_\Jcpx\big) \\ &=\mathscr{L}_Y\big(i\,\iota_{Z}\omega_\Jcpx\big)-i\,\iota_Z\mathscr{L}_Y\omega_\Jcpx \tag{Cartan's magic formula and $\omega_\Jcpx$ closed} \\ &=i\,\iota_{[Y,Z]}\omega_\Jcpx,
\end{align*} and this is equivalent to $\Kcpx\big([Y,Z]\big)=i\,[Y,Z].$
\end{proof}

\subsection{The hyperKähler structure of \texorpdfstring{$\mathcal{AF}(S)$}{AF(S)}}\label{sec: hyperKahler almost-fuchsian}

The almost-Fuchsian space $\mathcal{AF}(S)$ was defined in terms of isomorphism classes of embedding data of minimal surfaces in $\mathbb{H}^3$ with principal curvatures in $(-1,1)$. 
Moreover, as explained in Proposition~\ref{prop: AF and MGH are totally real submanifolds}, we obtain an inclusion $\mathcal{AF}(S)\subset\mathbb{M}_0(S)$ as a totally real full dimensional submanifold (see Definition~\ref{def: totally real submanifold}).

Let us briefly recall the definition of hyperKähler structure:

\begin{defn}\label{def: hyperkahler}
    A \emph{hyperKähler structure} on a smooth manifold $\mathbb{Y}$ is the datum of a 4-tuple $(\gpara,\Ipara,\Jpara,\Kpara)$, where $(\gpara,\Ipara)$, $(\gpara,\Jpara)$ and $(\gpara,\Kpara)$ are all Kähler structures, and  $\Ipara,\Jpara,\Kpara$ satisfy the quaternionic identities $\Ipara^2=\Jpara^2=\Kpara^2=\Ipara\Jpara\Kpara=-\mathrm{Id}$.
\end{defn}

The main result we recall here concerns the existence of a natural hyperK\"ahler structure on the space $\mathcal{AF}(S)$, originally introduced by Donaldson \cite{donaldson2003} and later studied also by Hodge \cite{hodge2005hyperkahler} and Trautwein \cite{trautwein2018infinite}.

\begin{thm}\label{thm:hyperkahler_almost_fuchsian}
Let $S$ be a closed oriented surface of genus $\ge 2$, then the almost-Fuchsian space $\mathcal{AF}(S)$ admits a $\mathrm{MCG}(S)$--invariant hyperK\"ahler structure $(\gpara,\Ipara,\Jpara,\Kpara)$. Moreover, the Fuchsian locus $\mathcal{F}(S)$ is an $\mathbf{I}$--holomorphic, totally geodesic submanifold and $(\gpara,\Ipara)$ restricts to a multiple of the Weil-Petersson K\"ahler metric on Teichm\"uller space.
\end{thm}


\subsubsection{Pull-backs from \texorpdfstring{$\mathbb{QD}(S)$}{QD(S)}} Let \(\pi_{\mathcal{AF}}:\mathcal{AF}(S)\to \mathbb{QD}(S)\) be the map obtained by associating to an almost-Fuchsian surface in $\mathbb{H}^3$ the equivalence class of the pair of induced complex structure and holomorphic quadratic differential $[J,q]\in\mathbb{QD}(S)$. This map $\pi_{\mathcal{AF}}$ is well-known to be a diffeomorphism onto its image, see also Remark~\ref{rmk: restriction of pi on AF and MGH}. We now recall some useful and interesting relations concerning Donaldson’s hyperK\"ahler metric that will be needed in the sequel.
\begin{prop}\label{prop: lemma I trautwein}The complex structure $\mathbf{I}$ in Donaldson's hyperKähler metric is obtained as the pull-back via $\pi_{\mathcal AF}$ of the complex structure on $\mathbb{QD}(S)$ acting as the multiplication by $i$ on elements $[\dot \mu,\dot q]$.
\end{prop}

The result above is essentially a consequence of what was explained in Section~\ref{sec: hq differentials}. More precisely, with our conventions, the complex structure \(\mathbf{I}\) differs by a sign from the one defined by Donaldson (see also Remark~\ref{rem: convention on I}).

\begin{prop}[{\cite[Theorems~3.1.1 and~4.2.5]{hodge2005hyperkahler}}]\label{prop: omega J pullback AF}
Let $\Omega^{\C}:=-\mathrm{d}\lambda^{\C}$ be the standard complex symplectic form on $T^*\mathcal T(S)\cong\mathbb{QD}(S)$ obtained from the complex Liouville $1$--form. Then, the complex structure $\mathbf{J}$ in Donaldson's hyperKähler metric satisfies:
 $$\omega_{\mathbf J}=2\pi_{\mathcal{AF}}^*\Im(\Omega^{\C}).$$
\end{prop}

\begin{remark}\label{rem: convention on I}
In the original paper by Donaldson \cite{donaldson2003} and subsequently in \cite{trautwein2018infinite}, the complex structure $\mathbf{I}$ when restricted to the Fuchsian locus, differs by a sign from the one introduced in Section~\ref{sec: beltrami differentials}, since it is compatible with the opposite orientation on Teichmüller space (see Proposition~\ref{prop: lemma I trautwein}). In fact, our convention is consistent with the one used by Hodge \cite{hodge2005hyperkahler}. To avoid sign propagation and to simplify notation, in this article we will denote by $\mathbf{I}$ the complex structure which, when restricted to $\mathcal{T}(S)$, coincides exactly with the one defined in Section~\ref{sec: beltrami differentials}
\end{remark}

\subsubsection{Pull-backs from \texorpdfstring{$\chi^{irr}\big(\pi_1(S),\PSL(2,\C)\big)$}{X(pi1(S),PSL(2,C))}}\label{subsec: pull backs from PSL(2,C)}
Recall that the restriction to $\mathcal{AF}(S)$ of the $\PSL(2,\C)$ holonomy map, denoted by \(\mathrm{hol}_{\mathcal{AF}}\), is a homeomorphism onto its image in the character variety $\chi^{irr}\big(\pi_1(S),\PSL(2,\C)\big)$. The latter is endowed with a natural complex structure $J_\chi$, which is defined by multiplication by $i$  on $T_{[\rho]}^{1,0}\chi^{irr}\big(\pi_1(S),\PSL(2,\C)\big)\cong H^1_{\mathrm{Ad}\rho}(\pi_1(S),\mathfrak{sl}(2,\C))$ (see subsection \ref{subsubsec:H^1}). We have the following relations:

\begin{prop}[{\cite[Propositions 5.2.6 and 5.3.5]{hodge2005hyperkahler}}]\label{prop: Hodge structure J and omega G}
The following identities hold for Donaldson's hyperKähler structure $(\gpara,\Ipara,\Jpara,\Kpara)$:
\begin{itemize}
    \item[(i)] $(\mathrm{hol}_{\mathcal{AF}})^* J_\chi=-\mathbf J$;
    \item[(ii)] $2(\mathrm{hol}_{\mathcal{AF}})^*\omega_G=\omega_{\mathbf I}+i\,\omega_{\mathbf K}$,
\end{itemize}where $\omega_G$ is the Goldman's symplectic form on $\chi^{irr}\big(\pi_1(S),\PSL(2,\C)\big)$ and $\omega_{\mathbf L}=\mathbf g(\mathbf L\cdot,\cdot)$ for $\mathbf L\in\{\mathbf I,\mathbf J,\mathbf K\}$.
\end{prop}

We will soon make use of these results, when dealing with the complex hyperKähler structure on $\mathbb{M}_0(S)$. We thus need to compare the different holonomy maps. From what has been explained in Subsection~\ref{sec: PSL(2,C)}, we have the following commutative diagram:

\begin{equation}\label{eq:comm diag AF}
\centering
\begin{tikzcd}
\mathcal{AF}(S) \arrow[r, ] \arrow[d, "\mathrm{hol}_{\mathcal{AF}}"'] & \mathbb{M}_0(S) \arrow[d, "\mathrm{hol}"] \\
\chi^{irr}\big(\pi_1(S),\PSL(2,\C)\big) \arrow[r, "\Delta^\dagger"]                & \chi^{*}\big(\pi_1(S),\PSL(2,\C)^2\big)
\end{tikzcd}
    \end{equation}
where the top arrow is the usual inclusion of $\mathcal{AF}(S)$ 
(Proposition~\ref{prop: AF and MGH are totally real submanifolds}). We recall that $\Delta^\dagger$ is the map defined by
$$\rho\mapsto (\rho,\overline{{}^t\rho^{-1}}).$$

\subsection{The para-hyperKähler structure of \texorpdfstring{$\mathcal{MGH}(S)$}{MGH(S)}}\label{sec:para_hyperkahler}

Similarly to what we discussed in Section~\ref{sec: hyperKahler almost-fuchsian}, the Anti-de Sitter locus $\mathcal{MGH}(S)$ parametrizes isomorphism classes of embedding data of equivariant spacelike maximal surfaces in $\Ads^3$.


The natural metric defined on $\mathcal{MGH}(S)$ is called para-hyperK\"ahler (\cite[\S~1.3]{MazSepTamb}), of which we now recall the definition. We first need to define a para-complex and para-Kähler structure.

\begin{defn}\label{def:para-complex structure}
    A \emph{para-complex structure} on a smooth manifold $\mathbb{Y}$ is the datum of a section $\mathbf{P} \in \Gamma(\mathrm{End}(T \mathbb{Y}))$ that satisfies the following properties:
    \begin{enumerate}
        \item $\mathbf{P}^2 = \mathrm{Id}$;
        \item the $(+1)$--eigenspace $V_y^+$ and the $(-1)$--eigenspace $V_y^-$ of $\mathbf{P}_y$ inside $T_y \mathbb{Y}$ have the same dimension for any $y \in \mathbb{Y}$;
        \item the distributions $\{V_y^+\}_{y \in \mathbb{Y}}$ and $\{V_y^-\}_{y \in \mathbb{Y}}$ are involutive. 
    \end{enumerate}
\end{defn}

By property (1), the endomorphism $\mathbf{P}_y$ is diagonalizable and has eigenvalues equal to $\pm 1$, so the eigenspaces $V_y^+$ and $V_y^-$ describe a direct sum decomposition $T_y \mathbb{Y} = V_y^+ \oplus V_y^-$ for any $y \in \mathbb{Y}$. Combining this with property (2), we deduce that any manifold $\mathbb{Y}$ that admits a para-complex structure is even-dimensional. Furthermore, property (3) is equivalent to requiring that the Nijenhuis tensor $N_{\mathbf{P}}$ of $\mathbf{P}$ vanishes identically (see e.g. \cite[Proposition~1]{MR2061551}). We refer to \cite[\S~2.1]{cruceanu1996survey} for an overview of the general theory of para-complex manifolds.

\begin{defn}
    A \emph{para-Kähler structure} is the datum of a pair $(\gpara,\mathbf P)$ where $\gpara$ is a pseudo-Riemannian metric of neutral signature, $\mathbf P$ is a para-complex structure, and $\gpara(\mathbf P\cdot,\cdot)$ is a symplectic form. 
\end{defn}

We can now recall the definition of para-hyperKähler structure.

\begin{defn}\label{def: para-hyperkahler}
    A \emph{para-hyperKähler structure} on a smooth manifold $\mathbb{Y}$ is the datum of a 4-tuple $(\gpara,\Ipara,\Jpara,\Kpara)$, where $(\gpara,\Ipara)$ is a Kähler structure, while $(\gpara,\Jpara)$ and $(\gpara,\Kpara)$ are para-Kähler structures, and  $\Ipara,\Jpara,\Kpara$ satisfy the para-quaternionic identities $\Ipara^2=-\mathrm{Id}$, $\Jpara^2=\Kpara^2=\Ipara\Jpara\Kpara=\mathrm{Id}$. 
\end{defn}

The main result of \cite{MazSepTamb} is then the following.
\begin{thm}[{\cite[Theorem A]{MazSepTamb}}]\label{thm:para_hyperkahler}
Let $S$ be a closed oriented surface of genus $\ge 2$, then $\mathcal{MGH}(S)$ admits a $\mathrm{MCG}(S)$--invariant para-hyperK\"ahler structure $(\gpara',\Ipara',\Jpara',\Kpara')$. Moreover, the Fuchsian locus $\mathcal{F}(S)$ is an $\mathbf{I}'$--holomorphic, totally geodesic submanifold and $(\gpara',\Ipara')$ restricts to twice the Weil-Petersson K\"ahler metric on Teichm\"uller space. 
\end{thm}

\begin{remark}
\label{rmk: convention I'}
The same situation as in Remark~\ref{rem: convention on I} also applies in this case, due to the convention adopted in \cite{MazSepTamb}, that is, the complex structure $\mathbf{I}'$ used in this exposition is equal to the opposite of the complex structure considered in \cite{MazSepTamb}. 
\end{remark}



All the tensors appearing in the para-hyperK\"ahler structure described above admit a geometric interpretation through the various parametrizations of the space $\mathcal{MGH}(S)$. We recall here those that will be needed in the next section in order to study some properties of the complex hyperK\"ahler metric that we will introduce.

\subsubsection{Pull-backs from \texorpdfstring{$\mathbb{QD}(S)$}{QD(S)}} While not strictly needed for the arguments in the rest of the paper, we briefly notice that the following two propositions, which are the analogue of Propositions \ref{prop: lemma I trautwein} and \ref{prop: omega J pullback AF}, hold for the map $\pi_{\mathcal{MGH}}:\mathcal{MGH}(S)\to \mathbb{QD}(S)$ that associates to an equivariant maximal surface the corresponding pair $[J,q]$.

\begin{prop}[{\cite[Theorem B]{MazSepTamb}}]The complex structure $\mathbf{I}$ compatible with the para-hyperKähler metric on $\mathcal{MGH}(S)$ is obtained as the pull-back via $\pi_{\mathcal{MGH}}$ of the complex structure on $\mathbb{QD}(S)$ acting as the multiplication by $i$ on elements $[\dot \mu,\dot q]$.
\end{prop}

\begin{prop}[{\cite[Theorem B]{MazSepTamb}}]\label{prop: omega J' and Omega MST}
Let $\Omega^{\C}=-\mathrm{d}\lambda^{\C}$ be the standard complex symplectic form on $T^*\mathcal T(S)\cong\mathbb{QD}(S)$ obtained from the complex Liouville $1$--form. Then, the complex structure $\mathbf{J'}$ compatible with the para-hyperKähler metric on $\mathcal{MGH}(S)$ satisfies:
 $$\omega_{\mathbf J'}=2\pi_{\mathcal{MGH}}^*\Re(\Omega^{\C})$$
\end{prop}

\subsubsection{Pull-backs from \texorpdfstring{$\chi^*\big(\pi_1(S),\PSL(2,\R)^2\big)$}{X(pi1(S),PSL(2,C))}} In this setting, the holonomy map from $\mathcal{MGH}(S)$ to $\chi^*\big(\pi_1(S),\PSL(2,\R)^2\big)$, denoted by \(\mathrm{hol}_{\mathcal{MGH}}\), is a homeomorphism onto its image, which is actually the connected component consisting of pairs of Fuchsian representations compatible with the orientation of $S$.


\begin{thm}[{\cite[Theorem C, Corollary H]{MazSepTamb}}]\label{thm: Thm C and Cor H}
The following identities hold for the para-hyperKähler metric $(\gpara',\Ipara',\Jpara',\Kpara')$ on $\mathcal{MGH}(S)$:
\begin{itemize}
    \item[(i)] $(\mathrm{hol}_{\mathcal{MGH}})^*(\mathrm{Id},-\mathrm{Id})=\mathbf J'$;\vspace{0.5em}
    \item[(ii)] $\displaystyle(\mathrm{hol}_{\mathcal{MGH}})^*\big(\pi_l^*\omega_{G}^\R+\pi_r^*\omega_G^\R\big)=\,\omega_{\mathbf I'}$. 
\end{itemize}
where $\omega_G^\R$ denotes the Goldman symplectic form on $\chi^{irr}\big(\pi_1(S),\PSL(2,\R)\big)$ and
\[
\pi_l, \pi_r : \chi^*\big(\pi_1(S),\PSL(2,\R)^2\big) \to \chi^{irr}\big(\pi_1(S),\PSL(2,\R)\big)
\]
denote the projections onto the first and second component of the representations. 
\end{thm}

Similarly to the previous section, from the totally geodesic inclusion of $\Ads^3$ into $\PSL(2,\C)$ (see subsection~\ref{subsubsec:tot geod}), we have a commutative diagram analogous to the almost-Fuchsian case:

\begin{equation}\label{eq:comm diag MGH}
\centering
\begin{tikzcd}
\mathcal{MGH}(S) \arrow[r, ] \arrow[d, "\mathrm{hol}_{\mathcal{MGH}}"'] & \mathbb{M}_0(S) \arrow[d, "\mathrm{hol}"] \\
\chi^*\big(\pi_1(S),\PSL(2,\R)^2\big) \arrow[r, "\iota_\R"]                & \chi^*\big(\pi_1(S),\PSL(2,\C)^2\big) 
\end{tikzcd}
    \end{equation}
where in this case \(\iota_\R\) comes from the obvious inclusion of $\PSL(2,\R)$ into $\PSL(2,\C)$.

\subsection{The tensors \texorpdfstring{$\mathbb I, \mathbb J, \mathbb K$}{I, J, K} on \texorpdfstring{$\mathbb M_0(S)$}{M0(S)}}\label{sec: tensors I,J,K}
We now construct the tensors that will constitute the complex hyperKähler structure on $\mathbb{M}_0(S)$ endowed with the complex structure $\mathbf{i}$ described in Theorem \ref{thm: Kernel projection equal to Lie derivatives}. Recall that the bundle $\mathbb{QD}(S)$ of holomorphic quadratic differentials over Teichmüller space carries a natural complex structure and we identified the holomorphic tangent bundle to Teichm\"uller space with the space of equivalence classes of Beltrami differentials on $S$ (see Section~\ref{sec: hq differentials}). In accordance with the definition of $\mathbb{M}_0(S)$ (Definition $\ref{def:M(S)}$), there is a natural projection
$\pi : \mathbb{M}_0(S) \to \mathbb{QD}(S) \times \mathbb{QD}(\overline{S})$,
which is a local biholomorphism (Theorem $\ref{thm: Kernel projection equal to Lie derivatives}$). This will allow us, for instance, to define one of the tensor of interest on $\mathbb{QD}(S) \times \mathbb{QD}(\overline{S})$ and pull it back to $\mathbb{M}_0(S)$ via the local biholomorphism $\pi$. Let us define the following endomorphism 
\begin{equation}\label{eq: def I in trivialization of pi}
    \begin{aligned}
    \mathbb I_0: T^{1,0}&(\mathbb{QD}(S) \times \mathbb{QD}(\overline{S}))\longrightarrow T^{1,0}(\mathbb{QD}(S) \times \mathbb{QD}(\overline{S})) \\ &([\dot \mu_1, \dot q_1],  [\dot{\bar{\mu_2}}, \dot{\overline{q_2}}])\longmapsto ([i\dot \mu_1, i\dot q_1], [-i\dot{\bar{\mu_2}},  -i\dot{\bar{q_2}}]).
    \end{aligned}
\end{equation}
\begin{lem}\label{lem: I_0 is holomorphic}
The tensor $\mathbb{I}_0$ is holomorphic with respect to the natural complex structure of $\mathbb{QD}(S) \times \mathbb{QD}(\overline{S})$ described above. 
\end{lem}
\begin{proof}
Set $\mathbb{U}:=\mathbb{QD}(S)\times \mathbb{QD}(\overline S)$. We show that \(\mathbb I_0\) is a holomorphic section of \(\operatorname{End}_{\C}(T^{1,0}\mathbb U)\). Since holomorphicity is a local property, let \((z_1,\dots,z_n,w_1,\dots,w_n)\) be local holomorphic coordinates on \(\mathbb U\), where \((z_1,\dots,z_n)\) are holomorphic coordinates on \(\mathbb{QD}(S)\), and \((w_1,\dots,w_n)\) are holomorphic coordinates on \(\mathbb{QD}(\overline S)\). Then $ \frac{\partial}{\partial z_1},\dots,\frac{\partial}{\partial z_n},
\frac{\partial}{\partial w_1},\dots,\frac{\partial}{\partial w_n}$
is a local holomorphic frame of \(T^{1,0}\mathbb{U}\). By definition of \(\mathbb I_0\), one has
$$ \mathbb I_0\left(\frac{\partial}{\partial z_j}\right)=i\frac{\partial}{\partial z_j}, \qquad \mathbb I_0\left(\frac{\partial}{\partial w_k}\right)=-i\frac{\partial}{\partial w_k}.$$
Therefore, with respect to the above holomorphic frame, the matrix of \(\mathbb I_0\) is $$\begin{pmatrix}i\,\mathrm{Id}_{n} & 0\\ 0 &  -i\,\mathrm{Id}_{n} \end{pmatrix}.$$
Its coefficients are constant, hence holomorphic. Thus \(\mathbb I_0\) is a holomorphic section of $\operatorname{End}_{\C}\big(T^{1,0}\mathbb U\big)$, as claimed.
\end{proof}

As a consequence of the previous result, we obtain a $\mathbf i$--holomorphic almost-complex structure on $\mathbb{M}_0(S)$ via
\begin{equation}\label{eq: Complex structure I}
\mathbb{I}:=\pi^*\mathbb{I}_0 : T^{1,0}\mathbb{M}_0(S)\to T^{1,0}\mathbb{M}_0(S).
\end{equation}
\begin{remark}\label{rmk: identification of real tangent bundle}
By Proposition~\ref{prop: AF and MGH are totally real submanifolds}, both $\mathcal{AF}(S)$ and $\mathcal{MGH}(S)$ are totally real submanifolds of $\mathbb M_0(S)$ (with respect to $\mathbf{i}$). There is a natural identification of $T\mathcal{AF}(S)$ and $T\mathcal{MGH}(S)$ as totally real subbundles of $T^{1,0}\mathbb{M}_{0}(S)$ defined by sending a tangent vector $v$ to $\frac{1}{2}(v-i\,\mathbf{i}(v))$. From now on, we will tacitly use these embeddings so that we can consider the action of $\mathbb{I}$ on $T \mathcal{AF}(S)$ and on $T\mathcal{MGH}(S)$. We observe that this is equivalent to considering the restriction of $\mathbb{I}_{\mathbb{R}}$ to the real tangent bundles to the submanifolds (see Remark~\ref{rmk:extending_tensors}). We will use this same strategy for the other tensors $\mathbb{J}$ and $\mathbb{K}$ that will be part of the complex hyperK\"ahler structure and will be defined later in this section. \end{remark}
We prove the following. 

\begin{prop}\label{prop: restriction of I}
    The operator $\mathbb I$ restricts to $\Ipara$ on $T \mathcal{AF}(S)$ and to $\Ipara'$ on $T\mathcal{MGH}(S)$. 
\end{prop}
\begin{proof}
We will carry out the proof for the tensor $\mathbf{I}$, since the argument for $\mathbf{I}'$ is analogous. The complex structure $\mathbf{I}$ is equal to the pull-back, via the map $\pi_{\mathcal{AF}}$ from Section~\ref{sec: hyperKahler almost-fuchsian}, of the endomorphism $[\dot\mu, \dot q]\mapsto [i\dot\mu, i\dot q]$ on the real tangent bundle $T\mathbb{QD}(S)$ (see Proposition~\ref{prop: lemma I trautwein} and Remark~\ref{rem: convention on I}). In accordance with the description of the almost-Fuchsian locus in~\eqref{eq:inclusion AF in M}, any $v \in T\mathcal{AF}(S)\subset T^{1,0}\mathbb M_0(S)$ can be written as $v=[\dot\mu,\dot q,\dot{\overline \mu},\dot{\overline q},\dot u]$. Since $\mathbb{I}$ was constructed as the pull-back of the complex structure $\mathbb{I}_0$ defined in~\eqref{eq: def I in trivialization of pi}, it follows that applying $\mathbb{I}$ to a vector $v \in T\mathcal{AF}(S)$ described as above yields exactly the action of $\mathbf{I}$ on the model $T\mathcal{AF}(S)$ inside $T^{1,0}\mathbb M_0(S)$.
  \end{proof}

We now define the second holomorphic (almost) complex structure that will enter the complex hyperk\"ahler metric on \(\mathbb{M}_0(S)\). To this end, recall that the natural complex structure $J_\chi$ on the $\PSL(2,\C)$--character variety induces the product complex structure $(\dot\rho_1,\dot\rho_2)\mapsto\big( J_\chi(\dot\rho_1),J_\chi(\dot\rho_2)\big)$ on $\chi^*\big(\pi_1(S),\PSL(2,\C)^2\big)$. Since the holomorphic tangent space to the $\PSL(2,\C)$--character variety, at a smooth point $\rho$, can be identified with $H^1(S,\mathfrak{sl}_2(\mathbb{C})_{\mathrm{Ad}\rho})$ (see subsection \ref{subsubsec:H^1}), any tangent vector $\dot\rho$ is represented by a closed $1$--form with values in $\mathfrak{sl}_2(\mathbb{C})_{\mathrm{Ad}\rho}$, that we can express as a finite sum of the form $\Sigma_k \sigma_k\otimes\phi_k$. In particular, the complex structure \(J_\chi\) sends any simple tensor \(\sigma \otimes \phi\) into \(\sigma \otimes (i\phi)\). Let us first define the endomorphism below at the level of simple tensors and then extend it by \(\C\)-linearity \begin{equation}\begin{aligned}\label{eq: definition J_0}
\mathbb J_0:T^{1,0}\chi^*\big(&\pi_1(S),\PSL(2,\C)^2\big)\to T^{1,0}\chi^*\big(\pi_1(S),\PSL(2,\C)^2\big) \\ &(\sigma_1\otimes\phi_1,\sigma_2\otimes\phi_2)\longmapsto(\sigma_1\otimes (-i\phi_1),\sigma_2\otimes i\phi_2).
\end{aligned}\end{equation}
An argument similar to that used in Lemma~\ref{lem: I_0 is holomorphic} yields the following result.
\begin{lem}\label{lem: J_0 holomorphic}
The tensor $\mathbb{J}_0$ is holomorphic with respect to the natural complex structure \((J_\chi,J_\chi)\) of $\chi^*\big(\pi_1(S),\PSL(2,\C)^2\big)$ described above.
\end{lem}
Now we define the tensor $\mathbb J \in \Gamma(\mathrm{End}_{\mathbb{C}}(T^{1,0}\mathbb M_0(S)))$ as follows:
\begin{equation}\label{eq:definition_J} 
\mathbb{J}:=(\mathrm{hol})^*\mathbb{J}_0 : T^{1,0}\mathbb{M}_0(S)\to T^{1,0}\mathbb{M}_0(S),
\end{equation}
which is $\mathbf{i}$--holomorphic by Theorem~\ref{thm: hol local biholo} and Lemma~\ref{lem: J_0 holomorphic}, being the holonomy map a local biholomorphism and hence \(\mathrm{hol}^*(J_\chi,J_\chi)=\mathbf i\).


\begin{prop}\label{prop: restriction of J}
    The tensor $\mathbb J$ restricts to $\Jpara$ on $\mathcal{AF}(S)$, and the tensor $\mathbf i\circ\mathbb J$ restricts to $\Jpara'$ on $\mathcal {MGH}(S)$.
\end{prop}

\begin{proof} 
Regarding the restriction to the almost-Fuchsian locus, it suffices to work at the level of smooth character varieties by diagram~\eqref{eq:comm diag AF}. In fact, 
by Proposition~\ref{prop: Hodge structure J and omega G} point $(i)$ it is enough to prove that the action of the endomorphism $\mathbb J_0$ on the subspace $(\Delta^\dagger)_*\big(H^1(S,\mathfrak{sl}_2(\mathbb{C})_{\mathrm{Ad}\rho})\big)$ is equal to the multiplication by $-i$.
This is a direct computation by using Equation $\eqref{eq: definition J_0}$\begin{align*}
        \mathbb J_0(\sigma\otimes\phi,-\sigma\otimes\overline{{}^t\phi})&=\big(-J_\chi(\sigma\otimes\phi),J_\chi(-\sigma\otimes\overline{{}^t\phi})\big) \\ &=(-\sigma\otimes i\phi,-\sigma\otimes i\, \overline{{}^t\phi}) .
    \end{align*}
    The last term in the equation belongs to the image of \(H^1(S,\mathfrak{sl}_2(\mathbb{C})_{\mathrm{Ad}\rho})\) via $\mathrm{d}\Delta^\dagger$; indeed, 
    \[
    -\overline{{}^t(-\sigma\otimes i\phi)} = - (-\sigma) \otimes (- i \, \overline{{}^t\phi}) = -\sigma\otimes i\, \overline{{}^t\phi} .
    \]
     On the other hand, regarding the restriction to \(\mathcal{MGH}(S)\), we observe that
     \[\mathbf i\circ\mathbb J=\mathrm{hol}^*(J_\chi,J_\chi)\circ\mathrm{hol}^*\mathbb J_0=\mathrm{hol}^*\big((J_\chi,J_\chi)\circ\mathbb J_0\big)=\mathrm{hol}^*(\mathrm{Id,-\mathrm{Id}}).\] Last chains of equalities directly implies the claim by diagram~\eqref{eq:comm diag MGH} and Theorem~\ref{thm: Thm C and Cor H}.
\end{proof}

The third endomorphism is simply obtained as the product
$$\mathbb{K} := \mathbb{I}\mathbb{J},$$
and, therefore, it is $\mathbf{i}$--holomorphic because $\mathbb{I}$ and $\mathbb{J}$ are. We are left to proving that $\mathbb K^2=-\mathrm{Id}$, which is equivalent to $\mathbb I\mathbb J=-\mathbb J\mathbb I$. 

\begin{prop}\label{prop: restriction of K}
   The tensor $\mathbb K$ satisfies $\mathbb K^2=-\mathrm{Id}$. Moreover, $\mathbb{K}$ restricts to $\mathbf K$ on $\mathcal{AF}(S)$ and $\mathbf i\circ\mathbb K$ restricts to $\mathbf K'$ on $\mathcal{MGH}(S)$.
\end{prop}
\begin{proof}
Since $\mathbf{IJ}=\mathbf{K}$ and $\mathbf I'\mathbf{J'}=\mathbf K'$, we easily get that $\mathbb K$ restricts to $\mathbf K$ on $\mathcal{AF}(S)$ and $\mathbf i\circ\mathbb K$ restricts to $\mathbf K'$ on $\mathcal{MGH}(S)$. Now, the tensors $\mathbb K^2$ and $-\mathrm{Id}$ are $\mathbf i$--holomorphic endomorphisms of $T^{1,0}\mathbb M_0(S)$ and they coincide on $\mathcal{AF}(S)$, which is a totally real submanifold of $\mathbb M_0(S)$ of full dimension. By Lemma~\ref{lem: fundamental lem analytic continuation} and Remark~\ref{rmk: fundamental_lemma_applied to tensors}, we conclude that $\mathbb K^2=-\mathrm{Id}$ everywhere on $\mathbb{M}_0(S)$. 
\end{proof}

\subsection{The holomorphic metric and the symplectic structures}\label{subsec: symplectic structures complex hyperkahler}

We can now introduce the symplectic forms and the holomorphic Riemannian metric on $\mathbb M_0(S)$. Recall that the holonomy map $\mathrm{hol}\colon \big(\mathbb M_0(S),\mathbf{i}\big)\to \big(\chi^{*}\big(\pi_1(S),\PSL(2,\C)^2\big),(J_\chi,J_\chi)\big)$ is a local biholomorphism by Theorem~\ref{thm: hol local biholo}. On $\chi^{*}\big(\pi_1(S),\PSL(2,\C)^2\big)$ consider the $(J_\chi,J_\chi)$--holomorphic symplectic form given by $ \pi_r^* \omega_G+ \pi_l^*\omega_G$, where $\omega_G$ is the Atiyah-Bott-Goldman symplectic form on $\chi^{irr}\big(\pi_1(S),\PSL(2,\C)\big)$ (see \cite{goldman1984,atiyah-bott}), while $\pi_r^*$ and $\pi_l^*$ denote the pull-back by the right and left projection. On $\mathbb M_0(S)$, we define
\begin{equation}\label{eq: omega_mathbb I}
    \omega_{\mathbb I}:= \mathrm{hol}^*\left(  \pi_r^* \omega_G+ \pi_l^*\omega_G\right).
\end{equation}
This is clearly a closed 2-form, and, by Theorem~\ref{thm: hol local biholo}, $\omega_{\mathbb I}$ is nowhere-degenerate and $\mathbf i$--holomorphic. It is only left to prove that $\omega_\mathbb I$ is indeed compatible with $\mathbb I$. In the following we will use the formula below for Goldman's symplectic form
\begin{equation}\label{eq: Goldman symplectic form explicit formula}
\omega_G\big(\sum_h\sigma_{1,h} \otimes \phi_{1,h}, \sum_k\sigma_{2,k} \otimes \phi_{2,k}\big)
= \sum_{h,k} \int_S \mathrm{tr}(\phi_{1,h} \cdot \phi_{2,k})\, \sigma_{1,h} \wedge \sigma_{2,k},
\end{equation}which is a consequence of Equation \eqref{eq: Goldman symplectic form} and the description of elements in $H^1(S,\mathfrak{sl}_2(\mathbb{C})_{\mathrm{Ad}\rho})$ (see subsection \ref{subsubsec:H^1}).
 
\begin{lem}\label{lem: compatibility_omega_I}
The complex symplectic form $\omega_\mathbb I$ restricts to $\omega_\mathbf I$ on the almost-Fuchsian locus $\mathcal{AF}(S)$. As a consequence,
\[
    \omega_{\mathbb I} (\mathbb I\cdot, \cdot)+\omega_{\mathbb I} (\cdot,\mathbb I\cdot)=0.
    \]
\end{lem}

\begin{proof}
The tensors $\omega_\mathbb I$ and $\mathbb I$ are both $\mathbf i$--holomorphic and $\mathbb I$ restricts to $\mathbf I$ on $T\mathcal{AF}(S)\subset T^{1,0}\mathbb M_0(S)$ (see Proposition~\ref{prop: restriction of I}). If we show that $\omega_{\mathbb{I}}\vert_{\mathcal{AF}(S)} = \omega_{\mathbf{I}}$, then we can conclude the proof by Lemma~\ref{lem: fundamental lem analytic continuation}, since $\omega_{\mathbf{I}}$ is compatible with $\mathbf{I}$ on $\mathcal{AF}(S)$ by Theorem~\ref{thm:hyperkahler_almost_fuchsian}. In fact, it suffices to prove that
\[
\omega_{\mathbb{I}}\vert_{\mathcal{AF}(S)} =2(\mathrm{hol}_{\mathcal{AF}})^*\Re(\omega_G),
\]
which in turn coincides with $\omega_{\mathbf{I}}$ by Proposition~\ref{prop: Hodge structure J and omega G} point $(ii)$. The differential of the map $$\Delta^\dagger : \chi^{irr}\big(\pi_1(S),\PSL(2,\C)\big) \to \chi^*\big(\pi_1(S),\PSL(2,\C)^2\big),$$  is given by $\dot\rho\mapsto (\dot\rho,-\overline{{}^t\dot{\rho}'})$. Therefore, to understand the restriction \(\omega_{\mathbb{I}}|_{\mathcal{AF}(S)}\), according to diagram~\eqref{eq:comm diag AF}, we need to compute the element
\begin{equation}\label{eq: formula in Lemma for Goldman}
 (\pi_r^* \omega_G)\big((\dot\rho_1,\dot\rho_2),(\dot\rho_1',\dot\rho_2')\big)+(\pi_l^*\omega_G)\big((\dot\rho_1,\dot\rho_2),(\dot\rho_1',\dot\rho_2')\big)
\end{equation}
under the conditions \(\dot{\rho}_2 = -{\overline{{}^t\dot{\rho}_1}}\) and \(\dot{\rho}_2' = -{\overline{{}^t\dot{\rho}_1'}}\). Hence, by using the description of Goldman's symplectic form from Equation~\eqref{eq: Goldman symplectic form explicit formula}, the term in~\eqref{eq: formula in Lemma for Goldman} is equal to
\begin{align*}
    \omega_G\big(-{\overline{{}^t\dot{\rho}_1}},-{\overline{{}^t\dot{\rho}_1'}}\big)+\omega_G(\dot\rho_1,\dot\rho_1')& =  \sum_{h,k} \int_S\Big(\mathrm{tr}(\overline{{}^t\phi_{1,h}} \cdot \overline{{}^t\phi_{1,k}'})+\mathrm{tr}(\phi_{1,h} \cdot \phi_{1,k}')\Big)\, \sigma_{1,h}\wedge\sigma_{1,k}' \\
  &=2\sum_{h,k} \int_S\Re\big(\mathrm{tr}(\phi_{1,h} \cdot \phi_{1,k}')\big)\, \sigma_{1,h}\wedge\sigma_{1,k}' \\
  &=2\Re(\omega_G)(\dot\rho_1,\dot\rho_1'),
\end{align*}
where we selected representatives $\sum_h \sigma_{1,h}\otimes \phi_{1,h}$ and $\sum_k \sigma_{1,k}'\otimes \phi_{1,k}'$ of $\dot{\rho}_1$ and $\dot{\rho}_1'$, respectively. This concludes the proof. 
\end{proof}

Following Definition~\ref{def: complex hyperkahler}, Remark~\ref{rmk: hyperkahler symplectic forms}, and the notations used there, we can define the remaining tensors of the complex hyperKähler structure on $\mathbb M_0(S)$ as follows:
\begin{equation}
\label{eq: g, omegaJ, omegaK}
    \begin{split}
        g&:= \omega_{\mathbb I} (\cdot, \mathbb I\cdot),\\
        \omega_{\mathbb J}&:= g(\mathbb J\cdot, \cdot),\\
         \omega_{\mathbb K}&:= g(\mathbb K\cdot, \cdot).
    \end{split}
\end{equation}
Using what we have proved so far, we know that \(g\) is an $\mathbf i$--holomorphic symmetric and non-degenerate bilinear form on \(T^{1,0} \mathbb M_0(S)\), which is compatible with \(\mathbb I\), in the sense that
\(
g(\mathbb I\cdot,\mathbb I\cdot)= g
\) by Lemma~\ref{lem: compatibility_omega_I}.
Moreover, \(\omega_{\mathbb J}\) and \(\omega_{\mathbb K}\) are indeed $\mathbf i$--holomorphic \(2\)-forms on \(\mathbb M_0(S)\). To obtain a complex hyperKähler metric, it remains to prove that \( g\) is compatible with \(\mathbb J\) and \(\mathbb K\), together with the fact that \(\omega_{\mathbb J}\) and \(\omega_{\mathbb K}\) are closed \(2\)-forms (compare with Lemma~\ref{lem:integrability_complex_hyperkahler}). In fact, it suffices to prove the compatibility between \( g\) and \(\mathbb J\), since the compatibility with \(\mathbb K\) will follow from the relation
\(
\mathbb K=\mathbb I\mathbb J=-\mathbb J\mathbb I .
\) We now provide a more explicit description of $\omega_{\mathbb J}$ and $\omega_{\mathbb K}$ from which it will directly follow that they are both closed \(2\)-forms.

\begin{prop}\label{prop: formula for omega_K}
    The complex symplectic form $\omega_{\mathbb K}$ satisfies:
    \[ \omega_{\mathbb K} =i\mathrm{hol}^*\left( \pi_r^* \omega_G- \pi_l^*\omega_G\right).\] In particular, $\omega_\mathbb K$ is a closed $2$--form.
\end{prop}
\begin{proof}
We begin by observing that $\omega_\mathbb{K}=\omega_\mathbb{I}(\mathbb{J}\cdot,\cdot)$, as the following computation shows:
$$
\omega_\mathbb{K}
= g(\mathbb{K}\cdot,\cdot)
= \omega_\mathbb{I}(\mathbb{K}\cdot,\mathbb{I}\cdot)
= -\omega_\mathbb{I}(\mathbb{I}\mathbb{K}\cdot,\cdot)
= \omega_\mathbb{I}(\mathbb{J}\cdot,\cdot).
$$ In particular, by Proposition~\ref{prop: restriction of J} and Lemma~\ref{lem: compatibility_omega_I}, this implies $\omega_{\mathbb K}|_{\mathcal{AF}(S)}=\omega_{\mathbf I}(\mathbf J\cdot,\cdot)=\omega_\mathbf K=2(\mathrm{hol}_{\mathcal{AF}})^*\big(\Im\omega_G\big)$ (see Proposition~\ref{prop: Hodge structure J and omega G} point $(ii)$). 
At this point, it is sufficient to prove the equality between 
$i\big(\pi_r^* \omega_G- \pi_l^*\omega_G\big)$ and $2\Im\omega_G$ on the image of the $\PSL(2,\C)$--character variety via the map $\Delta^\dagger$ (see diagram~\eqref{eq:comm diag AF}). Indeed, once the equality is established on the almost-Fuchsian locus, the claim follows by applying Lemma~\ref{lem: fundamental lem analytic continuation}. The computation is similar to the one performed in the proof of Lemma~\ref{lem: compatibility_omega_I}, namely
\begin{align*}
i(\pi_r^* \omega_G)\big((\dot\rho_1,-\overline{{}^t\dot{\rho}_1}),(\dot\rho_1',-\overline{{}^t\dot{\rho}_1'})\big)&-i(\pi_l^*\omega_G)\big((\dot\rho_1,-\overline{{}^t\dot{\rho}_1}),(\dot\rho_1',-\overline{{}^t\dot{\rho}_1'})\big)= \\
&=i\Big(\omega_G\big(\overline{{}^t\dot{\rho}_1},\overline{{}^t\dot{\rho}_1'}\big)-\omega_G(\dot\rho_1,\dot\rho_1')\Big) \\
&=i \sum_{h,k} \int_S \Big(\overline{\mathrm{tr}(\phi_{1,h} \cdot \phi_{1,k}')}-\mathrm{tr}(\phi_{1,h} \cdot \phi_{1,k}')\Big)\,\sigma_{1,h}\wedge\sigma'_{1,k} \\ 
&=-2i^2(\Im\omega_G)(\dot\rho_1,\dot\rho_1') \\ &=2(\Im\omega_G)(\dot\rho_1,\dot\rho_1').
\end{align*}
where we selected representatives $\sum_h \sigma_{1,h}\otimes \phi_{1,h}$ and $\sum_k \sigma_{1,k}'\otimes \phi_{1,k}'$ of $\dot{\rho}_1$ and $\dot{\rho}_1'$, respectively.
\end{proof}
\begin{prop}
The $2$--form $\omega_\mathbb J$ is compatible with $\mathbb J$, namely
    $$\omega_\mathbb J(\mathbb J\cdot,\cdot)+\omega_\mathbb J(\cdot,\mathbb J\cdot)=0.$$ 
\end{prop}
\begin{proof}
From Proposition $\ref{prop: restriction of I}$ and Lemma $\ref{lem: compatibility_omega_I}$ we know that the holomorphic metric $g=\omega_{\mathbb I}(\cdot,\mathbb I\cdot)$ restricts to the Riemannian metric $\mathbf g$ on $\mathcal{AF}(S)$. Moreover, by Proposition $\ref{prop: restriction of J}$ we have $\mathbb J|_{\mathcal{AF}(S)}=\mathbf J$. We can therefore deduce that \begin{align*}
\omega_{\mathbb J}(\mathbb J\cdot,\cdot)|_{\mathcal{AF}(S)}&=g(\mathbb J^2\cdot,\cdot)|_{\mathcal{AF}(S)} \\ &=-\mathbf g \\ &=-\mathbf{g}(\mathbf J\cdot,\mathbf J\cdot) \tag{$\mathbf g$ is compatible with $\mathbf J$} \\ &=-g(\mathbb J\cdot,\mathbb J\cdot)|_{\mathcal{AF}(S)} \\&=-\omega_{\mathbb J}(\cdot,\mathbb J\cdot)|_{\mathcal{AF}(S)}
\end{align*}and we conclude by Lemma~\ref{lem: fundamental lem analytic continuation}.
\end{proof}
Recall from Section~\ref{sec: beltrami differentials} that $\mathbb{QD}(S)$ can be identified with the holomorphic cotangent bundle of $\mathcal T(S)$. Let us denote by $\Omega^{\C}=-\mathrm{d}\lambda^{\C}$ the differential of the complex Liouville $1$--form on $\mathbb{QD}(S)$ and by $\overline{\Omega^{\C}}$ the induced 2-form on $\mathbb{QD}(\bar S)$ given by $\overline{\Omega^{\C}}([\dot{\overline \mu},\dot{\overline q}],[\dot{\overline \mu'},\dot{\overline q'}]):=\overline{\Omega^{\C}([\dot\mu,\dot q],[\dot\mu',\dot q'])}$. As a consequence, we have a natural symplectic structure on $\mathbb{QD}(S)\times \mathbb{QD}(\overline S)$ given by
\[
\Omega:=\text{pr}_l^*\Omega^{\C}+\text{pr}_r^*\overline{\Omega^{\C}},
\]
where $\text{pr}_l$ and $\text{pr}_r$ are the left and right projection from the product $\mathbb{QD}(S)\times \mathbb{QD}(\overline S)$ to $\mathbb{QD}(S)$ and $\mathbb{QD}(\bar S)$, respectively. By construction, the symplectic form $\Omega$ is $\mathbb{C}$--linear with respect to the complex structure on $\mathbb{QD}(S)\times \mathbb{QD}(\bar S)$ for which the projection map
\(\pi:\mathbb{M}_0(S)\to \mathbb{QD}(S)\times \mathbb{QD}(\bar S)\)
is a local biholomorphism (see Theorem~\ref{thm: Kernel projection equal to Lie derivatives}). The following result is a generalization of Proposition \ref{prop: omega J pullback AF} and \ref{prop: omega J' and Omega MST}.

\begin{prop}
    The complex symplectic form $\omega_{\mathbb J}$ satisfies:
    \[
    \omega_{\mathbb J}=-\pi^*\Omega(\mathbb{I}\cdot, \cdot).
    \]
    In particular, $\omega_\mathbb J$ is a closed $2$--form.   
\end{prop}
\begin{proof}
Let
\(\dot{\gamma},\dot{\gamma}'\)
be tangent vectors to $\mathbb{QD}(S)\times\mathbb{QD}(\bar S)$ in the image of $\mathrm{d}\pi|_{T\mathcal{AF}(S)}$. Then \begin{align*}
\Omega(\mathbb{I}\dot\gamma,\dot\gamma')&=\Omega([i\dot\mu,-i\overline{\dot\mu},i\dot q,-i\overline{\dot q}],[\dot\mu',\overline{\dot\mu'},\dot q',\overline{\dot q'}]) \\ &=\Omega^{\C}([i\dot\mu,i\dot q],[\dot\mu',\dot q'])+\overline{\Omega^{\C}([i\dot\mu,i\dot q],[\dot\mu',\dot q'])} \\ &=-2\Ima\big(\Omega^{\C}([\dot\mu,\dot q],[\dot\mu',\dot q'])\big).
\end{align*}
By the discussion preceding the statement, \(\pi^*\Omega\) is an \(\mathbf{i}\)-holomorphic symplectic form on \(\mathbb{M}_0(S)\), and therefore Lemma~\ref{lem: fundamental lem analytic continuation}, together with the identity \(2\pi_{\mathcal{AF}}^*\mathrm{Im}(\Omega^{\mathbb{C}})|_{\mathcal{AF}(S)}=\omega_{\mathbf{J}}\) from Proposition~\ref{prop: omega J pullback AF}, implies that
\(\omega_{\mathbb{J}}=-\pi^*\Omega(\mathbb{I}\cdot, \cdot)\) (see also Remark~\ref{rmk: restriction of pi on AF and MGH}).
\end{proof}
At this stage, we have defined three $\mathbf{i}$--holomorphic almost complex structures $\mathbb I$, $\mathbb J$, and $\mathbb K$ on $\mathbb M_0(S)$ such that $\mathbb K=\mathbb I\mathbb J=-\mathbb J\mathbb I$, and the $\mathbf{i}$--holomorphic metric $g$ is compatible with each of them. Moreover, the $2$--forms $\omega_{\mathbb L}:=g(\mathbb L\cdot,\cdot)$ are $\mathbf{i}$--holomorphic and closed. By Lemma
\ref{lem:integrability_complex_hyperkahler}, we obtain the following.
\begin{cor}
The $\mathbf i$--holomorphic complex structures $\mathbb I,\mathbb J,\mathbb K$ together with the $\mathbf i$--holomorphic metric $g$ defines a $\mathrm{MCG}(S)$--invariant complex hyperK\"ahler structure on $\mathbb M_0(S)$.
\end{cor}

The following corollary proves items (2) and (3) of Theorem 1.1.

\begin{cor}\label{cor: restriction of complex hyperkahler to AF e MSH}
The complex hyperK\"ahler structure on $\mathbb M_0(S)$ satisfies the following conditions: \begin{itemize}
    \item[(i)] the quadruple $(\mathbb I,\Jcpx,\Kcpx,g)$ restricts to Donaldson's hyperk\"ahler metric $(\mathbf I,\mathbf J,\mathbf K,\mathbf g)$ on the almost-Fuchsian locus $\mathcal{AF}(S)$;
    \item[(ii)] the quadruple $(\mathbb I,\mathbf i\circ\Jcpx,\mathbf i\circ\Kcpx,g)$ restricts to the para-hyperK\"ahler metric $(\Ipara',\Jpara',\Kpara',\mathbf g')$ of Mazzoli-Seppi-Tamburelli on the Anti-de Sitter locus $\mathcal{MGH}(S)$.
\end{itemize}
\end{cor}
\begin{proof}
Thanks to Propositions~\ref{prop: restriction of I},~\ref{prop: restriction of J} and~\ref{prop: restriction of K}, we already know that the three complex structures $\mathbb{I}, \mathbb{J}, \mathbb{K}$ restrict to the Donaldson's structures on $\mathcal{AF}(S)$, while $\mathbb{I}, \mathbf i{\circ}\mathbb{J}, \mathbf{i}\circ\mathbb{K}$ restrict to $\mathbf{I}'$, $\mathbf{J}'$, $\mathbf{K}'$ from \cite{MazSepTamb}. Moreover, in the proof of Lemma~\ref{lem: compatibility_omega_I} we showed that $\omega_{\mathbb{I}}$ restricts to $\omega_{\mathbf{I}}$ on the almost-Fuchsian locus, and this is sufficient to deduce that the metric $g$ and the other two symplectic forms $\omega_{\mathbb{J}}$ and $\omega_{\mathbb{K}}$ restrict to $\omega_{\mathbf{J}}$ and $\omega_{\mathbf{K}}$, thanks to formulas~\eqref{eq: g, omegaJ, omegaK}. This concludes the first item. Analogously, for the second item it suffices to prove that $\omega_{\mathbb{I}}$ restricts to the symplectic form $\omega_{\mathbf{I}'}$ on the locus $\mathcal{MGH}(S)$. To carry out the computation, similarly to Lemma~\ref{lem: compatibility_omega_I} we use the inclusion $\iota_{\mathbb{R}}$ introduced in diagram~\eqref{eq:comm diag MGH}. It suffices to evaluate
\[
(\pi_r^* \omega_G)\big((\dot\rho_1,\dot\rho_2),(\dot\rho_1',\dot\rho_2')\big)+(\pi_l^*\omega_G)\big((\dot\rho_1,\dot\rho_2),(\dot\rho_1',\dot\rho_2')\big) ,
\]
where $\dot{\rho}_i$ and $\dot{\rho}_i'$ are variations of representations into $\mathrm{PSL}(2,\mathbb{R})$ for $i=1,2$. 
Therefore, we obtain that \begin{align*}
(\pi_r^* \omega_G)\big((\dot\rho_1,\dot\rho_2),(\dot\rho_1',\dot\rho_2')\big)&+(\pi_l^*\omega_G)\big((\dot\rho_1,\dot\rho_2),(\dot\rho_1',\dot\rho_2')\big) = \\ &=\big(\omega_G^\R(\dot \rho_1,\dot\rho_1')+\omega_G^\R(\dot\rho_2,\dot\rho_2')\big) .
\end{align*} As a consequence of Theorem~\ref{thm: Thm C and Cor H} we deduce $\mathrm{hol}^*\big(\pi_r^*\omega_G+\pi_l^*\omega_G\big)|_{\mathcal{MGH}(S)}=\omega_{\mathbf I'}$.
\end{proof}

From what has been done so far, we have proved points (2) and (3) of Theorem~\ref{thm:main}, making essential use of the analytic continuation lemma for holomorphic tensors. As for point (1) of the same theorem, the same strategy cannot be applied, since the subspace to which we wish to restrict the complex hyperK\"ahler metric is not totally real, but rather a complex submanifold. For this reason, the next section is devoted to the proof of this remaining point.

\subsection{The complex bi-K\"ahler metric on \texorpdfstring{$\QF(S)$}{QF(S)}}

The main objective of this section is to prove point (1) of Theorem~\ref{thm:main}, namely the fact that the complex hyperK\"ahler metric introduced on \(\mathbb{M}_0(S)\) extends the already known complex bi-K\"ahler metric on the space of quasi-Fuchsian representations constructed by Loustau--Sanders \cite{loustau-sanders}. Point (1) needs a separate treatment from the other two. Indeed, the inclusion of $\QF(S)$ into $\mathbb M_0(S)$ is holomorphic, thus definitely not totally real. Hence, the analytic continuation arguments applied for $\mathcal{AF}(S)$ and $\mathcal{MGH}(S)$ does not apply in this setting.


In Section~\ref{sec: Bers metrics} we introduced the set of quasi-Fuchsian representations \(\mathcal{QF}(S)\) as a subset of the character variety \(\chi^{irr}\big(\pi_1(S),\PSL(2,\C)\big)\). Endowed with the natural complex structure induced by \(J_\chi\) on \(\chi^{irr}\big(\pi_1(S),\PSL(2,\C)\big)\) (see the beginning of subsection~\ref{subsec: pull backs from PSL(2,C)}), it is biholomorphic --- via the Bers parametrization map --- to the product space \(\mathcal{T}(S)\times\mathcal{T}(\overline{S})\), equipped with the (almost)-complex structure \(\mathbf{i}_{\mathcal{T}(S)}\times\mathbf{i}_{\mathcal{T}(\overline{S})}\) introduced in Section \ref{sec: beltrami differentials}.

In this setting, recalling from Subsection~\ref{subsubsec:tot geod} that quasi-Fuchsian representations arise as holonomies of minimal immersions into \(\mathrm{PSL}(2,\mathbb{C})\) whose image is contained in the totally geodesic subspace \(\mathbb{G}\), we have the following commutative diagram:

\begin{equation}\label{eq:comm diag TT}
\centering
\begin{tikzcd}
\mathcal{T}(S)\times\mathcal{T}(\overline{S}) \arrow[r, ] \arrow[d, ""'] & \mathbb{M}_0(S) \arrow[d, "\mathrm{hol}"] \\
\chi^{irr}\big(\pi_1(S),\PSL(2,\C)\big) \arrow[r, "\Delta"]                & \chi^*\big(\pi_1(S),\PSL(2,\C)^2\big) 
\end{tikzcd}
    \end{equation}
where the top arrow is, in accordance with Proposition~\ref{prop: dense open subset quasi-fuchsian},  given by
\([J_1,-J_2]\mapsto [J_1,0,-J_2,0,0]\), and we recall that $\Delta$ is the diagonal inclusion $\rho\mapsto (\rho,\rho)$. The left vertical arrow is the Bers parametrization map described above.

Loustau-Sanders defined a special metric structure on the quasi-Fuchsian space $\mathcal{QF}(S)$ using the fact that its natural complex structure can be interpreted as a complexification of the one coming from Teichm\"uller space (see \cite[\S~5.2]{loustau-sanders}). 


\begin{defn}\label{def:complex biKahler}
Let $(\mathbb X,\mathbf j)$ be a complex manifold with integrable almost-complex structure $\mathbf{j}:T\mathbb X\to T\mathbb X$. A \emph{complex bi-K\"ahler structure} is the data of two tensors $\mathcal{I}, \mathcal{P} \in \mathrm{End}_{\mathbb{C}}(T^{1,0}\mathbb X)$ and a complex Riemannian metric $g$ that satisfy the following properties:
\begin{itemize}
\item $\mathcal{I}^2 = -\mathcal{P}^2 = - \mathrm{Id}$ and $\mathcal{I}= - \mathbf{j} \circ \mathcal{P}$;
\item $\mathcal{I}$ and $\mathcal{P}$ are $\mathbf{j}$--holomorphic and integrable (that is, their Nijenhuis tensors $N_{\mathcal{I}}$ and $N_{\mathcal{P}}$ vanish identically on $T^{1,0} \mathbb{X}$);
\item the $(+1)$-- and $(-1)$--eigenspaces of $\mathcal{P}_x$ in $T^{1,0}_x \mathbb{X}$ have the same dimension for any $x \in \mathbb{X}$;
\item $g$ is $\mathbf{j}$--holomorphic and satisfies $g(\mathcal{P}\cdot,\mathcal{P}\cdot)=-g(\cdot,\cdot)$;
\item The pairing $\omega:=g(\mathcal{P}\cdot,\cdot)$ defines a symplectic form. \end{itemize}
\end{defn}
Following the same process outlined in Remark~\ref{rmk:extending_tensors}, there exist endomorphisms $\mathcal{I}_{\R}$ and $\mathcal{P}_{\R}$ on $T\mathbb{X}$ that commute with $\mathbf{j}$, and whose $\C$-linear extension restricts to $\mathcal I$ and $\mathcal P$ on $T^{1,0}\mathbb X$. It is not hard to check that $\mathcal{P}_{\R}$ provides $\mathbb{X}$ with a para-complex structure as in Definition~\ref{def:para-complex structure}. Moreover, by definition, the symplectic form
\(
\omega:=g(\mathcal{P}\cdot,\cdot) 
\)
is $\mathbf{j}$--holomorphic and satisfies $\omega(\mathcal{P}\cdot,\cdot)+\omega(\cdot,\mathcal{P}\cdot)=0$. Recall that the complex structures $\mathbf{i}_{\mathcal{T}(S)}$ and $\mathbf{i}_{\mathcal{T}(\overline{S})}$ on the Teichm\"uller spaces $\mathcal{T}(S)$ and $\mathcal{T}(\overline{S})$, respectively, are such that the map $[J]\mapsto [-J]$ is antiholomorphic. 

\begin{thm}[{\cite[Theorem 5.7 and 5.8]{loustau-sanders}}]\label{thm:loustau sanders}
The space $\mathcal{QF}(S)$ admits a $\mathrm{MCG}(S)$--invariant complex bi-K\"ahler metric $(g_{\mathcal{QF}},\mathbf{i}_{\mathcal{QF}},\mathbb{I}_{\mathcal{QF}},\mathcal P)$, where:
\begin{enumerate}
\item $\mathbf{i}_{\mathcal{QF}}$ is equal to the pullback of the endomorphism $\mathbf i_{\mathcal T(S)}\times\mathbf i_{\mathcal T(\overline S)}$ of $T(\mathcal{T}(S)\times \mathcal{T}(\overline{S}))$ via the Bers' double uniformization map;
\item $\mathbb{I}_{\mathcal{QF}}$ and $\mathcal{P}$ are equal to the pullback of the endomorphisms $\mathbf i_{\mathcal T(S)}\times\mathbf i_{\mathcal T(S)}$ and $(\mathrm{Id},-\mathrm{Id})$ of $T^{1,0}(\mathcal{T}(S)\times \mathcal{T}(\overline{S}))$, respectively, via the Bers' double uniformization map (where the $(1,0)$--part is taken with respect to $\mathbf{i}_{\mathcal{QF}}$);
\item The induced symplectic form $\omega_{\mathcal{QF}}$ coincides with the restriction on $\mathcal{QF}(S)$ of  Goldman's (complex) symplectic form on the character variety $\chi^{irr}\big(\pi_1(S),\PSL(2,\C)\big)$.
\end{enumerate}
\end{thm}

In accordance with point $(3)$ of the previous theorem, from now on we shall identify $\omega_{\mathcal{QF}}$ with the restriction of the Goldman symplectic form to \(\mathcal{QF}(S)\), denoted by \(\omega_G\). The holomorphic Riemannian metric $g_\mathcal{QF} = \omega_{G}(\cdot,\mathcal{P}\cdot)$ on $\mathcal{QF}(S)$ has been recently studied in \cite[\S~5]{ElE} along with some of its properties and its explicit description.
{\begin{remark}\label{rem: convention Loustau-Sanders}
In the original paper of Loustau--Sanders \cite{loustau-sanders}, the holomorphic Riemannian metric is defined as $\omega_G(\mathcal P\cdot, \cdot)$, hence it differs from \(g_{\mathcal{QF}}\) defined here by a sign. 
\end{remark}}

\begin{thm}\label{thm: restriction to quasi-Fuchsian}
The two complex structures $\mathbf i, \mathbb I$ and the symplectic form $\omega_\mathbb I$ defined on $\mathbb M_0(S)$ restrict to the triple $(\mathbf{i}_{\mathcal{QF}},\mathbb I_{\mathcal{QF}},2\omega_G)$ on $\mathcal{QF}(S)$. 
\end{thm}
\begin{proof}
By what was explained at the beginning of the section, the space of quasi-Fuchsian representations
\(\mathcal{QF}(S)\)
is diffeomorphic to a submanifold of \(\mathbb{M}_0(S)\), realized as \[
\mathcal{T}(S)\times\mathcal{T}(\overline S)=\{[J_1,q_1,-J_2,\overline{q_2},u]\in\mathbb{M}_0(S)\mid q_1=q_2=u=0\}.
\] By the definition of the complex structure \(\mathbf{i}\) on $\mathbb M(S)$, when restricted to
\(T^{1,0}\mathcal{T}(S)\times T^{1,0}\mathcal{T}(\overline{S})\),
it acts by multiplication by \(i\) on both the Beltrami differential \(\dot{\mu}_1\) and the conjugate Beltrami differential \(\dot{\overline{\mu}}_2\). In other words, it coincides with the product complex structure
\(\mathbf{i}_{\mathcal{T}(S)}\times\mathbf{i}_{\mathcal{T}(\overline{S})}\). 
By point (1) of Theorem~\ref{thm:loustau sanders}, this latter complex structure is biholomorphic to \(\mathbf{i}_{\mathcal{QF}}\) on $\mathcal{QF}(S)$. On the other hand, the endomorphism $\mathbb{I}$ defined in~\eqref{eq: def I in trivialization of pi} and~\eqref{eq: Complex structure I} restricts to \(T^{1,0}\mathcal{T}(S)\times T^{1,0}\mathcal{T}(\overline{S})\) and it acts as the multiplication by \(i\) on the Beltrami differential \(\dot{\mu}_1\) and by \(-i\) on the conjugate Beltrami differential \(\dot{\overline{\mu}}_2\). According to Theorem~\ref{thm:loustau sanders} point (2), its restriction coincides with $\mathbb{I}_{\mathcal{QF}}$ on $\mathcal{QF}(S)$. The final statement concerning the symplectic form $\omega_{\mathbb{I}}$ follows from its definition in~\eqref{eq: omega_mathbb I}, by a computation entirely analogous to that carried out in the proof of Lemma~\ref{lem: compatibility_omega_I}, where in this case the inclusion at the level of representations is the diagonal embedding.
\end{proof}

The following Corollary \ref{cor: restriction of complex hyperkahler to QF(S)} completes the proof of the remaining part of Theorem~\ref{thm:main}, namely item (1).

\begin{cor}\label{cor: restriction of complex hyperkahler to QF(S)}
The complex hyperK\"ahler metric on $\mathbb M_0(S)$ extends the Loustau-Sanders complex bi-K\"ahler metric on $\mathcal{QF}(S)$, up to some constants. More precisely, the holomorphic Riemannian metric \(g\) restricts to $-2ig_{\mathcal{QF}}$.
\end{cor}
\begin{proof}
The complex bi-K\"ahler structure on $\mathcal{QF}(S)$ (Theorem \ref{thm:loustau sanders}) is completely determined by the triple $(\mathbf{i}_{\mathcal{QF}},\mathbb{I}_{\mathcal{QF}},\omega_G)$, since $\mathcal{P}=\mathbf{i}_{\mathcal{QF}}\circ \mathbb{I}_{\mathcal{QF}}$ and $g_{\mathcal{QF}}=\omega_G(\mathcal{P}\cdot,\cdot)$. In particular, 
\begin{align*}
    g|_{\mathcal{QF}(S)}&=\omega_{\mathbb I}(\cdot,\mathbb I\cdot)|_{\mathcal{QF}(S)}  \\ &=2\omega_G(\cdot,\mathbb I_{\mathcal{QF}}\cdot) \tag{Theorem \ref{thm: restriction to quasi-Fuchsian}}  \\ &=-2\omega_G(\cdot,\mathbf{i}_{\mathcal{QF}}\circ\mathcal P\cdot) \tag{Definition \ref{def:complex biKahler}}  \\ &=-2i\omega_G(\cdot,\mathcal P\cdot) \tag{$\omega_G$ is $\mathbf{i}_{\mathcal{QF}}$-holomorphic} \\ &=-2ig_{\mathcal{QF}}.
\end{align*}
\end{proof}

\section{The infinitesimal \texorpdfstring{$\mathbb C^*$}{C*} action and totally geodesic submanifolds}\label{sec:C*}

In this section, we prove further properties of the complex hyperKähler structure on $\mathbb M_0(S)$. First, we discuss the existence, at least in a local sense, of a $\mathbb C^*$ action on $\mathbb M(S)$, and see that some objects of the complex hyperKähler structure are preserved by it. We then observe that the area functional descends to a holomorphic map on $\mathbb M(S)$, and, extending analogue results for $\mathcal{AF}(S)$ and $\mathcal{MGH}(S)$, its differential can be seen as the contraction of $\omega_{\mathbb I}$ with the infinitesimal generator of the $\C^*$-action. Finally, we show that $\mathcal{AF}(S), \mathcal{MGH}(S),$ and $\QF(S)$ are all totally geodesic for the holomorphic Riemannian metric $g$ on $\mathbb M(S)$.

\subsection{The infinitesimal \texorpdfstring{$\C^{*}$}{C*}--action}
Recall that the spaces $\mathrm M(S)$ and $\mathrm M^*(S)$ have been defined in Section \ref{subsec:regular min imm}.

 Consider the \(\mathbb C^*\)-action on $\mathrm M(S)$ given by: \begin{equation}\label{eq: C* action}\begin{aligned}\mathbb C^*\times \mathrm M(S)&\longrightarrow\mathrm M(S), \\\big(\zeta,(J_1,q_1,-J_2,\overline{q_2},u)\big)&\longmapsto
(J_1,\zeta q_1,-J_2,\zeta^{-1}\overline{q_2},u),\end{aligned}\end{equation}
which is well-defined since it preserves the differential equation satisfied by \(u\) (see Equation~\eqref{eq: Gauss semplif in u}).  
Observe that this $\mathbb C^*$-action also preserves the property of being infinitesimaly rigid (see Definition~\ref{def:regular immersion}), since Equation \eqref{eq: linearized Gauss} only depends on $q_1$ and $\overline{q_2}$ through their product. 

Unfortunately, it does not seem clear whether the action of $\mathbb C^*$ preserves the condition of having holonomy within $\chi^*(\pi_1(S), \PSL(2,\C)^2)$, hence we cannot say that $\mathbb C^*$ acts on $\mathrm M^*(S)$. Nevertheless, since $\mathrm M^*(S)$ is an open subset of $\mathrm M(S)$, we can say the following.

\begin{lem}
\label{lem: C* su mathbb M}
For every $\gamma_0$ in $\mathrm M^*(S)$, there exist open subsets $U\subset \mathrm M^*(S)$ and $V\subset \mathbb C^*$ with $\gamma_0\in U$ and $1\in V$ such that, if $\gamma_1, \gamma_2\in U$ are identified in $\mathbb M(S)$ and $\zeta\in V$, then $\zeta\cdot \gamma_1$ and $\zeta\cdot \gamma_2$ both lie in $\mathrm M^*(S)$ and are identified in $\mathbb M(S)$. 

Denoting by $[U]$ the image of $U$ through the projection $\mathrm M^*(S)\to \mathbb M(S)$, the induced maps 
\begin{align*}
    V\times U\to \mathrm M^*(S) &\qquad (\zeta, \gamma)\mapsto \zeta\cdot \gamma\\
    V\times [U]\to \mathbb M(S) & \qquad (\zeta, [\gamma])\mapsto [\zeta\cdot \gamma]
\end{align*}
are both holomorphic.
\end{lem}
\begin{proof}
Recall the commutative diagram in Equation \eqref{eq: diagramma proiezioni} and that $\pi$ and $\pi_0$ are local biholomorphisms. The $\mathbb C^*$ action on $\mathrm M(S)$ commutes through the projection $\mathrm M(S)\to \mathrm{QD}(S)\times \mathrm{QD}(\overline S)$ with the $\C^*$-action on $\mathrm{QD}(S)\times \mathrm{QD}(\overline S)$ given by $\zeta\cdot (c_1, \overline{c_2}, q_1, \overline{q_2})= (c_1, \overline{c_2}, \zeta q_1, \zeta^{-1}\overline{q_2})$; this action is trivially holomorphic and descends to a well-defined holomorphic action of $\mathbb C^*$ on $\mathbb{QD}(S)\times\mathbb{QD}(\overline S)$. 

Now, since topologically we have $\mathrm{QD}(S)\times\mathrm{QD}(\overline S)\cong\mathbb{QD}(S)\times \mathbb{QD}(\overline S)\times \Diff_0(S)\times \Diff_0(S)$, every element in $\mathrm{QD}(S)\times\mathrm{QD}(\overline S)$ admits a complete system of neighborhoods such that the intersection with every $\Diff_0(S)\times \Diff_0(S)$-orbit is connected. Using this remark, together with the facts that $\pi$ and $\pi_0$ are local homeomorphisms and that the $\mathbb C^*$-action on $\mathrm M(S)$ is continuous, we can say that, for all element $\gamma_0\in \mathrm M^*(S)$, there exist $U$ and $U_0$ neighborhoods of $\gamma$ and $V\subset \mathbb C^*$ neighborhood of $1$ such that: \begin{itemize}
\item for all $\zeta\in V$ and $\gamma\in U$, $\zeta\cdot \gamma\in U_0$;
\item
the restrictions of $\pi_0$ on $U$ and $U_0$ and the restrictions of $\pi$ on $[U]$ and $[U_0]$ are all injective, hence biholomorphisms with their respective images;
    \item the intersection of $U$ with every fiber of the projection to $\mathrm M^*(S)\to \mathbb M(S)$ is connected.
\end{itemize}

Now, let $F$ be the intersection of $U$ with a fiber of $\mathrm M^*(S)\to \mathbb M(S)$. By construction of $\mathbb M(S)$, $[\pi_0(F)]\in \mathbb{QD}(S)\times\mathbb{QD}(\overline S)$ is a point.  The assumption on $U$ implies that, for all $\zeta\in V$, $\zeta\cdot F\subset U_0$ is connected, and it is mapped through $[\pi_0]$ to the element $\zeta\cdot [\pi_0(F)]$. We therefore conclude that $[\zeta\cdot F]$ is one point in $\mathbb M(S)$. 

From the commutativity of the diagram
 \begin{center}
\label{eq: diagramma restrizioni}
   \begin{tikzcd}
 V\times U\arrow[r, "\mathit{action}"] \arrow[d, "\mathit{id}\times \pi_0"]  & U_0 \arrow[d, "\pi_0"]  \\
V\times \pi_0(U) \arrow[r, "\ \mathit{action}\ "] & \pi_0(U_0)
\end{tikzcd}
\end{center}
we get that the top horizontal map is holomorphic as a consequence of the holomorphicity of the bottom one and of the fact that the vertical maps are biholomorphisms. Similarly for $V\times[U]\to [U_0]$.
\end{proof}

While we might not have a global $\C^*$-action on $\mathbb M(S)$, Lemma \ref{lem: C* su mathbb M} allows to introduce the notion of ``infinitesimal generator" of the action, as explained in the following proposition.

\begin{prop}
\label{prop: infintiesimal generator well-def}
    The derivative
    \[
   {\xi}|_{[\gamma]}:= \frac{\mathrm{d}}{\mathrm{d}z}\bigg|_{z=0} [e^{iz}\cdot \gamma]
    \]
provides a well-defined $\mathbf i$-holomorphic vector field on $\mathbb M(S)$, hence an element $\xi\in H^0(T^{1,0}(\mathbb M(S))$.
\end{prop}
\begin{proof}
    By Lemma \ref{lem: C* su mathbb M}, for every element $\gamma_0\in \mathrm M_0(S)$ there exists a neigbourhood $U$ such that, for all $\gamma\in U$ with $[\gamma]=[\gamma_0]$, the maps $z\mapsto [e^{iz}\cdot \gamma]$ and $z\mapsto [e^{iz}\cdot \gamma_0]$, defined for sufficiently small $z\in \C$, coincide and are holomorphic. Hence, they determine the same vector 
    $$
    \frac{\mathrm{d}}{\mathrm{d}z}\bigg|_{z=0} [e^{iz}\cdot \gamma]\in T_{[\gamma]}^{1,0}\mathbb M(S)~.
    $$
    This defines a locally constant map from the fiber over $[\gamma]$ of $\mathrm M^*(S)\to \mathbb M(S)$ to $T_{[\gamma]}^{1,0}\mathbb M(S)$. Since the fibers are connected, we conclude that the map is well defined. 

    Finally, the holomorphicity in Lemma \ref{lem: C* su mathbb M} implies the holomorphicity of the vector field $\xi$.
\end{proof}

The existence of the infinitesimal generator of the $\C^*$-action allows to talk about infinitesimal variations of holomorphic tensors on $\mathbb M(S)$ through the $\mathbb C^*$-action, as we explain now. 

With a little abuse of notation, for all $\zeta\in \mathbb C^*$ we denote by $\zeta$ also the induced map on $\mathrm M(S)$, and use $\zeta^*$ to generically denote the pull-back of tensors throught this map.

Let ${\upalpha}$ be a holomorphic tensor on $\mathbb M(S)$. 
As explained in the proof of Proposition \ref{prop: infintiesimal generator well-def}, the map $z\mapsto [e^{iz}\gamma_0]$, defined for sufficiently small $z$, does not change if we replace $\gamma_0$ with other nearby representatives of $[\gamma_0]$; moreover, the map $(z, [\gamma])\mapsto [e^{iz}\gamma]$ defines a holomorphic function for $z$ in a neighborhood of $1$ and $[\gamma]$ in a neighborhood of $[\gamma_0]$; we still denote such local map on $\mathbb M(S)$ by $e^{iz}$. As a consequence, we can define the holomorphic tensor 
\begin{equation}
\label{eq: infinitesimal deformation of tensors through C*}
\xi\cdot \upalpha := \frac{\mathrm{d}}{\mathrm{d}z} \bigg|_{z=0}\left(  (e^{-iz})^*\upalpha \right)
\end{equation}
on $\mathbb M(S)$. By identifying canonically the holomorphic vector field $\xi$ with its real part in the smooth tangent bundle, we can see $\xi\cdot \upalpha$ as the Lie derivative of $\upalpha$ in the direction corresponding to $\xi$.

\begin{remark}
    \label{rmk: definition of infinitesimal tensors under C*}
Observe that $\xi\cdot \upalpha=0$ if and only if $(e^{-iz})^*\upalpha\equiv \upalpha$ on sufficiently small neighborhoods of any point in $\mathbb M(S)$ and for sufficiently small $|z|$. Indeed, in a local description, $\xi\cdot \upalpha=0$ implies that, for $z_0$ sufficiently small,
\[
\frac{\mathrm{d}}{\mathrm{d}z}\bigg|_{z=z_0}  (e^{-{iz}})^*\upalpha = \frac{\mathrm{d}}{\mathrm{d}z}\bigg|_{z=z_0}  (e^{-iz_0})^* (e^{{-iz+iz_0}})^*\upalpha = (e^{-iz_0})^* (\xi\cdot \upalpha)=0,
\]
which means that, for small $z$, $(e^{-iz})^*\upalpha$ does not depend on $z$, hence $(e^{-iz})^*\upalpha\equiv \upalpha$.
\end{remark}

We are now able to discuss the interplay between the infinitesimal $\C^*$-action and some elements of the complex hyperKähler structure on $\mathbb M_0(S)$.

First, let us observe that the $\C^*$ action extends the $S^1$ action on \(\mathcal{AF}(S)\) and \(\mathcal{MGH}(S)\) defined in \cite{donaldson2003, trautwein2018infinite} for \(\mathcal{AF}(S)\) and in \cite{MazSepTamb} for \(\mathcal{MGH}(S)\). Indeed, it is straightforward to check that for all $\theta\in \R$, the action times $e^{i\theta}$ preserves both $\mathcal{AF}(S)$ and $\mathcal{MGH}(S).$ We use the notation as above and denote by $e^{i\theta}$ also the induced map from each of these spaces to itself.

\begin{thm}[{\cite[Lemma 4.5.9]{trautwein2018infinite}, \cite[Theorem E]{MazSepTamb}}]\label{thm: circle action MGH}
The circle action on \(\mathcal{AF}(S)\) satisfies: \[\label{eq: relation for the circle action on AF} ({e^{i\theta}})^*\mathbf g=\mathbf g,\quad ({e^{i\theta}})^*\mathbf I=\mathbf I,\quad ({e^{i\theta}})^*\omega_{\mathbf I}=\omega_\mathbf {I}. 
\]
Similarly, the circle action on $\mathcal{MGH}(S)$ satisfies:
 \[\label{eq: relation for the circle action on MGH} ({e^{i\theta}})^*\mathbf g'=\mathbf g',\quad ({e^{i\theta}})^*\mathbf I'=\mathbf I',\quad ({e^{i\theta}})^*\omega_{\mathbf I'}=\omega_\mathbf {I'}. 
\]

\end{thm}

These results extend in our setting in the following sense.

\begin{thm}\label{thm: C* action on M0(S)}
Recall the tensors $g, \mathbb I,$ and $\omega_{\mathbb I}$ from the complex hyperKähler structure of $\mathbb M_0(S)$. Then,
\[
\xi\cdot g =0, \qquad\qquad \xi\cdot \mathbb I= 0, \qquad\qquad \xi\cdot \omega_{\mathbb I}=0.
\]
\end{thm}
\begin{proof}
Let $\gamma \in \mathrm{MGH}(S)$ (see Equation \eqref{eq:inclusion MGH in M}).
By Theorem \ref{thm: circle action MGH} and the third item of Theorem \ref{thm:main}, we know that, at $[\gamma]$ we have the identities 
\[\xi\cdot g |_{[\gamma]}= \frac{\mathrm{d}}{\mathrm{d}z} \bigg|_{z=0} (e^{iz})^* g = \frac{\mathrm{d}}{\mathrm{d}\theta}\bigg|_{\theta=0} (e^{i\theta})^*g= \frac{\mathrm{d}}{\mathrm{d}\theta}\bigg|_{\theta=0} (e^{i\theta})^*\mathbf g'= 0
\]
and similar for $\mathbb I$ and $\omega_{\mathbb I}$. Since they are holomorphic tensors that vanish along $\mathcal{MGH}(S)$, a totally real submanifold of maximal dimension, they vanish on the whole $\mathbb M_0(S)$ (Lemma \ref{lem: fundamental lem analytic continuation} and Remark \ref{rmk: fundamental_lemma_applied to tensors}).
\end{proof}

\subsection{The interplay with the area functional}

Let us now show that also the following statement can be extended in our context. 

\begin{thm}[{\cite[Theorem 4.6.9]{trautwein2018infinite},\cite[Theorem E]{MazSepTamb}}]\label{thm: hamiltonian MGH}
The circle action on \(\mathcal{AF}(S)\) (resp. \(\mathcal{MGH}(S)\)) is Hamiltonian with respect to \(\omega_{\mathbf I}\) (resp. \(\omega_{\mathbf I'}\)), with Hamiltonian function given by the area functional, assigning to each equivariant minimal (resp. maximal) immersion the area of the first fundamental form, seen as a Riemannian metric on $S$.
\end{thm}
In order to extend this result, we first need to introduce an appropriate notion of area functional on $\mathbb M_0(S)$, that extends the area functional on \(\mathrm{AF}(S)\) and \(\mathrm{MGH}(S)\).
Define the complex area functional by  
\[ 
\hat{\mathcal A}\colon \mathrm M^*(S)\longrightarrow \mathbb C, \qquad \hat{\mathcal A}(\gamma):=\int_S \mathrm{d}A_{\Is_\gamma}, 
\] 
where \(\Is_\gamma\) denotes the first fundamental form associated with \(\gamma\), and \( \mathrm{d}A_{\Is_\gamma} = \Is_\gamma(\mathcal J_\gamma\,\cdot,\cdot) \) is its complex area form, with \(\mathcal J_\gamma\) being the induced bicomplex structure.
\begin{lem}\label{lem: Area functional on M*(S) is holomorphic}
The complex area functional \(\hat{\mathcal A}:\mathrm{M}^*(S)\to\mathbb C\) is holomorphic.
\end{lem}

\begin{proof}
We will prove the stronger statement that the map
\[
    \mathrm{M}^*(S) \longrightarrow \Omega^2(S,\mathbb C),
    \qquad
    \gamma \longmapsto \mathrm{d}A_{\mathtt{I}_\gamma},
\]
is holomorphic.

We first consider finite-dimensional families. Let $\Lambda$ be a
finite-dimensional smooth manifold and let
$\lambda \mapsto \gamma_\lambda \in \mathrm{M}^*(S)$ be a smooth map.
Set
\[
    \mathtt{I}_\lambda := \mathtt{I}_{\gamma_\lambda},
    \qquad
    \mathcal J_\lambda := \mathcal J_{\gamma_\lambda},
    \qquad
    \omega_\lambda := \mathrm{d}A_{\mathtt{I}_\lambda}
    = \mathtt{I}_\lambda(\mathcal J_\lambda\,\cdot,\cdot).
\]
By Lemma \ref{lem: from M* to I II}, the tensor $\mathtt{I}_\lambda$ depends smoothly on $(\lambda,p)\in\Lambda\times S$.

The same is true of $\mathcal J_\lambda$ and $\omega_\lambda$.
Indeed, the two isotropic lines of $\mathtt{I}_\lambda$ are distinct
roots of the quadratic equation
$\mathtt{I}_\lambda(v,v)=0$ in the projective tangent space, and
therefore vary smoothly with $(\lambda,p)$.
Since $\mathtt{I}_\lambda$ is positive, these lines lie in distinct
components of $\mathbb P_{\mathbb C}(\mathbb C T_pS) \setminus \mathbb P_{\mathbb R}(T_pS)$, and the orientation of $S$ consistently distinguishes them. Since these lines depend smoothly on $\lambda$, also the tensor $\mathcal J_{\lambda}$, defined as $i$ on one line and $-i$ on the other one, depends smoothly on $\lambda$. 

Suppose now that $\Lambda$ is a complex manifold and that
$\lambda\mapsto\gamma_\lambda$ is holomorphic.
Fix a complex structure $J_0$ compatible with the orientation
of $S$, and let $z$ be a local $J_0$-holomorphic coordinate
on an open subset $U\subset S$. 
Keeping this coordinate fixed as $\lambda$ varies, write
\[
    \mathtt{I}_\lambda
    =
    E_\lambda\,\mathrm{d}z\otimes \mathrm{d}z
    +F_\lambda
       \bigl(\mathrm{d}z\otimes \mathrm{d}\bar z+\mathrm{d}\bar z\otimes \mathrm{d}z\bigr)
    +G_\lambda\,\mathrm{d}\bar z\otimes \mathrm{d}\bar z,
    \qquad
    \omega_\lambda=a_\lambda\,\mathrm{d}z\wedge \mathrm{d}\bar z.
\]
By Lemma \ref{lem: from M* to I II}, the coefficients $E_\lambda,F_\lambda,G_\lambda$
are holomorphic in $\lambda$ at every point of $U$, since
the frame $(\partial_z,\partial_{\bar z})$ is independent
of $\lambda$. By a standard computation, we have $    a_\lambda^2=E_\lambda G_\lambda-F_\lambda^2$,
thus the area form $\omega_\lambda$ depends holomorphically on $\lambda$ at each point of $S$, as well. Together with smoothness in $(\lambda,p)$, this shows that
\[
    \Lambda\longrightarrow\Omega^2(S,\mathbb C),
    \qquad
    \lambda\longmapsto\omega_\lambda,
\]
is holomorphic.

We now pass to the Fr\'echet setting.
The argument above, applied to smooth curves,
shows that the map
$\gamma\mapsto \mathrm{d}A_{\mathtt{I}_\gamma}$
maps smooth curves to smooth curves in $\Omega^2(S,\mathbb C)$.
It is therefore smooth by Theorem \ref{thm: boman}.
Moreover, its restriction to every holomorphic disk is holomorphic
by the preceding computation, thus its differential is
$\mathbb C$-linear.
Hence $\gamma\mapsto \mathrm{d}A_{\mathtt{I}_\gamma}$ is holomorphic.

Finally, integration over $S$ defines a continuous
$\mathbb C$-linear map $\Omega^2(S,\mathbb C)\rightarrow\mathbb C$,
therefore
\[
    \widehat{\mathcal A}(\gamma)
    =\int_S \mathrm{d}A_{\mathtt{I}_\gamma}
\]
is holomorphic.
\end{proof}

\begin{prop}\label{prop:area functional well defined}
The complex area functional descends to the complex manifold \(\mathbb M(S)\), thereby defining a holomorphic function
\[
\mathcal A:\mathbb M(S)\longrightarrow\mathbb C.
\]
\end{prop}
\begin{proof}
First, observe that \(\hat{\mathcal A}\) is invariant under the action of \(\Diff_0(S)\).  Indeed, for every \(\varphi\in\Diff_0(S)\), one has \( \Is_{\varphi^*\gamma}=\varphi^*\Is_\gamma \) and consequently \( \mathrm{d}A_{\Is_{\varphi^*\gamma}} = \varphi^*\mathrm{d}A_{\Is_\gamma}. \) Given that \(\varphi\) is orientation-preserving, it follows that \( \hat{\mathcal A}(\varphi^*\gamma)= \hat{\mathcal A}(\gamma).\) As a consequence of the \(\mathrm{Diff}_0(S)\)-invariance and Lemma~\ref{lem: invariance of map under Diff-action}, the functional \(\hat{\mathcal A}\) descends to a well-defined holomorphic function \(\mathcal A:\mathbb M(S)\longrightarrow\mathbb C.\) 
\end{proof}

Despite the $\C^*$-action is not properly defined on $\mathbb M_0(S)$, we are now able to show that Theorem \ref{thm: hamiltonian MGH} extends nicely to $\mathbb M_0(S)$.

\begin{thm}
The following identity holds on $\mathbb M_0(S)$:
\[
\mathrm{d}\mathcal A= \omega_{\mathbb I}(\xi, \cdot) 
\]

\end{thm}
\begin{proof}

Both terms in the equality are holomorphic 1-forms, so it suffices to show that they coincide along the totally real submanifold $\mathcal{MGH}(S)$. This is indeed the case: along $\mathcal{MGH}(S)$, $\mathcal A$ restricts to the usual area functional, $\xi$ corresponds to $\frac{\mathrm{d}}{\mathrm{d}\theta} e^{i\theta}$, and $\omega_{\mathbb I}$ restricts to $\omega_{\mathbf I'}$.
\end{proof}

\subsection{ \texorpdfstring{$\mathcal{QF}(S)$}{QF(S)}, 
\texorpdfstring{$\mathcal{AF}(S)$}{AF(S)}, and \texorpdfstring{$\mathcal{MGH}(S)$}{MGH(S)} are totally geodesic} 

In this subsection, we prove that three recurring loci of $\mathbb M_0(S)$, namely $\mathcal{AF}(S)$, $\mathcal{MGH}(S)$, and $\QF(S)$, are totally geodesic for the holomorphic Riemannian metric $g$ on $\mathbb M_0(S)$. 

The fact that the spaces $\mathcal{AF}(S)$ and $\mathcal{MGH}(S)$ are totally geodesic follows from a very general result about totally real submanifolds in holomorphic manifolds with non-degenerate induced Riemannian metric.

\begin{prop}\label{prop: totally geodesic holomorphic metric}
Let $(\mathbb{X}, \mathbf{i})$ be a complex manifold of complex dimension $n$, equipped with a holomorphic Riemannian metric $g$. Let $N\subset \mathbb{X}$ be a maximal totally real submanifold and assume that the restriction $g|_{TN}$
is real-valued and non-degenerate. Then $N$ is a totally geodesic sumbanifold of $\mathbb{X}$.
\end{prop}
\begin{proof}
Denote $h=g|_{TN}$. Let $D$ denote the holomorphic Levi-Civita connection of $g$. Recall that, according to Definition~\ref{def:homorphic levi-civita}, $D$ is the $\mathbb{C}$-linear extension of an affine connection of $\mathbb{X}$. In particular, $D_V$ keeps the real tangent bundle $T \mathbb{X}$ invariant whenever $V \in \Gamma(T \mathbb{X})$. As in the classical Riemannian setting, the submanifold $N$ is totally geodesic if and only if, for any pair of tangent vector fields $X, Y \in \Gamma(T N)$, $D_X Y$ is tangent to $N$. Since $N$ is maximal totally real, for every $p\in N$ one has a direct sum decomposition $T_p\mathbb{X}=T_pN\oplus \mathbf{i}(T_pN)$ as real vector spaces (see Definition \ref{def: totally real submanifold}).
Let $X,Y,Z$ be smooth vector fields tangent to $N$. The Koszul formula gives
\[
\begin{aligned}
2g(D_XY,Z) = & \, X\bigl(g(Y,Z)\bigr) + Y\bigl(g(X,Z)\bigr) - Z\bigl(g(X,Y)\bigr)\\
& + g([X,Y],Z) - g([X,Z],Y) - g([Y,Z],X) .
\end{aligned}
\]
Every term on the right-hand side is real along $N$. Hence $g(D_XY,Z)\in\mathbb R$ for every vector field $Z$ tangent to $N$.
Fix $p\in N$ and write $D_XY=A+\mathbf{i}B$, for some $A,B\in T_pN$. Since $g$ is complex bilinear,
\[
g(D_XY,Z)=h(A,Z)+i\,h(B,Z)
\]
for every $Z\in T_pN$. The left-hand side is real, so $h(B,Z)=0$
for every $Z\in T_pN$. Since by assumption $h$ is non-degenerate, we conclude that $B=0$ and, hence, $D_XY\in T_pN$. Since $p \in N$ was arbitrary, we conclude that $D_X Y$ is tangent to $N$, as desired.
\end{proof}
\begin{cor}
The submanifolds $\mathcal{AF}(S)$ and $\mathcal{MGH}(S)$ are both totally geodesic in $(\mathbb M_0(S),g)$. Since their intersection coincides with the Teichmüller
space $\mathcal T(S)$, the latter is also totally geodesic in $(\mathbb M_0(S),g)$.
\end{cor}
\begin{proof}
By Corollary~\ref{cor: restriction of complex hyperkahler to AF e MSH}, the holomorphic metric \(g\) restricts to \(\mathbf g\) on \(\mathcal{AF}(S)\) and to \(\mathbf g'\) on \(\mathcal{MGH}(S)\), and both restricted metrics are real and non-degenerate. Moreover, Proposition~\ref{prop: AF and MGH are totally real submanifolds} shows that both spaces are totally real submanifolds of maximal dimension. The first part of the statement therefore follows from Proposition~\ref{prop: totally geodesic holomorphic metric}. As for the second part, Theorem~\ref{thm:hyperkahler_almost_fuchsian} implies, for instance, that \(\mathcal{T}(S)\) is totally geodesic in \(\mathcal{AF}(S)\). The desired conclusion then follows.
\end{proof}

To prove that $\mathcal{QF}(S)$ is totally geodesic inside $\mathbb M_0(S)$ we will need a different approach. In fact, we will use again the properties of the infinitesimal $\C^*$ action.

\begin{cor}
The quasi-Fuchsian space \(\mathcal{QF}(S)\) is a totally geodesic submanifold of \(\big(\mathbb M_0(S),g\big)\). 
\end{cor}
\begin{proof}
Consider any point $[\gamma_0]$ in $\mathcal{QF}(S)$ and take a neighborhood $W_0$ inside $\mathbb M_0(S)$ over which the projection $\pi$ to $\mathbb{QD}(S)\times\mathbb{QD}(\overline S)$ gives a homeomorphism with the image. Up to shrinking it, we can assume that an element $[\gamma]\in W_0$ lies in $\mathcal {QF}(S)$ if and only if $\pi([\gamma])$ lies in the zero section $\TSTS\subset\mathbb{QD}(S)\times\mathbb{QD}(\overline S)$.

    By Lemma \ref{lem: C* su mathbb M}, there exists $\eps>0$ and a smaller neighborhood $W\subset W_0$ of $[\gamma]$ such that, for all $|z|<\eps$, the action times $e^{iz}$ can be seen as a map from $W$ to $W_0$.
    By Theorem \ref{thm: C* action on M0(S)} and by Remark \ref{rmk: definition of infinitesimal tensors under C*}, we have that $(e^{-iz})^*g=g$ on $W$, hence every $e^{-iz}$ determines locally an isometry with the image. On the other hand, since $\pi$ conjugates the $\C^*$-actions, an element $[\gamma]$ in $W$ is fixed by $e^{-iz}$ for all $|z|<\eps$ if and only if $\pi([\gamma])\in \TSTS$, hence if and only if $[\gamma]\in W_0\cap \QF(S)$.
By a standard argument (e.g. see  \cite[Chapter~II, Theorem~5.1]{kobayashi2012transformation}), every geodesic tangent to $\QF(S)$ starting at $[\gamma_0]$ is then preserved by the action times each $e^{iz}$, hence it is contained inside $\QF(S)$, proving that $\QF(S)$ is totally geodesic.
\end{proof}

\begin{remark}\label{rmk:holo involution}
   Similarly to the proof that $\mathcal{QF}(S)$ is totally geodesic as a fixed point of an isometric, holomorphic action, also  $\mathcal{AF}(S)$ and $\mathcal{MGH}(S)$ can be realized as fixed loci of a local \emph{antiholomorphic} involution. To see this, observe that the involutions of $\mathrm M(S)$ given by
     \begin{align*}
       \nu\colon (J_1, q_1, -J_2, \overline{q_2}, u)&\mapsto (J_2, q_2, -J_1, \overline{q_1}, \overline u)\qquad \text{and}\\
       \nu'\colon (J_1, q_1, -J_2, \overline{q_2}, u)&\mapsto (J_2, -q_2, -J_1, -\overline{q_1}, \overline u) \ ,
   \end{align*}   
 fix pointwise $\mathrm{AF}(S)$ and $\mathrm{MGH}(S)$ respectively. 
These involutions descend to suitable neighborhoods of $\mathcal{AF}(S)$ and $\mathcal{MGH}(S)$ in $\mathbb M_0(S)$ that project injectively to $\mathbb{QD}(S)\times\mathbb{QD}(\overline S)$, over which they correspond to $ [J_1, q_1, -J_2, \overline{q_2}]\mapsto [J_2, \pm q_2, -J_1, \pm \overline{q_1}]$, hence they are antiholomorphic involutions. By an analytic continuation argument, using that the restriction of $g$ to $\mathcal{AF}(S)$ and $\mathcal{MGH}(S)$ is real, one can show that these involutions pull-back the metric $g$ to its complex conjugate, hence they are isometries for the real part of $g$; as a consequence, one could deduce by an alternative argument that $\mathcal{AF}(S)$ and $\mathcal{MGH}(S)$ are totally geodesic.

We actually expect $\nu$ and $\nu'$ to descend to global anti-holomorphic maps on $\mathbb M(S)$. However, although it is easily verified that $\nu$ and $\nu'$ preserve the infinitesimal rigidity, it is less evident that they preserve  the condition for which the holonomy lies in $\chi^*(\pi_1(S), \PSL(2,\C)^2)$.

\end{remark}

\section{A geometric interpretation of \texorpdfstring{$\mathbb{K}$}{K}} \label{sec:geometric interpretation of K}

In this section we provide a more geometric description of $\mathbb K$ using an approach analog to \cite[\S~5.3]{MazSepTamb}. In the Anti-de Sitter setting, Krasnov-Schlenker (\cite{kirill07}) introduced a geometric construction yielding a diffeomorphism
\(
\mathcal{C}:\mathcal{MGH}(S)\to \mathcal{T}(S)\times\mathcal{T}(S)
\)
which differs from the parametrization using the holonomy map. More precisely, given a point $[\gamma]\in\mathcal{MGH}(S)$, the associated pair of hyperbolic metrics is obtained as the first fundamental forms of two equivariant immersions into $\Ads^3$ having constant Gaussian curvature and canonically determined by the geometry of $[\gamma]$. This construction was subsequently used to prove the identity $\mathbf{K}'=\mathcal{C}^*(\mathrm{Id},-\mathrm{Id})$ (\cite[Theorem D]{MazSepTamb}). We now show that an analogous construction can be carried out in the present setting, after restricting to a suitable open subset of $\mathbb{M}_0(S)$.
\\ \\
Let $\gamma\in \mathrm M(S)$, denote its immersion data by $\Is$ and $\II=\Is(B\cdot, \cdot)$, and the normal vector field by $N$. Consider the two maps $\sigma'_{\pm}\colon \widetilde S\to \PSL(2,\C)$ as $\sigma'_{\pm}(x)= \exp_{\sigma(x)}( \pm i\frac \pi 4 N(x))$, where $\sigma$ is the immersion associated with $\gamma$. An explicit computation (see \cite{BEE}), shows that, taking local lifts of $\sigma$ in $\SL(2,\C)$, one has 
\[
\exp_{\sigma(x)}( z N)= \cosh(z)\sigma +\sinh(z) N, 
\]
hence 
\[
\sigma_{\pm}'(x)= \frac{\sqrt 2} 2  \sigma(x) \pm i \frac{\sqrt{2}}2 N(x).
\]

\begin{remark}\label{rmk:constant curvature surfaces}
    The above definition is inspired by the  construction of equidistant surfaces in the Anti-de Sitter space (\cite{kirill07}, see also \cite{barbot2011, bonsante-seppi2018}). To explain this, assume here that $\sigma$ is an immersion in $\Ads^3\cong \PSL(2,\R)$, seen as a totally geodesic subspace of $\PSL(2,\C)$. Then $\sigma$ has a  normal vector field in $\Ads^3$ of norm identically $-1$. The normal vector field $N$ in $\PSL(2,\C)$ is defined by the condition $\langle N, N\rangle =1$, so it equals $\pm i$ times the $\Ads^3$ normal vector field. Therefore $\sigma_{\pm}'$ are obtained as exponential maps \emph{inside} $\Ads^3$, and this recovers the construction in \cite{kirill07}.
\end{remark}
\begin{prop} In the notations above, 
\[
\Is'_{\pm}:=(\sigma_{\pm}')^*\inners_{\PSL(2,\C)} =  \frac 12 \Is ((id\mp iB)\cdot , (id\mp iB)\cdot).
\]
In particular,  $\Is'_{\pm}$ is a complex metric if and only if $\mathrm K^{\Is}$ vanishes nowhere. When this happens, $\Is'_{\pm}$ both have constant curvature $-2$.
\end{prop}
\begin{proof}
Differentiating the expression for $\sigma_{\pm}$, we have that
\begin{align*}
\Is'_{\pm} =& \inner{\mathrm{d}\sigma_{\pm}, \mathrm{d}\sigma_{\pm}}= \frac 12 \inner{\mathrm{d}\sigma \pm i dN, \mathrm{d}\sigma \pm i dN}=\\
=&\frac 12 \inner{\mathrm{d}\sigma \mp \mathrm{d}\sigma\circ (iB), \mathrm{d}\sigma \mp \mathrm{d}\sigma\circ (iB) }= \\
=&\frac 12 \Is( id\mp iB, id\mp iB).
\end{align*}
For the second part of the statement, just observe that, by the Gauss equation~\eqref{eq: Gauss PSL(2,C)}, $\Is'_{\pm}$ is a complex metric if and only if $$0\ne \det(id\pm iB) =1-\det(B) = -\mathrm K^{\Is}.$$ Since $id\pm iB$ is Codazzi for $\Is$, $\det(id\pm iB)\ne 0$ implies that the Levi-Civita connection for $\Is'_{\pm}$ is given by $$\nabla^{\Is'_{\pm}}_X Y= (id\pm iB)^{-1}  (\nabla^{\Is}_X ( (id\pm iB)(Y) ) ). $$
Hence, by a standard computation (see for instance the proof of \cite[Proposition 3.9]{ElE}), one  gets that
    $$\mathrm K^{\Is_{\pm}'} =2\frac{K_{\Is}}{\det(id\pm iB)} = 2\frac{-1+\det(B)}{1-\det(B)}=-2. $$
\end{proof}
Explicitly, if $\gamma=(J_1, q_1, -J_2, \overline{q_2}, u)\in \mathrm M(S)$, then 
\begin{equation}
    \label{eq: metriche equidistanti}
\Is'_{\pm}= \frac 1 2 \left( 1- e^{-4u}\frac{q_1\overline {q_2}}{(\hpair)^2}\right) e^{2u}\hpair \mp \frac{i}{2} (q_1+\overline{q_2})\ .
\end{equation}

Notice that the map \begin{equation}
\label{eq: metriche equidistanti }
    \begin{split}
     \Lambda:=  \Is_{\pm}' \colon \mathrm M^*(S)  &\to \mathrm{Sym}_2(\C TS)\times \mathrm{Sym}_2(\C TS)\\
        \gamma=(J_1, -J_2, q_1, \overline{q_2}, u) &\mapsto (\Is'_+(\gamma), \Is'_-(\gamma))
    \end{split}
\end{equation}
is holomorphic. Indeed, using Proposition~\ref{prop: Bers metrics holo dependence}, we can see that it clearly extends to a holomorphic map from $\QD(S)\times\QD(S)\times C^{\infty}(S,\C)$, hence its restriction to the complex Fréchet submanifold $\mathrm M^*(S)$ is holomorphic. 
Moreover, both the maps $\Is_{\pm}'$ are equivariant for the action of $\mathrm{Diff}_0(S)$. Recall from Section \ref{sec: Bers metrics} that Bers metrics determine a connected component of the space of complex hyperbolic metrics with constant curvature $-1$ (\cite[Theorem 6.15]{BEE}). As a consequence, we can define the open subset of $\mathrm M^*(S)$ given by 
\[
    \mathrm W(S):=\left\{ \gamma\in \mathrm M^*(S)\ |\ \text{ $(2\cdot \Is'_{\pm}(\gamma))$ are both Bers metrics}
    \right\}.
\]

\begin{remark}
Observe that $\mathrm W(S)$ contains $\mathrm{MGH}(S)$: if $\gamma\in\mathrm{MGH}(S) $ then $\Is'_{\pm}$ are hyperbolic Riemannian metrics, which are Bers metrics. \end{remark}
\begin{lem}
The space $\mathrm{AF}(S)$ is contained in $\mathrm W(S)$.
\end{lem}
\begin{proof}
Let $\gamma=(J,q,-J, \overline q, u)\in \mathrm{AF}(S)$, let  $z$ be a local holomorphic coordinate for $J$ and write $q=\phee \mathrm{d}z^2$ and $h=\lambda \mathrm{d}z\mathrm{d}\overline z$ for the hyperbolic metric conformal to $J$. Then, the isotropic directions of $\Is'_{\pm}$ are given by the complex lines generated by \[
    X_{\pm}= \partial_z \pm \frac{\phee}{e^{2u}\lambda}\partial_{\overline z} \qquad \text{and} \qquad  Y_{\pm}=\partial_{\overline z} \pm \frac{\overline{\phee}}{e^{2u}\lambda}\partial_{ z}.
    \]
    A simple computation shows that a vector $\alpha \partial_z+\beta\partial_{\overline z}\ne 0$ determines a point in a connected component of $\mathbb P(\C TS)\setminus \mathbb P(TS)$ according to the sign of $Re(\frac{\alpha-\beta}{\alpha+\beta})$ or equivalently of $Re(\frac{\alpha+\beta}{\alpha-\beta})$ (or just one of the two when the other is not defined). In our setting, it is straightforward to check that $$\frac{\phee\overline{\phee}}{e^{4u}\lambda^2}<1 $$ implies that $X_{\pm}$ and $Y_{\pm}$ lie in opposite connected components, implying that $\Is'_{\pm}$ are both positive complex metrics. Finally, $\mathrm{AF}(S)$ is connected, $\Is'_{\pm}$ is clearly a Bers metric when $q=0$, and Bers metrics are connected components of the space of positive complex metrics of curvature $-1$: we conclude that $2\Is'_\pm$ are both Bers metrics.
\end{proof}
By looking at the holonomy of $\Is'_{\pm}$, we get a holomorphic map
\[
\widehat \Theta\colon \mathrm W(S)\to \QF(S)\times \QF( S),
\] where the target space is equipped with the product complex structure $(J_\chi,J_\chi)$ (see Section~\ref{sec: tensors I,J,K}).
\begin{lem}\label{lem: Theta holomorphic}
    There exists an open connected subset $\mathbb W(S)\subset \mathbb M_0(S)$ containing both $\mathcal{AF}(S)$ and $\mathcal{MGH}(S)$ such that $\widehat \Theta$ descends to a holomorphic map
    \[
    \Theta\colon \big(\mathbb W(S),\mathbf i\big) \to  \big(\QF(S)\times \QF(S), (J_\chi,J_\chi)\big),
    \]
 namely $\Theta([\gamma])=\widehat\Theta(\gamma)$ for an open subset of representatives $\gamma$ of $[\gamma]$. 
\end{lem}
\begin{proof}
Denote $Z= \mathrm{AF}(S)\cup \mathrm{MGH}(S)$ and $\mathcal Z=\mathcal{AF}(S)\cup \mathcal{MGH}(S)\subset \mathbb M(S)$. 

Since $\widehat \Theta$ is $\mathrm{Diff}_0(S)$--invariant, it is constant within each fiber of the projection ${Z}\to \mathcal Z$, defining continuously the map $\Theta$ on $Z$.
Moreover, since $\widehat \Theta$ is also holomorphic, we get that, for all $\gamma\in \mathrm W(S)$ and $X\in \Gamma(\C TS)$, 
$$
\mathrm{d}\widehat \Theta(\mathscr L_X\gamma)=\mathrm{d}\widehat \Theta(\mathscr L_{Re (X)}\gamma)+i \mathrm{d}\widehat \Theta(\mathscr L_{Im (X)}\gamma)= 0~.
$$ 
Hence, by Proposition \ref{prop: tangent fiber M*}, $\mathrm{d}\widehat \Theta$ is constantly zero on the kernel of the differential map of the projection $ \mathrm W(S)\to \mathbb M(S)$.

Since the map $\pi\colon \mathrm M^*(S)\to \mathrm{QD}(S)\times\mathrm{QD}(S)$ is a local biholomorphism, for all $\gamma$ in $\mathrm Z$ we can find a neighborhood $U\subset \mathrm W(S)$ such that the fibers of the projection $p\colon \mathrm U\to \mathbb M(S)$ are connected. The map $\widehat \Theta|_{U}$ is therefore constant among the fibers of $U\to \mathbb M(S)$, hence $\widehat \Theta|_U$ descends to a well-defined map $\Theta$ in a neighborhood of $[\gamma]$.
It is simple to check that this local definition of $\Theta$ is holomorphic. Indeed, using the local biholomorphisms $\mathrm M^*(S)\to \mathrm{QD}(S)\times  \mathrm{QD}(\overline S)$ and $\mathbb M(S)\to \mathbb{QD}(S)\times\mathbb{QD}(\overline S)$, we can construct local charts for $\mathrm M^*(S)$ of the form by $U_1\times U_2$ with $U_1\subset \mathbb{M}(S)$ and $U_2\subset \mathrm{Diff}_0(S)\times \mathrm{Diff}_0(S)$, over which $\widehat \Theta(u_1,u_2)= \Theta(u_1)$: the holomorphicity of $\Theta$ then follows by the holomorphicity of $\widehat \Theta$.

We are left to prove that these local descriptions of $\Theta$ along $Z$ match. This follows from a variation of standard topological argument (see the proof of Theorem 6.20 in \cite{ESrank2}, or the suggested proof for the \emph{Generalized inverse function theorem} in Exercise 1.8.14 in \cite{GuillPoll} and \cite{GIFT}).

Let $m$ be a Riemannian metric on a neighborhood of $\mathcal Z$ inside $\mathbb M(S)$ such that $\mathcal Z_1:=\mathcal{AF}(S)$, $\mathcal Z_2:=\mathcal{MGH}(S)$ and $\mathcal Z_3:=\mathcal{AF}(S)\cap \mathcal{MGH}(S)$ are all totally geodesic. This can be done for instance by restricting $\pi$ a global biholomorphism from a neighborhood of $\mathcal Z$ inside $\mathbb M(S)$ to an open subset of $\mathbb{QD}(S)\times\mathbb{QD}(\overline S)$ as in Theorem 6.20 in \cite{ESrank2}, and by pulling back a product metric coming from a Riemannian metric on $\mathbb{QD}(S)$ for which $[c,q]\mapsto [c,-q]$ is an isometry. 

For each $\mathcal Z_i$, one can construct a locally finite cover $\{V_i^j\}$ of open subsets of $\mathbb M(S)$ such that:
\begin{itemize}
    \item the construction above defines, on each $V_i^j$, a local holomorphic extension $\Theta_i^j$ coinciding with $\Theta|_{\mathcal Z_i}$ on $V_i^j\cap \mathcal Z_i$,
    \item each $\overline {V_i^j}$ is a geodesically convex ball,
    \item every point of $V_i^j$ has a unique closest point projection to $\mathcal Z_i$, with the segment realizing the distance being contained in $V_i^j$.
\end{itemize} 
With these assumptions, it is clear that $V_i^j$ and $V_i^k$ intersect if and only if $V_i^j\cap V_i^k\cap \mathcal Z_i$ is connected and nonempty, therefore $\Theta_i^j$ and $\Theta_i^k$ coincide on $V_i^j\cap V_i^k$ by uniqueness of the holomorphic continuation of $\Theta|_{\mathcal Z_i}$. As a result, we obtain a globally defined $\Theta_i$ on the open subset $\mathbb W_i=\cup_j V_i^j$.

Finally, up to shrinking the $\mathbb W_i$'s, we can assume that $\mathbb W_1\cap \mathbb W_2= \mathbb W_3$. Clearly, $\Theta_1$ and $\Theta_3$ must coincide on $\mathbb W_3$ by uniqueness of the complex analytic extension on $\mathbb W_3\cap \mathcal Z_1$. By applying the same argument for $\mathbb W_2\cap \mathbb W_3$, we get that on $\mathbb W(S):=\mathbb W_1\cup \mathbb W_2$ we can define $\Theta=\Theta_1\cup \Theta_2$.
\end{proof}

\begin{thm}\label{prop: iK on W_0}
On the open subset $\mathbb W(S)\subset\mathbb M_0(S)$, the tensor $\mathbb K$ coincides with $\Theta^*(-J_\chi,J_\chi)$.
\end{thm}
\begin{proof}
By construction and by Proposition~\ref{prop: restriction of K} the tensor \(\mathbf i\circ\mathbb K\) is holomorphic on \(\mathbb{W}_0(S)\) and its restriction to \(\mathcal{MGH}(S)\) is equal to \(\mathbf{K}'\). Since the tensor \(\mathbf i\circ\Theta^*(-J_\chi,J_\chi)\) is also holomorphic on \(\mathbb{W}_0(S)\), if we show that its restriction to the anti-de Sitter locus coincides with \(\mathbf{K}'\), then we conclude by the analytic continuation argument (Lemma~\ref{lem: fundamental lem analytic continuation}). This follows from the following chain of equalities:
\begin{align*}
   \mathbf i\circ\Theta(-J_\chi,J_\chi)|_{\mathcal{MGH}(S)}&= \Theta^*(J_\chi,J_\chi)\circ\Theta^*(-J_\chi,J_\chi)|_{\mathcal{MGH}(S)} \tag{Lemma \ref{lem: Theta holomorphic}} \\ & =\Theta^*(\mathrm{Id},-\mathrm{Id})|_{\mathcal{MGH}(S)} \\ &=\mathcal{C}^*(\mathrm{Id},-\mathrm{Id}) \tag{Remark \ref{rmk:constant curvature surfaces}} \\ &=\mathbf{K}' \ . \tag{\cite[Theorem D]{MazSepTamb}}
\end{align*}
\end{proof}

\section{Bi-complex geometry}\label{sec:bi-complex}

\subsection{Bi-complex numbers}\label{sec:bicomplex_numbers}
\noindent The algebra of bi-complex numbers $\mathbb{C}_\tau$ (\cite{bi_complexbook}) is the commutative algebra over $\mathbb{R}$ generated by $\{1, \tau, i, j:=\tau i=i\tau\}$ where $\tau^2=1$ and $i^2=j^2=-1$. In particular, any element $w\in\C_\tau$ can be written as $w=x_1+\tau y_1+ix_2+jy_2$, for $x_i,y_i\in\R$. The algebra $\C_\tau$ contains two copies of the complex numbers, one with respect to $i$ and one with respect to $j$. The first is obtained by imposing $y_{1}=y_{2}=0$, whereas the second arises by imposing $y_{1}=x_{2}=0$. There is another distinguished two-dimensional subalgebra contained in $\C_\tau$, generated by $\{1,\tau\}$. In the literature this is known as the algebra of para-complex (or split-complex) numbers, and we shall denote it by $\R_\tau$.
 \\ \\
\noindent In what follows, we shall write a bi-complex number as $w = z_1 + \tau z_2$, where $z_k = x_k + iy_k$ are complex numbers, for $k = 1,2$. In other words, $\C_\tau=\R_\tau\otimes_\R\mathbb C$, where $\C$ is the copy of complex numbers with respect to $i$. Given a bi-complex number $w=z_{1}+\tau z_{2}$, we will use the standard notation $\overline{w}$ to denote complex conjugation, so that $\overline{w}=\overline{z}_{1}+\tau\overline{z}_{2}$. We then define
\[
    \Re(w):=\frac{w+\overline{w}}{2}  \ \ \ \ \  \text{and} \ \ \ \ \ \Im(w):=\frac{w-\overline{w}}{2i} \ . 
\]
These are para-complex numbers with the property that $w=\Re(w)+i\Im(w)$. In a similar way, we denote by $w^{\dagger}$ the conjugate of $w$ with respect to the para-complex unit $\tau$, in other words $w^{\dagger}:=z_{1}-\tau z_{2}$. We then define
\[
    \Re_{\tau}(w):=\frac{w+w^{\dagger}}{2}  \ \ \ \ \  \text{and} \ \ \ \ \ \Im_{\tau}(w):=\frac{w-w^{\dagger}}{2\tau} \ 
\]
as the $\tau$--real and $\tau$--imaginary part of $w$. Note that, despite the name, $\Re_{\tau}(w)$ and $\Im_{\tau}(w)$ are complex numbers. \\

\noindent The $\tau$--norm of a bi-complex number $w=z_{1}+\tau z_{2}$ is the complex number 
\[
    |w|_{\tau}:=w{w}^{\dagger}=z_{1}^2-z_{2}^{2} \ .
\]
The non-vanishing of the $\tau$--norm characterizes invertible bi-complex numbers. 

\noindent For what follows, it will be useful to write a bi-complex number $w=z_{1}+\tau z_{2}$ with respect to the basis of idempotents $\{e_{+}, e_{-}\}$, where
\[
    e_{+}:= \frac{1+\tau}{2} \ \ \ \ \ \text{and} \ \ \ \ \ e_{-}:= \frac{1-\tau}{2} \ .
\]
By a straightforward computation, we can rewrite $w=z^{+}e_{+}+z^{-}e_{-}$ with 
\(
    z^{+}=z_{1}+z_{2} \ \text{and} \  z^{-}=z_{1}-z_{2}    
\), so that we get an isomorphism of algebras between $\C_\tau$ and $\C\oplus\C$.

\subsection{Para-complex coordinates}  Let $(M, \mathbf{i})$ be a complex manifold. Recall that a holomorphic almost para-complex structure $\mathbf{P}$ on $M$ is a holomorphic ($\C$--linear) bundle endomorphism $\mathbf{P}:T^{1,0}M\rightarrow T^{1,0}M$ such that $\mathbf{P}^{2}=\id$ and the eigenspaces of $\mathbf{P}$ relative to the eigenvalues $\pm1$ have the same dimension (see Definition~\ref{def:para-complex structure}). \\

Let $(M, \mathbf{i}, \mathbf{P})$ and $(N, \mathbf{i'},\mathbf{P}')$ be two holomorphic almost para-complex manifolds. A holomorphic map $f:M\rightarrow N$ is bi-holomorphic if $\mathbf{P}'\circ \mathrm{d}f=\mathrm{d}f\circ \mathbf{P}$. \\

An almost holomorphic para-complex structure $\mathbf{P}$ on $M$ is said to be integrable if the eigen-distributions of $\mathbf{P}$ are involutive. This is equivalent to the existence of local charts $\phi_{\alpha}:U_{\alpha}\subset M \rightarrow \C_{\tau}^{n}$ such that change of coordinates are bi-holomorphic functions (cfr. \cite[Section 2.6]{loustau-sanders}). If $z^{i}$ is a $\C_{\tau}$--valued coordinate, its $\tau$--real and $\tau$--imaginary parts define holomorphic $\C$--valued coordinates on $M$. 

\begin{defn} A \emph{bi-holomorphic manifold} is a triple $(M, \mathbf{i},\mathbf{P})$, where $(M, \mathbf{i})$ is a complex manifold and $\mathbf{P}$ is a holomorphic integrable para-complex structure.    
\end{defn}

The bi-complexified tangent bundle of a bi-complex manifold $(M,\mathbf{i}, \mathbf{P})$ is $\C_{\tau}TM=\C TM\otimes_{\C}\C_{\tau}$. 
Extending the endomorphisms $\mathbf{i}$ and $\mathbf{P}$ on $\C_{\tau}TM$ by $\C_{\tau}$--linearity, since $\mathbf{i}$ and $\mathbf{P}$ commute, we have an eigenspace decomposition 
\[
    \C_{\tau}TM=T^{1,0}_{+}M\oplus T^{1,0}_{-}M \oplus T^{0,1}_{+}M \oplus T^{0,1}_{-}M \ ,
\]
where $T^{1,0}_{\pm}M$ denotes tangent vectors that are eigenvectors for $\mathbf{i}$ relative to the eigenvalue $i$ and eigenvectors for $\mathbf{P}$ relative to the eigenvalue $\pm \tau$. The definition of $T^{0,1}_{\pm}M$ is analogous but now tangent vectors are eigenvectors for $\mathbf{i}$ relative to $-i$.

In local $\C_{\tau}$--valued bi-holomorphic coordinates $\xi^{i}=z^{i}+\tau w^{i}$, the vector spaces $T^{1,0}_{+}M$ and $T^{1,0}_{-}M$ are spanned by the vector fields
\[
    \frac{\partial}{\partial \xi^{i}}=\frac{1}{2}\left( \frac{\partial}{\partial z^{i}}+\tau \frac{\partial}{\partial w^{i}}\right) \ \ \ \ \ \frac{\partial}{\partial {\xi^{i}}^{\dagger}}=\frac{1}{2}\left( \frac{\partial}{\partial z^{i}}-\tau \frac{\partial}{\partial w^{i}}\right) \ .
\]
Similarly, the induced para-complex structure $\mathbf{P}^{*}$ on the bi-complex cotangent bundle $\C_{\tau}T^{*}M=\C T^{*}M \otimes_{\C} \C_{\tau}$ induces a decomposition into eigenspaces 
\[
    \C_{\tau}T^{*}M=\Lambda^{1,0}_{+}M\oplus \Lambda^{1,0}_{-}M\oplus \Lambda^{0,1}_{+}M\oplus \Lambda^{0,1}_{-}M .
\]
In local $\C_{\tau}$--valued bi-holomorphic coordinates $\xi^{i}=z^{i}+\tau w^{i}$, the vector spaces $\Lambda^{1,0}_{+}(M)$ and $\Lambda^{1,0}_{-}(M)$ are spanned by the $\C_{\tau}$--valued holomorphic 1-forms
\[
    \mathrm{d}\xi^{i}=\mathrm{d}z^{i}+\tau \mathrm{d}w^{i}  \ \ \ \text{and} \ \ \ \mathrm{d}{{\xi}^{i}}^{\dagger}=\mathrm{d}z^{i}-\tau \mathrm{d}w^{i} \ .
\]
More in general, we have a decomposition of $\C_{\tau}$--valued $n$--forms on $M$ into types. In particular a $\C_{\tau}$--valued $2$--form $\Omega$ can be written uniquely as sum 
\[
    \Omega^{2,0}_{++}+\Omega^{2,0}_{+-}+\Omega^{2,0}_{--}+\Omega^{1,1}_{++}+\Omega^{1,1}_{+-}+\Omega^{1,1}_{--}+\Omega^{0,2}_{++}+\Omega^{0,2}_{+-}+\Omega^{0,2}_{--},
\]
where upper indices express the behavior with respect to $\mathbf{i}$ and lower indices are relative to $\mathbf{P}$, for example
\begin{align*}
    \Omega^{1,1}_{+-}(\mathbf{P}X,Y)&=-\Omega^{1,1}_{+-}(X,\mathbf{P}Y) \\
    \Omega^{1,1}_{+-}(\mathbf{i}X,Y)&=-\Omega^{1,1}_{+-}(X,\mathbf{i}Y)
\end{align*}
for all bi-complex vector fields $X$ and $Y$. 
\begin{defn}\label{def:bicomplex_sympl} Let $(M,\mathbf{i}, \mathbf{P})$ be a bi-holomorphic manifold. A \emph{bi-holomorphic symplectic form} on $M$ is a $\C_{\tau}$--valued $2$--form $ \Omega=\Omega_{1}+\tau \Omega_{2}$ that is non-degenerate on $T^{1,0}_{+}M$ such that
\begin{enumerate}
    \item $\Omega(X,\mathbf{i}Y)=\Omega(\mathbf{i}X,Y)=i\Omega(X,Y)$ for all bi-complex vector fields $X$ and $Y$;
    \item $\Omega(X,\mathbf{P}Y)=\Omega(\mathbf{P}X,Y)=\tau\Omega(X,Y)$ for all bi-complex vector fields $X$ and $Y$;
    \item $\mathrm{d} \Omega=0$.
\end{enumerate}
\end{defn}

The first two conditions imply that $\Omega$ is of type $(2,0)$ with respect to $\mathbf{i}$ and of type $(++)$ with respect to $\mathbf{P}$, hence in local bi-holomorphic coordinates $\xi^{i}=z^{i}+\tau w^{i}$, it can be written as
\[
    \Omega=\sum_{i<j} f_{i,j} \mathrm{d}\xi^{i}\wedge \mathrm{d}\xi^{j} , 
\]
for some $\C_{\tau}$--valued smooth functions $f_{i,j}$. The third condition implies that $f_{i,j}$ are bi-holomorphic functions (\cite[Theorem 7.6.5]{bi_complexbook}), i.e. 
\[
    \frac{\partial f_{i,j}}{\partial \bar{\xi}^{k}}=\frac{\partial f_{i,j}}{\partial \bar{\xi}^{k\dagger}}=\frac{\partial f_{i,j}}{\partial \xi^{k \dagger}}=0 \ \ \ \  \text{for all \ $i,j,k$.}
\]

\subsection{Re-interpretation of our results}\label{subsec: bi-complex Goldman}
The group $G_{\tau}=\PSL(2,\C_\tau)$ is isomorphic to $\PSL(2,\C)\times \PSL(2,\C)$. The smooth locus of the $G_{\tau}$--character variety has a natural bi-complex structure, given by multiplication by $i$ and $\tau$ on $\mathfrak{g}_{\tau}$--valued co-cycles.
Let $\omega_{G_{\tau}}$ denote Goldman's symplectic form on the $G_{\tau}$--character variety, which is thus $\C_{\tau}$--valued. Let $\omega_{G}$ denote Goldman's symplectic form on the $\PSL(2,\C)$--character variety. By the identification $G_{\tau}\cong \PSL(2,\C)^2$, it makes sense to consider the open dense subset $\chi^*(\pi_1(S), G_{\tau})$ corresponding to pairs of irreducible representations into $\PSL(2,\C)$, which comes with natural projections $\pi_{r,l}$ to $\chi^{irr}(\pi_1(S), \PSL(2,\C))$.   It is easy to verify that 
\begin{equation}\label{eq:omega_G_tau}
    \omega_{G_{\tau}}=\frac{1}{2}(\pi_{l}^{*}\omega_{G}+\pi_{r}^{*}\omega_{G})+\frac{\tau}{2}(\pi_{l}^{*}\omega_{G}-\pi_{r}^{*}\omega_{G})
\end{equation}
The same computations as in \cite[Section 5.4.2]{MazSepTamb}, with the only difference that the group $\PSL(2,\R_{\tau})$ is replaced by $\PSL(2,\C_{\tau})$, show that $\omega_{G_{\tau}}$ satisfies Condition $(1)$ and $(2)$ in Definition~\ref{def:bicomplex_sympl}. Since it is closed, we obtain the following:

\begin{thm} Goldman's symplectic form on $\chi^*(S,G_{\tau})$ is bi-holomorphic.
\end{thm}

We can also describe global bi-holomorphic Darboux coordinates for $\omega_{G_{\tau}}$. By work of Wolpert (\cite{Wolpert1982FenchelNielsen}, \cite{Wolpert1983SymplecticGeometry}), after fixing a pants decomposition $\{\alpha_{j}\}_{j}$ on $S$, Fenchel-Nielsen length-twist parameters $(\ell_{\rho}(\alpha_{i}), \theta_{\rho}(\alpha_{i}))$ are global Darboux coordinates for the Weil-Petersson, and hence for the Goldman, symplectic form on Teichm\"uller space. These coordinates were later extended holomorphically to the space of quasi-Fuchsian representations $\mathcal{QF}(S)$ by Kourouniotis (\cite{Kourouniotis1991Bending}, \cite{Kourouniotis1992Geometry}, \cite{Kourouniotis1994Complex}) and Tan (\cite{Tan1994ComplexFenchelNielsen}). Platis (\cite{Platis2001ComplexSymplectic}) then developed a complex version of Wolpert's theory and Loustau (\cite{Loustau2014ComplexSymplectic}) observed that these complex Fenchel-Nielsen coordinates $(\ell_{\rho}^{\C}(\alpha_j), \theta_{\rho}^{\C}(\alpha_j))$ are Darboux for the restriction to $\mathcal{QF}(S)$ of the Goldman symplectic form on $\chi^{irr}(\pi_{1}(S), \PSL(2,\C))$. On the other hand, Mazzoli-Seppi-Tamburelli (\cite{MazSepTamb}) defined para-complex length and twist functions $(\ell_{\rho}^{\R_\tau}(\alpha_j), \theta_{\rho}^{\R_\tau}(\alpha_j))$ to parameterize $\mathcal{MGH}(S)$ and found that they provide a global Darboux frame for the Goldman symplectic form on $\chi(\pi_{1}(S), \PSL(2, \mathbb{R}_{\tau}))$. Moreover, these para-complex Fenchel-Nielsen coordinates are closely related to the classical ones via the isomorphism $\PSL(2,\R_{\tau}) \cong \PSL(2,\R)\times \PSL(2,\R)$: indeed, if $\rho \in \mathcal{MGH}(S)$ corresponds to the pair $(\rho_{+}, \rho_{-}) \in \mathcal{T}(S) \times \mathcal{T}(S)$ via Mess' parameterization, then
\[
 \ell^{\R_\tau}_{\rho}(\alpha_{j})=\left(\frac{ \ell_{\rho_{+}}(\alpha_{j})+\ell_{\rho_{-}}(\alpha_{j})}{2}\right)+\tau\left(\frac{ \ell_{\rho_{+}}(\alpha_{j})-\ell_{\rho_{-}}(\alpha_{j})}{2}\right) \ 
\]
and
\[
    \theta^{\R_{\tau}}_{\rho}(\alpha_{j})=\left(\frac{ \theta_{\rho_{+}}(\alpha_{j})+\theta_{\rho_{-}}(\alpha_{j})}{2}\right)+\tau\left(\frac{ \theta_{\rho_{+}}(\alpha_{j})-\theta_{\rho_{-}}(\alpha_{j})}{2}\right) \ .
\]
It follows immediately that, if now $\rho=(\rho_{+}, \rho_{-})$ is a representation in $\PSL(2,\C)\times \PSL(2,\C)$ with $\rho_{\pm}$ quasi-Fuchsian, the coordinates
\[
 \ell^{\C_\tau}_{\rho}(\alpha_{j})=\left(\frac{ \ell_{\rho_{+}}^{\C}(\alpha_{j})+\ell^{\C}_{\rho_{-}}(\alpha_{j})}{2}\right)+\tau\left(\frac{ \ell^{\C}_{\rho_{+}}(\alpha_{j})-\ell^{\C}_{\rho_{-}}(\alpha_{i})}{2}\right) \ 
\]
and
\[
    \theta^{\C_{\tau}}_{\rho}(\alpha_{j})=\left(\frac{ \theta^{\C}_{\rho_{+}}(\alpha_{j})+\theta^{\C}_{\rho_{-}}(\alpha_{j})}{2}\right)+\tau\left(\frac{ \theta^{\C}_{\rho_{+}}(\alpha_{j})-\theta^{\C}_{\rho_{-}}(\alpha_{j})}{2}\right) \ .
\]
are well-defined, bi-holomorphic and, by Lemma~\ref{lem: fundamental lem analytic continuation}, are global Darboux coordinates for $\omega_{G_{\tau}}$, in other words
\[
    \omega_{G_{\tau}}=\frac{1}{2}\sum_{j=1}^{n} \mathrm{d}\ell^{\C_{\tau}}_{\rho}(\alpha_{j})\wedge \mathrm{d}\theta^{\C_{\tau}}_{\rho}(\alpha_j) \ .
\]

In this language of bi-holomorphic symplectic geometry, we can rephrase the results of Section~\ref{sec:complex_hyperkahler} as follows:

\begin{thm}\label{thm: bi-complex reinterpretations} On $\mathbb M_0(S)$ we have
\begin{equation}\label{eq:hol_omega_g_tau}
    2\mathrm{hol}^{*}\omega_{G_{\tau}}=\omega_{\mathbb{I}}+i\tau\omega_{\mathbb{K}},  
\end{equation}
hence $\mathrm{hol}^{*}\omega_{G_{\tau}}$ is bi-holomorphic with respect to $(\mathbf{i}, \mathbf{i}\circ\mathbb{J})$. Moreover, on the open subset $\mathbb W_0(S)\subset\mathbb M_0(S)$ of Theorem~\ref{prop: iK on W_0}, the following identity holds
\begin{equation}\label{eq:Theta_omega_G_tau}
    2\Theta^{*}\omega_{G_{\tau}}=\omega_{\mathbb{I}}+i\tau\omega_{\mathbb{J}}. 
\end{equation}
In particular, $\Theta^{*}\omega_{G_{\tau}}$ is bi-holomorphic with respect to $(\mathbf{i}, -\mathbf{i}\circ\mathbb{K})$.
\end{thm}
\begin{proof} Equation (\ref{eq:hol_omega_g_tau}) follows immediately combining Equation~\eqref{eq:omega_G_tau} with our definition of $\omega_{\mathbb{I}}$ in~\eqref{eq: omega_mathbb I} and Proposition~\ref{prop: formula for omega_K}. The right-hand side is bi-holomorphic with respect to $(\mathbf{i}, \mathbf{i}\circ\mathbb{J})$. \\
Equation~\eqref{eq:Theta_omega_G_tau} holds on the maximal totally real submanifold $\mathcal{MGH}(S)$ by \cite[Theorem D]{MazSepTamb}, hence it is true everywhere on the open subset $\mathbb W_0(S)\subset\mathbb M_0(S)$ which contains $\mathcal{MGH}(S)$. Again, the right-hand side is manifestly bi-holomorphic with respect to $(\mathbf{i}, -\mathbf{i}\circ\mathbb{K})$.
\end{proof}

\appendix

\section{Analytic continuation}\label{appendix:extension}

\begin{defn}\label{def: totally real submanifold} Let $(\mathbb X,\mathbf i)$ be a complex manifold. A real submanifold $Z\subset \mathbb X$ is said to be \emph{totally real} if \[ T_xZ\cap \mathbf i(T_xZ)=\{0\} \] for every $x\in Z$. It is said to be \emph{totally real of full dimension} (or \emph{maximal totally real}) if, in addition, \( \dim_{\mathbb R}Z=\dim_{\mathbb C}\mathbb X. \) Equivalently, in this case, \[ T_x\mathbb X=T_xZ\oplus \mathbf i(T_xZ) \] as real vector spaces, for every $x\in Z$. \end{defn}


\begin{lem}
\label{lem: fundamental lem analytic continuation}
    Let $\mathbb M$ and $\mathbb N$ be complex manifolds with $\mathbb M$ connected, and let $Z\subset \mathbb M$ be a totally real submanifold of full dimension ($\text{dim}_\R Z=\text{dim}_{\C}\mathbb{M}$). If $f_1, f_2\colon \mathbb M\to \mathbb N$ are holomorphic maps that coincide along $Z$, then they coincide everywhere.
\end{lem}
\begin{proof}
    By complex analiticity, it is sufficient to prove that if $f_1$ and $f_2$ coincide on an open subset of $\mathbb M$, then they coincide everywhere because $\mathbb M$ is connected. Fix a point $x\in Z$ and, up to composing with local charts around $x$, see $f_1, f_2$ as holomorphic maps between open subsets $\mathbb C^m\supset U\to V\subset \mathbb  C^n$, and see $Z\cap U$ as a totally real submanifold of $U$. 
    
   It is now simple to check that $f_1$ and $f_2$ have the same power series in the point $x\in U$. Let $X_1, \dots X_m $ be a local frame of linearly independent vector fields on $Z$ in a neighbrohood of $x$ that, up to shrinking the chart, we can assume to be the whole $U\cap Z$. On $U\cap Z$, we clearly have that $\partial_{X_k} f_1 =\partial_{X_k}f_2$, hence, by holomorphicity of $f_1$ and $f_2$, we also have that $\partial_{\mathbf J X_k} f_1 =i\partial_{X_k} f_1 =i \partial_{X_k}f_2 =\partial_{\mathbf JX_k}f_2$: as a consequence, in every point of $U\cap Z$, the maps $f_1$ and $f_2$ have the same 1-jet as functions on $U$. In particular the first order partial derivatives $\partial_{z_k} f_j$ induced by the holomorphic coordinates $(z_1, \dots, z_n)$ of $\C^n$ coincide on $U'\cap Z$ for $j=1$ and $j=2$. Applying the same argument by induction on $k$ to the holomorphic functions $\partial^k_{z_{i_1}, \dots, z_{i_k}} f_1$ and $\partial^k_{z_{i_1}, \dots, z_{i_k}} f_2$, we conclude that $f_1$ and $f_2$ share the same partial derivatives of any order in $x$. As a consequence, their power series at $x$ have the same coefficients, hence $f_1$ and $f_2$ coincide in an open neighborhood of $x$ inside $\C^n$.
\end{proof}
\begin{remark}\label{rmk: fundamental_lemma_applied to tensors} The same argument applies to holomorphic sections of a holomorphic vector bundle. Indeed, after choosing a local holomorphic trivialization, the claim reduces component-wise to Lemma~\ref{lem: fundamental lem analytic continuation}. Consequently, two holomorphic tensors on $\mathbb M$ that coincide along $Z$ must coincide everywhere on $\mathbb M$. \end{remark}

\section{Infinitesimal variation of mean curvature}\label{appendix:inf mean curv}

The purpose of this appendix is to prove the formula for the variation of the mean curvature of a minimal immersion in $\PSL(2,\C)$ by a normal variation, generalizing the well-known formulae in the context of $\Hyp^3$ and $\Ads^3$.

\begin{lem}\label{lemma: formula variation H}
    Let $\Sigma$ be an oriented surface, let $\sigma:\Sigma\to \PSL(2,\C)$ be an admissible minimal immersion in $\PSL(2,\C)$, let $f\in C^\infty(\Sigma,\C)$  and set $$\sigma_t(p):=\exp_{\sigma(p)}(tf(p)N(p))~.$$
    Then, if $H_\sigma$ denotes the mean curvature of $\sigma$, then
  
    $$\left. \frac{\mathrm{d}}{\mathrm{d}t} \right|_{t=0} H_\sigma =\Delta_{I}f +(\tr(B^ 2)-2)f$$
    where $\Is$ denotes the complex first fundamental form of $\sigma$ and $B$ the complex shape operator.
\end{lem}

In the proof, we will locally lift $\sigma_t$ to the double covering $\SL(2,\C)$, which will allow us to use the following matrix model. Realizing $\SL(2,\C)$ as an embedded submanifold in $\mathrm{Mat}(2,\C)$, the space of complex 2-by-2 matrices, the holomorphic Riemannian metric of $\SL(2,\C)$ (that is, the pull-back of the holomorphic Riemannian metric of $\PSL(2,\C)$) is equal to the restriction of the symmetric bilinear form on $\mathrm{Mat}(2,\C)$ whose associated quadratic form is $-\det$ (see \cite[\S~2.5]{BEE}). Explicitly, this bilinear form has the expression

\begin{equation}
\inner{ M,  N}= \frac 1 2 \left(\tr( M N)-\tr( M)\tr( N)\right)~,
\end{equation}
for $M,N$ in $\mathrm{Mat}(2,\C)$. 

Since $\SL(2,\C)=\{M\in \mathrm{Mat}(2,\C)|\inner{M,M}=-\det(M)=-1\}$, this has the practical advantage that we can perform the computation analogously to how one does in the models of the sphere, the hyperbolic space, the Anti-de Sitter space and so on, as quadrics inside pseudo-Euclidean spaces.

For instance, we have the following expression for the holomorphic Levi-Civita connection of $\SL(2,\C)$.

\begin{lem}
    Let $d$ denote the flat connection of the vector space $\mathrm{Mat}(2,\C)$, and $D$ denote the holomorphic Levi-Civita connection of $\SL(2,\C)$. Given $A$ in $\SL(2,\C)$ and two vector fields $X,Y$ on $\SL(2,\C)$, 
    \begin{equation} \label{eq: ambient connections SL dentro Mat} 
d_X Y|_A = D_X Y|_A + \inner{X,Y}\cdot A \in T_A\mathrm{Mat}(2,\C)~.
\end{equation}
\end{lem}

\begin{proof}
Decomposing  $d_X Y |_A \in T_A \mathrm{Mat}(2,\C)$ with respect to the direct sum decomposition $T_A \mathrm{Mat}(2,\C) \cong T_A \SL(2,\C) \oplus \mathrm{Span}_{\C}(A)$, the component along $\mathrm{Span}_{\C}(A)$ is characterized by the following identity:
$$
\inner{d_X Y|_A, A} = X(\inner{Y,\mathrm{Id}})|_A - \inner{Y, d_X \mathrm{Id}}|_A = -\inner{Y, d_X\mathrm{Id}}|_A= -\inner{X,Y} .
$$
Moreover, the tangent component to $\SL(2,\C)$ is equal to $D_X Y|_A$, where $D$ is the holomorphic Levi-Civita connection of $\mathrm{SL(2,\C)}$.
\end{proof}




\begin{proof}
We locally lift $\sigma_t$ to $\mathrm{SL}(2,\mathbb C)$, and we denote  by $\Is_t$, $\II_t$, and $N_t$ the corresponding fundamental forms and normal vector fields. Recalling that
$$\dot\sigma=\left.\frac d {dt}\right|_{t=0}\sigma_t=fN~,$$
and observing that, as a consequence of Equation~\eqref{eq: ambient connections SL dentro Mat}, for every $X \in\Gamma(\C T S)$,
$$
d_X N:=d_{\mathrm{d}\sigma(X)} N= D_{\mathrm{d}\sigma(X)} N =-\mathrm{d}\sigma(B(X)) ,
$$
we compute:

\begin{align*}
    \dot \Is(X,Y) &=\left.\frac d {dt}\right|_{t=0} \inner{d_X\sigma_t, d_Y\sigma_t} \\ &= \inner{d_X\dot \sigma, d_Y\sigma}+ \inner{d_X\sigma, d_Y\dot\sigma} \\
    &=\inner{\mathrm{d}f(X)N + fd_X N, d_Y\sigma}+ \inner{d_X\sigma, \mathrm{d}f(Y)N +fd_Y N}\\
    &=f\inner{d_X N, d_Y\sigma}+ f\inner{d_X\sigma,d_Y N}\\
    &=f\Is(-B(X),Y)+f\Is(X,-B(Y)) \\
    &=-2f \II(X,Y)~.
\end{align*}    


To compute the infinitesimal variation of $\II_t$, we first show the following formula for the infinitesimal variation of $N_t$: 
$$\dot N= f\sigma - \mathrm{grad}f~,$$
where $\mathrm{grad} f$ denotes the complex gradient vector field associated with the function $f$. 
To see this, we differentiate the identities $\langle N_t,N_t\rangle=1$, $\langle N_t,\sigma_t\rangle=0$ and $\langle N_t,\mathrm{d}\sigma_t(X)\rangle=0$, to obtain $\inner{\dot N, N}=0$, $\inner{\dot N, \sigma}=-\inner{N, \dot \sigma}= -f$ and 
$$\inner{\dot N, \mathrm{d}\sigma(X)}= -\inner{N, \mathrm{d}\dot \sigma (X)}= - \inner{N, \mathrm d f(X) N + f D_X N }= -\mathrm df(X)~.$$ 
This implies the formula.

Therefore, the infinitesimal variation of the second fundamental form is given by the following computation:
        \begin{align*}
        \dot{\II}(X,Y) & = - \left. \frac{d}{dt} \right|_{t=0} \inner{d_X\sigma_t, D_Y N_t} \\
        & = - \left. \frac{d}{dt} \right|_{t=0} \inner{d_X \sigma_t, d_Y N_t} \\
        & = - \inner{d_X \dot \sigma, d_Y N} - \inner{d_X \sigma, d_Y \dot N} \\
        &=- \inner{\mathrm{d}f(X) \cdot N + f d_X N, d_Y N} - \inner{d_X \sigma, -d_Y \mathrm{grad} f + \mathrm{d}f(Y) \cdot \sigma + f d_Y\sigma} \\
        &=- \inner{fd_X N, d_Y N} + \inner{d_X \sigma, \mathrm{d}\sigma(\nabla_Y^{\Is} \mathrm{grad} f)} - \inner{d_X \sigma,f d_Y \sigma} \\
        &=-f \Is(B(X), B(Y))+(\Hess_{\Is} f)(X,Y) - f \Is(X,Y) ,
    \end{align*}
    where, in the second-to-last equality, we observed that $\langle N, d_Y N \rangle = Y(\langle N, N \rangle) = 0$. By taking the trace of the last expression with respect to the first fundamental form $\Is$, we also obtain the identity:
    \[
    \mathrm{tr}_{\Is}(\dot{\II}) = - f \mathrm{tr}(B^2) + \Delta_{\Is} f - 2 f .
    \]
    
    To compute the infinitesimal variation of the mean curvature $\dot H_\sigma$, we observe that (formally, we treat $\Is_t$ and $\II_t$ as matrices here, with respect to a local frame, since the result does not depend on the frame)

     \[
    \dot H_\sigma =\left.\frac{\mathrm{d}}{\mathrm{d}t}\right|_{t=0} \mathrm{tr}_{\Is_t}\II_t= \mathrm{tr}_{\Is} \dot \II + \mathrm{tr}\left(\dot{(\Is^{-1})}\II\right). 
    \]
    
  By the above computation, $\dot \Is=-2f\II$, and since $\dot{(\Is^{-1})}=-(\Is^{-1})\dot\Is (\Is^{-1})$, we obtain
$$\mathrm{tr}\left(\dot{(\Is^{-1})}\II\right)=2f\tr\left((\Is^{-1}\II)\right)=2f\tr(B^2)~.$$
   In conclusion, $ \dot H_\sigma= f \cdot \mathrm{tr}(B^2) -2f +\Delta_{\Is} f$ as claimed.
     \end{proof}
    

 We conclude by applying Lemma~\ref{lemma: formula variation H} to a positive minimal immersion, rewriting the formula in terms of the quantities of our interest as in Theorem~\ref{thm: II sum of quadratic diffs} and Corollary~\ref{cor: GC semplificato}.

\begin{cor}
   Let $\Sigma$ be an oriented surface, let $\sigma:\Sigma\to \PSL(2,\C)$ be a positive minimal immersion in $\PSL(2,\C)$, such that the first fundamental form is of the form $\Is:=\varrho h$ for $h$ a positive complex metric of constant curvature $-1$ and $\varrho\colon \Sigma\to \C^*$ a smooth function, and the second fundamental form is of the form $\II=(q_1+\overline{q_2})/2$ for $q_1\in \QD(\Sigma,J_1)$ and $\overline{q_2}\in \QD(\Sigma,-J_2)$.
   
   Given $f\in C^\infty(\Sigma,\C)$, if we set $$\sigma_t(p):=\exp_{\sigma(p)}(tf(p)N(p))~.$$
    then the derivative of the mean curvature of $\sigma$ satisfies
  
    $$\left. \frac{\mathrm{d}}{\mathrm{d}t} \right|_{t=0} H_\sigma =\varrho^{-1}\Delta_{h}f +  \frac{2q_1\overline {q_2}}{\varrho^2h^2}f -2f $$

\end{cor}
\begin{proof}
    It suffices to apply Lemma~\ref{lemma: formula variation H}, taking into account that $\Delta_{\Is}=\varrho^{-1}\Delta_h$ and, by the Cayley-Hamilton theorem $B^2=-\det(B)\mathrm{id}$ so that $\tr(B^2)=-2\det(B)=-2\det_{\Is}\II=2(q_1\overline{q_2})/(\varrho^2 h^2)$, where the last identity has been verified in the proof of Corollary \ref{cor: GC semplificato}.
\end{proof}

In particular, when $\Sigma=\widetilde S$ is the universal cover of a closed surface $S$, all quantities of interest descend to $S$, and choosing $h=\hpair$ and $\varrho=e^{2u}$ we obtain the following formula on $S$:

\begin{equation}\label{eq:linearized H on S}
    \left. \frac{\mathrm{d}}{\mathrm{d}t} \right|_{t=0} H_\sigma = e^{-2u}\Delta_{\hpair}f +  \frac{2q_1\overline {q_2}}{e^{4u}(\hpair)^2}f -2f = e^{-2u} \partial_u \mathscr G(f)
\end{equation}
where we have used the expression of $\partial_u \mathscr G$ from Equation~\eqref{eq: linearized Gauss}.

\printbibliography

\end{document}